\documentclass{article}

\PassOptionsToPackage{numbers, compress}{natbib}

\usepackage[preprint]{neurips_2022}

\usepackage[utf8]{inputenc} 
\usepackage[T1]{fontenc}    
\usepackage{hyperref}       
\usepackage{url}            
\usepackage{booktabs}       
\usepackage{amsfonts}       
\usepackage{nicefrac}       
\usepackage{microtype}      
\usepackage{xcolor}         

\usepackage{caption}
\usepackage{graphicx}
\usepackage{float}
\usepackage{subfigure}
\usepackage{babel,blindtext}
\usepackage[ruled,vlined,linesnumbered]{algorithm2e}
\usepackage{amsthm}
\usepackage{amssymb}
\usepackage{amsmath} 
\usepackage[short]{optidef}

\usepackage{algorithmic}
\usepackage{array}
\usepackage{lipsum}
\usepackage{hyperref}
\usepackage{CJK}
\usepackage{verbatim}
\usepackage{enumitem}

\newtheorem{theorem}{Theorem}[section]
\newtheorem{lemma}[theorem]{Lemma}
\newtheorem{assumption}[theorem]{Assumption}

\DeclareMathOperator*{\argmin}{arg\,min}
\DeclareMathOperator{\Tr}{Tr}

\title{Decentralized Bilevel Optimization}

\author{%
  Xuxing Chen\thanks{Department of Mathematics, University of California, Davis} \\
  \texttt{xuxchen@ucdavis.edu} \\
   \And
   Minhui Huang\thanks{Department of Electrical and Computer Engineering, University of California, Davis} \\
   \texttt{mhhuang@ucdavis.edu} \\
   \And
   Shiqian Ma
   \footnotemark[1]\\
   \texttt{sqma@ucdavis.edu} \\
}

\begin{document}

\maketitle

\begin{abstract}
  Bilevel optimization has been successfully applied to many important machine learning problems. Algorithms for solving bilevel optimization have been studied under various settings. In this paper, we study the nonconvex-strongly-convex bilevel optimization under a decentralized setting. We design decentralized algorithms for both deterministic and stochastic bilevel optimization problems. Moreover, we analyze the convergence rates of the proposed algorithms in difference scenarios including the case where data heterogeneity is observed across agents. Numerical experiments on both synthetic and real data demonstrate that the proposed methods are efficient. 
\end{abstract}

    \section{Introduction}
Bilevel optimization provides a framework for solving problems arising from meta learning \cite{snell2017prototypical, bertinetto2018meta, rajeswaran2019meta}, hyperparameter optimization \cite{pedregosa2016hyperparameter, franceschi2018bilevel}, reinforcement learning\cite{hong2020two}, etc. It aims at minimizing an objective in the upper level under a constraint given by another optimization problem in the lower level, and has been studied intensively in recent years \cite{franceschi2018bilevel, bertinetto2018meta, ghadimi2018approximation,ji2021bilevel,hong2020two,chen2021closing}. Mathematically, it can be formulated as:
\begin{mini}
    {x\in\mathbb{R}^{p}}{\Phi(x) = f(x,y^*(x)),\quad \text{(upper level)}}
    {\label{eq: bilevel_opt}}{}
    \addConstraint{y^*(x) = \argmin_{y\in\mathbb{R}^{q}}g(x,y),\quad \text{(lower level)}}
\end{mini}
where $g$ is the lower level function which is usually assumed to be strongly convex with respect to $y$ for all $x$, and $f$ is the upper level function which is possibly non-convex. Designing a bilevel optimization algorithm typically consists of two parts: the lower loop and the upper loop. In the lower loop one can run algorithms like gradient descent on $y$ to find a good estimation of the global minimum of $g$ given $x$, which is guaranteed by the strong convexity of $g$. In the upper loop one can run gradient descent on $x$, which requires the estimation of the hypergradient $\nabla \Phi(x)$.

Decentralized optimization aims at solving the finite-sum problem:
\begin{mini}
    {x\in\mathbb{R}^{p}}{\frac{1}{n}\sum_{i=1}^{n}f_i(x),}
    {}{}
\end{mini}
where the $i$-th agent only has access to the information related to $f_i$, and each agent communicates with neighbors so that they can cooperate to solve the original problem. There is no central server collecting local updates. Decentralized algorithms are better choices in certain scenarios \cite{lian2017can, tang2018d}. Since decentralized training has been proved to be efficient, it is natural to ask:
\begin{center}
    \textit{Can we design an algorithm to solve bilevel optimization problems in a decentralized regime?}
\end{center}
We will see later that the answer is affirmative. Our contributions can be summarized as follows.
\begin{itemize}[leftmargin=*]
    \item We propose a novel algorithm that estimates the hypergradient in different cases.
    \item We design a decentralized bilevel optimization (DBO) algorithm and analyze its convergence rate. We also analyze the convergence results for the stochastic version of DBO. To the best of our knowledge, our paper is the first work proposing provably convergent decentralized bilevel optimization algorithms in the presence of data heterogeneity.
    \item We study the effect of gradient tracking in the deterministic decentralized bilevel optimization and analyze the convergence rates.
    \item We conduct numerical experiments on several hyperparameter optimization problems. The results demonstrate the efficiency of our algorithms.
\end{itemize}

\subsection{Related work}
\textbf{Bilevel optimization} can be dated back to \cite{bracken1973mathematical}. Due to its great success in solving problems in meta learning \cite{snell2017prototypical, bertinetto2018meta, rajeswaran2019meta} and hyperparameter optimization \cite{pedregosa2016hyperparameter, franceschi2018bilevel}, there is a flurry of work proposing and analyzing bilevel optimization algorithms. The major challenge in bilevel optimization is the estimation of the hypergradient $\nabla \Phi(x)$ in \eqref{eq: bilevel_opt}, which includes Jacobian and Hessian matrices in the closed-form expression. There are several strategies to overcome this: approximate implicit differentiation (AID) \cite{domke2012generic, pedregosa2016hyperparameter, gould2016differentiating, ghadimi2018approximation, grazzi2020iteration, ji2021bilevel}, iterative differentiation (ITD) \cite{domke2012generic, maclaurin2015gradient, franceschi2018bilevel, grazzi2020iteration, ji2021bilevel} and Neumann series-based approach \cite{ghadimi2018approximation, chen2021closing,hong2020two,ji2021bilevel}. Provably convergent algorithms (both deterministic and stochastic cases) include BSA \cite{ghadimi2018approximation}, TTSA \cite{hong2020two}, stocBiO \cite{ji2021bilevel}, ALSET \cite{chen2021closing}, to name a few.

\textbf{Decentralized optimization} plays a key role in distributed optimization. Under a decentralized setting, the data is distributed to different agents, and each agent communicates with neighbors to solve a finite-sum minimization problem. As opposed to centralized optimization, decentralized optimization aims at solving the problem without a central server that collects iterates from local agents. The main challenge is the data heterogeneity across agents, which should be mitigated by communications. It has been proved that decentralized algorithms have their own advantages \cite{lian2017can} under certain circumstances, for example, low network bandwidth will greatly hinder the communication with the central server if the algorithm is designed to be centralized. An important approach to accelerate the decentralized algorithms is gradient tracking, which has been proved to be efficient \cite{xu2015augmented, di2016next, qu2017harnessing, nedic2017achieving}. We refer the interested readers to \cite{xin2020decentralized}, which provides a comprehensive review of decentralized optimization in a unified variance reduction framework.

There exist recent works considering bilevel optimization under distributed setting. Bilevel optimization under a federated setting has received some attention recently \cite{tarzanagh2022fednest, fedbilevel}, and so does min-max optimization under various distributed settings \cite{xian2021faster, luo2022decentralized, sharma2022federated}. However, none of these papers considers bilevel optimization under the decentralized setting. We notice that there is a concurrent work \cite{lu2022decentralized} also studying decentralized bilevel optimization. However, it aims at solving decentralized bilevel optimization problems under a personalized learning setting, in a sense that the lower level problems are different among agents. In Section \ref{sec: algorithms} we will see that our problem is substantially different from \cite{lu2022decentralized}. 
    \section{Preliminaries}
In this paper we consider the following decentralized optimization problem:
\begin{mini}
    {x\in\mathbb{R}^{p}}{\Phi(x) = \frac{1}{n}\sum_{i=1}^{n}f_i(x,y^*(x)),\quad \text{(upper level)}}
    {}{}
    \addConstraint{y^*(x) = \argmin_{y\in\mathbb{R}^{q}}g(x,y):= \frac{1}{n}\sum_{i=1}^{n}g_i(x,y),\quad \text{(lower level)}}
\end{mini}
where $x\in \mathbb{R}^{p}, y\in\mathbb{R}^{q}$. $f_i$ is possibly nonconvex and $g_i$ is strongly convex in $y$. Here $n$ denotes the number of agents, and agent $i$ only has access to information related to the local objective $f_i$:
\[
    f_i(x,y) = \mathbb{E}_{\phi\sim \mathcal{D}_{f_i}}\left[F(x,y;\phi)\right],\quad g_i(x,y) = \mathbb{E}_{\xi\sim \mathcal{D}_{g_i}}\left[G(x,y;\xi)\right].
\]
In practice we can replace the expectation by empirical loss,
\[
    f_i(x,y) = \frac{1}{n_{f_i}}\sum_{j=1}^{n_{f_i}}F(x,y;\phi_{ij}),\quad g_i(x,y) = \frac{1}{n_{g_i}}\sum_{j=1}^{n_{g_i}}G(x,y;\xi_{ij}),
\]
and then use mini-batch or full batch gradient descent in the updates. When we use mini-batch gradient descent, we call it {\bf "stochastic case"}, and when we use full batch gradient descent, we call it {\bf "deterministic case"}. 
We will study the convergence rates under these two cases in Section \ref{sec: algorithms}.

{\bf Notation.} We denote by $\nabla f(x,y)$ and $\nabla^2f(x,y)$ the gradient and Hessian matrix of $f$, respectively. We use $\nabla_x f(x,y)$ and $\nabla_y f(x,y)$ to represent the gradients of $f$ with respect to $x$ and $y$, respectively. Denote by $\nabla_{xy}f(x,y)\in\mathbb{R}^{p\times q}$ the Jacobian matrix of $f$ and $\nabla_y^2f(x,y)$ the Hessian matrix of $f$ with respect to $y$. $\|\cdot\|$ denotes the $\ell_2$ norm for vectors and spectral norm for matrices, unless specified. $\mathbf{1}_n$ is the all one vector in $\mathbb{R}^{n}$, and $J_n = \mathbf{1}_n\mathbf{1}_n^{\mathsf{T}}$ is the $n\times n$ all one matrix.

The following assumptions will be used, which are standard in bilevel optimization \cite{chen2021closing,ghadimi2018approximation,hong2020two,ji2021bilevel} and decentralized optimization literature \cite{qu2017harnessing, nedic2017achieving, lian2017can, tang2018d}.

\begin{assumption}\label{assump: lip}
    For any $i$, functions $f_i$, $\nabla f_i$, $\nabla g_i$, $\nabla^2g_i$ are $L_{f,0}, L_{f,1}, L_{g,1}, L_{g,2}$ Lipschitz continuous respectively. Moreover, we define $L = \max\{L_{f,1}, L_{g,1}\}$, and $\kappa = \frac{L}{\mu}$.
\end{assumption}

\begin{assumption}\label{assump: strong_convexity}
    Function $g_i$ is $\mu$-strongly convex in $y$ for all $i$.
\end{assumption}


\begin{assumption}\label{assump: W}
    The weight matrix $W = (w_{ij})\in \mathbb{R}^{n\times n}$ is symmetric and doubly stochastic, i.e.:
    \[
        W = W^{\mathsf{T}},\quad W\mathbf{1_n} = \mathbf{1_n},\quad w_{ij}\geq 0,\forall i,j,
    \]
    and its eigenvalues satisfy $1 = \lambda_1> \lambda_2\geq ...\geq \lambda_n$ and $\rho := \max\{|\lambda_2|, |\lambda_n|\} < 1$.
\end{assumption}

\begin{assumption}\label{assump: data_similarity}
    (Data homogeneity on $g$) Assume the data associated with $g_i$ is independent and identically distributed, i.e., $\mathcal{D}_{g_i} = \mathcal{D}_{g}$. (Note that we do not require data homogeneity in the upper level function.)
\end{assumption}


For the stochastic algorithm in Section \ref{sec: algorithms}, we assume:
\begin{assumption}\label{assump: stoc_derivatives}
    (Bounded variance) The stochastic derivatives $\nabla f_i(x,y;\phi)$, $\nabla g_i(x,y;\phi)$, $\nabla^2 g_i(x,y;\phi)$ are unbiased with bounded variances $\sigma_f^2$, $\sigma_{g,1}^2$, $\sigma_{g,2}^2$, respectively.
\end{assumption}

    \section{Our Algorithms}\label{sec: algorithms}


Here we first discuss the main challenge when there is data heterogeneity, i.e., when Assumption \ref{assump: data_similarity} does not hold. 
Recall that in the outer loop, each node needs to estimate the hypergradient given by:
\begin{equation}\label{eq: global_phi}
    \nabla \Phi_i(x) = \nabla_x f_i(x,y^*(x)) - \nabla_{xy} g(x,y^*(x))\left[\nabla_y^2g(x,y^*(x))\right]^{-1}\nabla_y f_i(x,y^*(x)).
\end{equation}
Note that node $i$ only has access to $\nabla f_i$ and $\nabla g_i$ but it does not have access to $\nabla_{xy}g(x,y^*(x))\left[\nabla_y^2g(x,y^*(x))\right]^{-1}$ which requires the global information about $g$. One natural idea is to use \eqref{local_hyper} as a surrogate (here $y_i^*(x):=\argmin_y g_i(x,y)$):
\begin{equation}\label{eq: local_phi}
    \nabla f_i(x, y_i^*(x)) = \nabla_x f_i(x,y_i^*(x)) - \nabla_{xy} g_i(x,y_i^*(x))\left[\nabla_y^2g_i(x,y_i^*(x))\right]^{-1}\nabla_y f_i(x,y_i^*(x)).
\end{equation}
Unfortunately, this leads to the estimation for the error bound of
\[
    \left\|\left[\nabla^2g_i(x,y)\right]^{-1} - \left[\nabla^2g(x,y)\right]^{-1}\right\|,
\]
which is not a diminishing term in the theoretical analysis without assuming the data homogeneity (Assumption \ref{assump: data_similarity}). Mathematically, we would like to compute $$\left[\sum_{i=1}^{n}\nabla_{xy}g_i(x,y)\right]\left[\sum_{i=1}^{n}\nabla_y^2g_i(x,y)\right]^{-1}$$ on node $i$ for any given $(x,y)$. In the next section we design a novel oracle to solve this subproblem with heterogeneous data at the price of the computation of Jacobian matrices.

\subsection{Jacobian-Hessian-Inverse Product oracle}\label{subsec: JHIP}
We first introduce the Jacobian-Hessian-Inverse Product (JHIP) oracle, which is essentially a decentralized algorithm. Denote by $H_i\in\mathbb{S}_{++}^{q\times q}$ and $J_i\in \mathbb{R}^{p\times q}$ the Hessian and Jacobian of $g_i$. Every agent would like to find $Z\in\mathbb{R}^{q\times p}$ such that:

\begin{equation}\label{eq: JHI_eq}
    \sum_{i=1}^{n}H_iZ = \sum_{i=1}^{n}J_i^{\mathsf{T}}\quad \text{or equivalently, } Z^{\mathsf{T}} = \left[\sum_{i=1}^{n}J_i\right]\left[\sum_{i=1}^{n}H_i\right]^{-1}.
\end{equation}

Notice that this is exactly the optimality condition of:
\begin{equation}\label{eq: min_trace}
    \min_{Z\in\mathbb{R}^{p\times q}}\frac{1}{n}\sum_{i=1}^{n}h_i(Z),\quad \text{where  }h_i(Z)=\frac{1}{2}\Tr{(Z^{\mathsf{T}}H_iZ)} - \Tr{(J_iZ)}.
\end{equation}
The objective in \eqref{eq: min_trace} is strongly convex since each $H_i$ is symmetric positive definite. Hence we can design a decentralized algorithm with gradient tracking so that all the agents can collaborate to solve for \eqref{eq: JHI_eq} without a central server. The algorithm is described in Algorithm \ref{algo: JHI_oracle}, where we use the bold texts to highlight the different updates when the problem is deterministic and stochastic. 
\begin{algorithm}
    \caption{Jacobian-Hessian-Inverse Product oracle}\label{algo: JHI_oracle}
    \SetAlgoLined
    \KwIn{$Z_i^{(0)}\in\mathbb{R}^{q\times p}$, stepsizes $\{\gamma_t\}_{t=0}^{\infty}$, iteration number $N$, and initialization 
    \begin{itemize}
    \item \mbox{ \bf{if deterministic, }} $Y_i^{(0)} = H_iZ_i^{(0)} - J_i^{\mathsf{T}}, \mbox{ and } G_i^{(0)} = H_iZ_i^{(0)},$
    \item \mbox{ \bf{if stochastic, }} $Y_i^{(0)} = \hat{H}_i^{(0)}Z_i^{(0)}-\left[\hat{J}_i^{(0)}\right]^{\mathsf{T}}, \mbox{ and } G_i^{(0)} = \hat{H}_iZ_i^{(0)} - \left[\hat{J}_i^{(0)}\right]^{\mathsf{T}}.$
    \end{itemize}
    }
    
    \KwData{$H_i\in\mathbb{S}_{++}^{q\times q}, J_i\in\mathbb{R}^{p\times q}$ accessible only to agent $i$ (deterministic).\\
    \quad\quad$(\hat{H}_i^{(t)}, \hat{J}_i^{(t)}), t\in\{0,1,...,N-1\}$ accessible only to agent $i$ (stochastic).
    }
    \For{$t=0,1,...,N-1$}{
    \For{$i=1,2,...,n$}{
        $Z_i^{(t+1)} = \sum_{j=1}^{n}w_{ij}Z_j^{(t)} - \gamma_t Y_i^{(t)}$, \\
        $G_i^{(t+1)} = {H_iZ_i^{(t+1)}}$ {\bf if deterministic}, or  ${=\hat{H}_i^{(t+1)}Z_i^{(t+1)}-\left[\hat{J}_i^{(t+1)}\right]^{\mathsf{T}}}$ {\bf if stochastic}, \\
        $Y_i^{(t+1)} =\sum_{i=1}^{n}w_{ij}Y_j^{(t)} + G_i^{(t+1)} - G_i^{(t)}$. \\
    }
    }
    \KwOut{$Z_i^{(N)} $ on each node.}
\end{algorithm}

Note that for the deterministic case we can just maintain $G_{i}^{(t+1)} = H_iZ_i^{(t+1)}$ instead of the gradient $H_iZ_i^{(t+1)} - J_i^{\mathsf{T}}$ because we only use $G_i^{(t+1)}$ in line 5 -- the gradient tracking step, and the constant term $J_i^{\mathsf{T}}$ will be cancelled if we set $G_i^{(t+1)}$ as the gradient. We use $\hat{H}_i^{(t)}Z_i^{(t)}-\left[\hat{J}_i^{(t)}\right]^{\mathsf{T}}$ to represent the stochastic gradient of $h_i(Z)$ at $Z_i^{(t+1)}$. Although this oracle does not require Hessian computation for $H_i$, it still requires Jacobian oracle for $J_i$, which is more expensive than Jacobian-vector product oracle. Moreover, we have the convergence rates that have been well understood \cite{qu2017harnessing, pu2021distributed, xin2020decentralized}.
In general, one can also design other oracles (e.g., decentralized ADMM \cite{mota2013d, chang2014multi, shi2014linear, aybat2017distributed, makhdoumi2017convergence}) to solve \eqref{eq: min_trace}. The convergence rates of this algorithm are summarized in Lemma \ref{lem: JHI_error} in supplementary materials.

\subsection{Hypergradient estimate}\label{subsec: HE}

\begin{algorithm}
    \caption{Hypergradient estimate}\label{algo: Hypergrad}
    \SetAlgoLined
    \KwIn{$\nabla f_i(x,y), \nabla_{xy} g_i(x,y), \nabla_y^2g_i(x,y), N$.}
    \KwData{Samples $\phi = (\phi_i^{0},...,\phi_i^{(N+1)})$ on node $i$.}
    \eIf{Assumption \ref{assump: data_similarity} holds}{
        \If{Deterministic case}{
            Run $N$-step conjugate gradient method on $\left[\nabla_y^2g_i(x,y) \right]v = \nabla_y f_i(x,y)$ to get $v^{N}$.\\
            Set $\hat{\nabla}f_i = \nabla_xf_i(x,y) - \nabla_{xy}g_i(x,y)v^{N}$\quad \text{ -- AID based approach}.
        }
        \If{Stochastic case}{
            Set $\hat{\nabla} f_i = \nabla_xf_i(x, y;\phi^{(0)}) - \nabla_{xy}g_i(x, y;\phi^{(1)})\cdot H\cdot\nabla_yf_i(x, y;\phi^{(0)})$, where $H = \left[\epsilon M\prod_{n=1}^{M'}(I - \epsilon\nabla_y^2g_i(x, y;\phi^{(n+1)}))\right]$, and $M'$ is drawn from $\{1,2,...,M\}$ uniformly at random\quad \text{ -- Neumann series based approach}.
        }
    }
    {
    \If{Deterministic case}{
        Run $N$-step deterministic Algorithm \ref{algo: JHI_oracle} with $H_i = \nabla_y^2g_i(x,y), J_i = \nabla_{xy}g_i(x,y)$ and stepsize $\gamma_t =\gamma$ to get $Z_i^{(N)}$. \\
        Set $\hat{\nabla}f_i = \nabla_x f_i(x,y) -\left[Z_i^{(N)}\right]^{\mathsf{T}}\nabla_y f_i(x,y)$\quad \text{ -- deterministic JHIP}.
    }
    \If{Stochastic case}{
        Run $N$-step stochastic Algorithm \ref{algo: JHI_oracle} with $\hat{H}_i^{(t)} = \nabla_y^2g_i(x,y;\phi_i^{(t)}), \hat{J}_i^{(t)} = \nabla_{xy}g_i(x,y;\phi_i^{(t)})$ and $\gamma_t = \mathcal{O}(\frac{1}{t})$ to get $Z_i^{(N)}$. \\
        Set $\hat{\nabla}f_i = \nabla_x f_i(x,y;\phi_i^{(0)}) -\left[Z_i^{(N)}\right]^{\mathsf{T}}\nabla_y f_i(x,y;\phi_i^{(0)})$\quad \text{ -- stochastic JHIP}.
    }}
    \KwOut{$\hat{\nabla}f_i$ on node $i$.}
\end{algorithm}
Before we propose our algorithms, we first introduce hypergradient estimates (Algorithm \ref{algo: Hypergrad}) under different cases.
\begin{itemize}[leftmargin=*]
    \item {\bf Case 1: Deterministic + homogeneous data}. In this case the hypergradient is estimated based on the AID approach \cite{ji2021bilevel}, which is essentially utilizing conjugate gradient method, and only requires Hessian-vector product oracles instead of explicit Hessian matrix computation. We adopt the approximation error (Lemma 3 in \cite{ji2021bilevel}) in Lemma \ref{lem: aid_approx}.
    \item {\bf Case 2: Stochastic + homogeneous data}. In this case the hypergradient is estimated based on the Neumann series approach \cite{ghadimi2018approximation}. Gradients $\nabla f_i$ and $\nabla g_i$ are replaced by their corresponding first order stochastic oracles (i.e., stochastic gradients). We have the error estimation in Lemma \ref{lem: neumann_error}.
    \item {\bf Case 3: Heterogeneous data}. In this case we compute the global Jacobian-Hessian-Inverse product by using the JHIP oracle (Algorithm \ref{algo: JHI_oracle}). The error estimation results are given in Lemma \ref{lem: inner_error}.
\end{itemize}

\subsection{Deterministic decentralized bilevel optimization}\label{subsec: DBO}

\begin{algorithm}
    \caption{(Deterministic) Decentralized Bilevel Optimization}\label{algo:DBO}
	\SetAlgoLined
	\KwIn{$W, N, K, T, \eta_x, \eta_y, x_{i,0}, y_i^{(0)}$.}
	\For{$k=0,1,...,K-1$}{
	    $y_{i,k}^{(0)} = y_{i,k-1}^{(T)} \text{ if } k>0 \text{ otherwise }  y_{i,k}^{(0)} = y_i^{(0)}$.\\
		\For{$t=0,1,..., T-1$}{
		    \For{$i=1,2,...,n$}{
		        \eIf{Assumption \ref{assump: data_similarity} holds}{
		            $y_{i,k}^{(t+1)} = y_{i,k}^{(t)} -\eta_y \nabla_y g_i(x_{i,k},y_{i,k}^{(t)})$.
		        }{
		            $v_{i,k}^{(t)} = \sum_{j=1}^{n}w_{ij}v_{j,k}^{(t-1)} + \nabla_y g_i(x_{i,k}, y_{i,k}^{(t)}) - \nabla_y g_i(x_{i,k}, y_{i,k}^{(t-1)})$. \\
		             $y_{i,k}^{(t+1)} = \sum_{j=1}^{n}w_{ij}y_{j,k}^{(t)} - \eta_y v_{i,k}^{(t)}$.
		        }
		    }
        }
        \For{$i=1,2,...,n$}{
            Run Algorithm \ref{algo: Hypergrad} (with "deterministic case" option) and set the output as $\hat{\nabla}f_i(x_{i,k},y_{i,k}^{(T)})$. \\
    		$x_{i,k+1} = \sum_{j=1}^{n}w_{ij}x_{j,k} -\eta_x \hat{\nabla} f_i(x_{i,k},y_{i,k}^{(T)})$.
        }
	}
\end{algorithm}

We propose the decentralized bilevel optimization algorithm (DBO) in Algorithm \ref{algo:DBO}. 
In the inner loop (lines 3-9) each agent performs local gradient descent updates for variable $y$ in parallel (line 4-8). When Assumption \ref{assump: data_similarity} holds, we can simply run local gradient descent without communication because in the lower level local distribution already captures the global function information. When Assumption \ref{assump: data_similarity} does not hold, then the data distribution is substantially heterogeneous across agents, so we add weighted averaging steps (line 9) to reach consensus and gradient tracking steps (line 8) to reduce the complexity. In the outer loop (lines 11-14) each agent communicates with neighbors and then performs gradient descent for variable $x$. For simplicity we use constant stepsize $\eta_x$ in the outer loop. Similar results can be obtained for diminishing stepsizes. Using proper parameters, we have the following sublinear convergence results of Algorithm \ref{algo:DBO} for solving the deterministic problem.

\begin{theorem}\label{thm:dbo}
    In Algorithm \ref{algo:DBO}, suppose Assumptions \ref{assump: lip}, \ref{assump: strong_convexity}, and \ref{assump: W} hold. If Assumption \ref{assump: data_similarity} holds, we set $\eta_x = \Theta(K^{-\frac{1}{3}}n^{-\frac{1}{3}}\kappa^{-\frac{8}{3}}),\ \eta_y\leq \frac{2}{\mu + L},\ T=\Theta(\log(\kappa)),\ N=\Theta(\log(\kappa))$. If Assumption \ref{assump: data_similarity} does not hold, we set $\eta_x = \Theta(K^{-\frac{1}{3}}n^{-\frac{1}{3}}\kappa^{-\frac{8}{3}}),\ \eta_y^{(t)} = \mathcal{O}(1),\ T = \Theta(\log K),\ N=\Theta(\log K)$. In both cases, we have:
    $
        \frac{1}{K+1}\sum_{j=0}^{K}\|\nabla\Phi(\overline{x_j})\|^2 = \mathcal{O}\left(\frac{n^{\frac{1}{3}}\kappa^{\frac{8}{3}}}{K^{\frac{2}{3}}}\right),
    $ where $\overline{x_j} = \frac{1}{n}\sum_{i=1}^{n}x_{i,j}$.
\end{theorem}

This convergence rate matches the previous results in nonconvex decentralized optimization \cite{lian2017can,tang2018d}. We conjecture that the bound $N=\Theta(\log K)$ in the second case is not tight, and can be improved by using a technique similar to the case when Assumption \ref{assump: data_similarity} holds. We leave this as a future work.

\subsection{Deterministic decentralized bilevel optimization with gradient tracking}\label{subsec: DBOGT}

In this section we study the effect of gradient tracking in decentralized bilevel optimization. We propose the Decentralized Bilevel Optimization with Gradient Tracking (DBOGT) Algorithm \ref{algo:DBOGT}. We introduce $u$ and $v$ to serve as the update direction. Gradient tracking technique has been widely used in distributed optimization literature \cite{di2016next, qu2017harnessing, nedic2017achieving}. The effect of adding such a direction accelerates the algorithm in a sense that one can choose a constant stepsize that is independent of the iteration number. For Algorithm \ref{algo:DBOGT} we have the following theorem.

\begin{algorithm}[H]
    \caption{(Deterministic) Decentralized Bilevel Optimization with Gradient Tracking}\label{algo:DBOGT}
	\SetAlgoLined
	\KwIn{$W, N, K, T, \eta_x, \eta_y, x_i^{(0)}, y_i^{(0)}$.}
	\For{$k=0,1,...,K-1$}{
	    $y_{i,k}^{(0)} = y_{i,k-1}^{(T)} \text{ if } k>0 \text{ otherwise }  y_{i,k}^{(0)} = y_i^{(0)}$. \\
		\For{$t=0,1,..., T-1$}{
		\For{$i=1,2,...,n$}{
		    \eIf{Assumption \ref{assump: data_similarity} holds}{
		        $y_{i,k}^{(t+1)} = y_{i,k}^{(t)} -\eta_y \nabla_yg_i(x_{i,k},y_{i,k}^{(t)})$.
		    }{
		        $v_{i,k}^{(t)} = \sum_{j=1}^{n}w_{ij}v_{j,k}^{(t-1)} + \nabla_y g_i(x_{i,k}, y_{i,k}^{(t)}) - \nabla_y g_i(x_{i,k}, y_{i,k}^{(t-1)})$. \\
		        $y_{i,k}^{(t+1)} = \sum_{j=1}^{n}w_{ij}y_{j,k}^{(t)} - \eta_y v_{i,k}^{(t)}$.
		    }
		}
        }
        \For{$i=1,2,...,n$}{
        Run Algorithm \ref{algo: Hypergrad} (with "deterministic case" option) and set the output as $\hat{\nabla}f_i(x_{i,k},y_{i,k}^{(T)})$. \\
	    $u_{i,k} = \sum_{j=1}^{n}w_{ij}u_{j,k-1} +  \hat{\nabla} f_i(x_{i,k},y_{i,k}^{(T)}) - \hat{\nabla}f_i(x_{i,k-1},y_{i,k-1}^{(T)})$. \\
	    $x_{i,k+1} = \sum_{j=1}^{n}w_{ij}x_{j,k} -\eta_x u_{i,k}$.
        }
	}
\end{algorithm}
\begin{theorem}\label{thm:dbogt}
    In Algorithm \ref{algo:DBOGT}, suppose Assumptions \ref{assump: lip}, \ref{assump: strong_convexity}, and \ref{assump: W} hold. If Assumption \ref{assump: data_similarity} holds, we set $\eta_x=\Theta(\kappa^{-3}),\ \eta_y\leq \frac{2}{\mu + L},\ T = \Theta(\log(\kappa)),\ N = \Theta(\log(\kappa))$.
     If Assumption \ref{assump: data_similarity} does not hold, we set $\eta_x=\Theta(\kappa^{-3}),\ \eta_y=\Theta(1),\ N = \Theta(\log K),\ T = \Theta(\log K)$.
    In both cases, we have
    $  
        \frac{1}{K+1}\sum_{j=0}^{K}\|\nabla\Phi(\overline{x_j})\|^2 = \mathcal{O}\left(\frac{1}{K}\right),
    $ where $\overline{x_j} = \frac{1}{n}\sum_{i=1}^{n}x_{i,j}$.
\end{theorem}
Note that this result implies that in DBOGT we can set $\eta_x$ as a constant that is independent of the total number of iterations $K$, which matches the results in gradient tracking literature \cite{qu2017harnessing, nedic2017achieving, koloskova2021improved}.

\subsection{Decentralized stochastic bilevel optimization}\label{subsec: DSBO}

Our stochastic version of the DBO algorithm: Decentralized Stochastic Bilevel Optimization (DSBO), is described in Algorithm \ref{algo:DSBO}. Its convergence rate is given in Theorem \ref{thm:dsbo}.
\begin{algorithm}[H]
    \caption{Decentralized Stochastic Bilevel Optimization}\label{algo:DSBO}
	\SetAlgoLined
	\KwIn{$W, N, K, T, \eta_x, \eta_y, x_i^{(0)}, y_i^{(0)}$.}
	\For{$k=0,1,...,K-1$}{
	    $y_{i,k}^{(0)} = y_{i,k-1}^{(T)} \text{ if } k>0 \text{ otherwise }  y_{i,k}^{(0)} = y_i^{(0)}$. \\
		\For{$t=0,1,..., T-1$}{
		\For{$i=1,2,...,n$}{
		    \eIf{Assumption \ref{assump: data_similarity} holds}{
		        $y_{i,k}^{(t+1)} = y_{i,k}^{(t)} -\eta_y \nabla_y g_i(x_{i,k},y_{i,k}^{(t)}; \xi_{i,k}^{(t)})$.
		    }{
                $y_{i,k}^{(t+1)} = \sum_{j=1}^{n}w_{ij}y_{j,k}^{(t)} -\eta_y^{(t)} \nabla_y g_i(x_{i,k},y_{i,k}^{(t)}; \xi_{i,k}^{(t)})$.
		    }
		}
        }
        \For{$i=1,2,...,n$}{
        Run Algorithm \ref{algo: Hypergrad} (with "stochastic case" option) and set the output as $\hat{\nabla}f_i(x_{i,k},y_{i,k}^{(T)};\phi_{i,k})$. \\ 
		$x_{i,k+1} = \sum_{j=1}^{n}w_{ij}x_{j,k} -\eta_x \hat{\nabla} f_i(x_{i,k},y_{i,k}^{(T)};\phi_{i,k})$.
        }
	}
\end{algorithm}
\begin{theorem}\label{thm:dsbo}
    In Algorithm \ref{algo:DSBO}, suppose Assumptions \ref{assump: lip}, \ref{assump: strong_convexity}, and \ref{assump: W} hold. If Assumption \ref{assump: data_similarity} holds, we set $M = \Theta(\log(K) ,\ T = \Theta(\log(\kappa)),\ \epsilon<\frac{1}{L},\ \eta_x = \Theta(K^{-\frac{1}{2}}),\ \eta_y = \Theta(K^{-\frac{1}{2}})$.
    If Assumption \ref{assump: data_similarity} does not hold, we set $\eta_x = \Theta(K^{-\frac{1}{2}}),\ \eta_y = \Theta(K^{-\frac{1}{2}}),\ T=\Theta(K^{\frac{1}{2}}),\ N=\Theta(\log K)$.
    In both cases, we have
    $
        \frac{1}{K+1}\sum_{j=0}^{K}\mathbb{E}\left[\|\nabla\Phi(\overline{x_{j}})\|^2\right] = \mathcal{O}(\frac{1}{\sqrt{K}}),
    $ where $\overline{x_j} = \frac{1}{n}\sum_{i=1}^{n}x_{i,j}$.
\end{theorem}
    \section{Numerical experiments}\label{sec: exp}
In this section we conduct several experiments on hyperparameter optimization problems in the decentralized setting, which can be formulated as:
\begin{mini}
    {\lambda\in\mathbb{R}^{p}}{\Phi(\lambda) = \frac{1}{n}\sum_{i=1}^{n}f_i(\lambda, \tau^*(\lambda)),}
    {\label{HPO}}{}
    \addConstraint{\tau^*(\lambda) = \argmin_{\tau\in\mathbb{R}^q} \frac{1}{n}\sum_{i=1}^{n}g_i(\lambda, \tau).}
\end{mini}
Here $f_i$ and $g_i$ denote the validation loss and training loss on node $i$, respectively. The goal is to find the best hyperparameter $\lambda$ under the constraint that $\tau^{*}(\lambda)$ is the optimal model parameter of the lower level problem. For each experiment, we set our network topology as a special ring network, where $W= (w_{i,j})$ and the only nonzero elements are given by:
\[
    w_{i,i} = a,\ w_{i,i+1}=w_{i,i-1} = \frac{1-a}{2},\ \text{ for some }a\in (0,1).
\]
Here we overload the notation and set $w_{n,n+1}= w_{n,1}, w_{1,0}=w_{1,n}$. Note that $a$ is the unique parameter that determines the weight matrix and will be specified in each experiment.

\subsection{Synthetic data}
\begin{figure*}[ht]
	\centering  
	\subfigure[Logistic regression on synthetic data.]{\label{figure:synthetic_data}\includegraphics[width=0.45\textwidth]{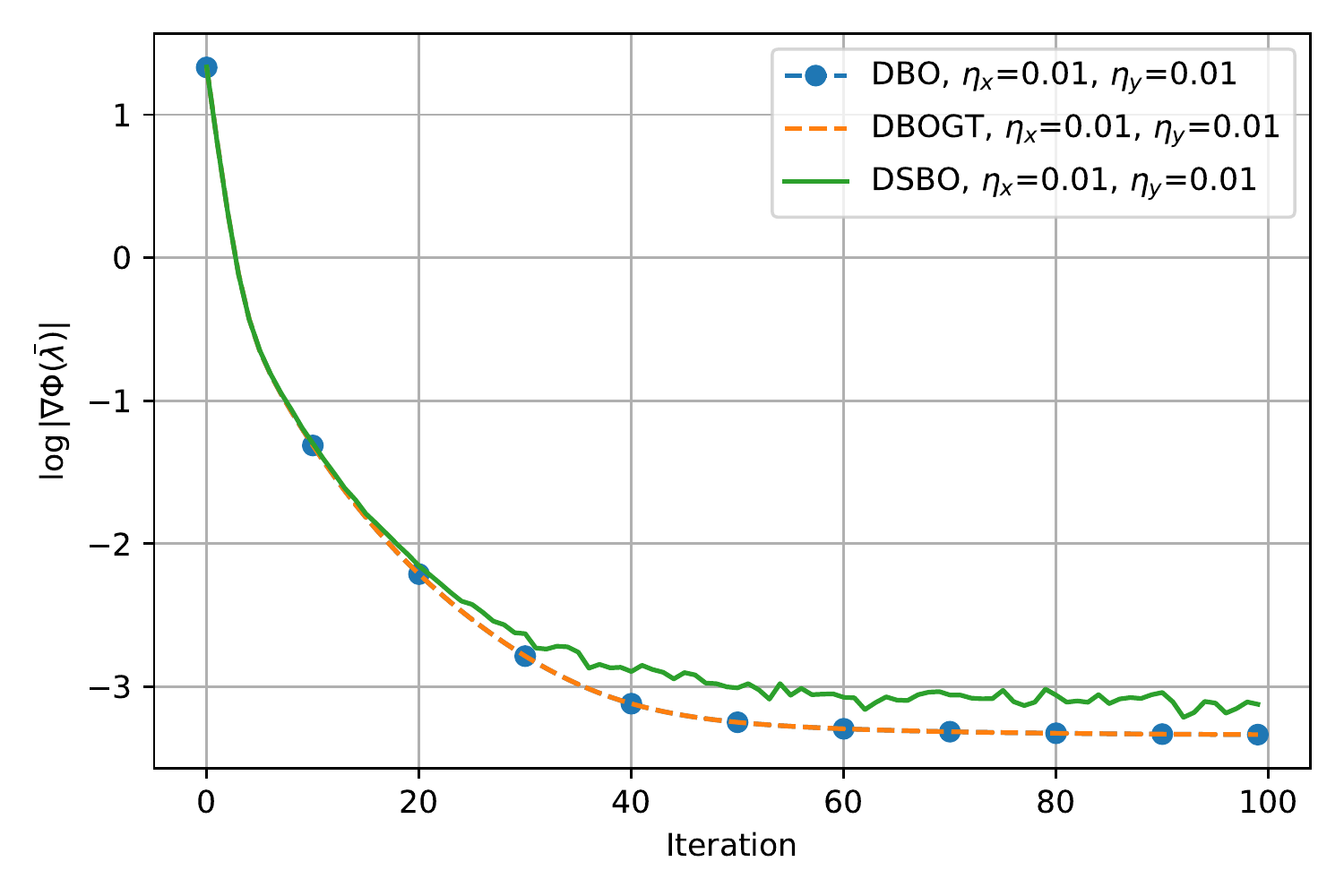} }
	\subfigure[Logistic regression on 20 Newsgroup dataset.]{\label{figure:DSBO}\includegraphics[width=0.45\textwidth]{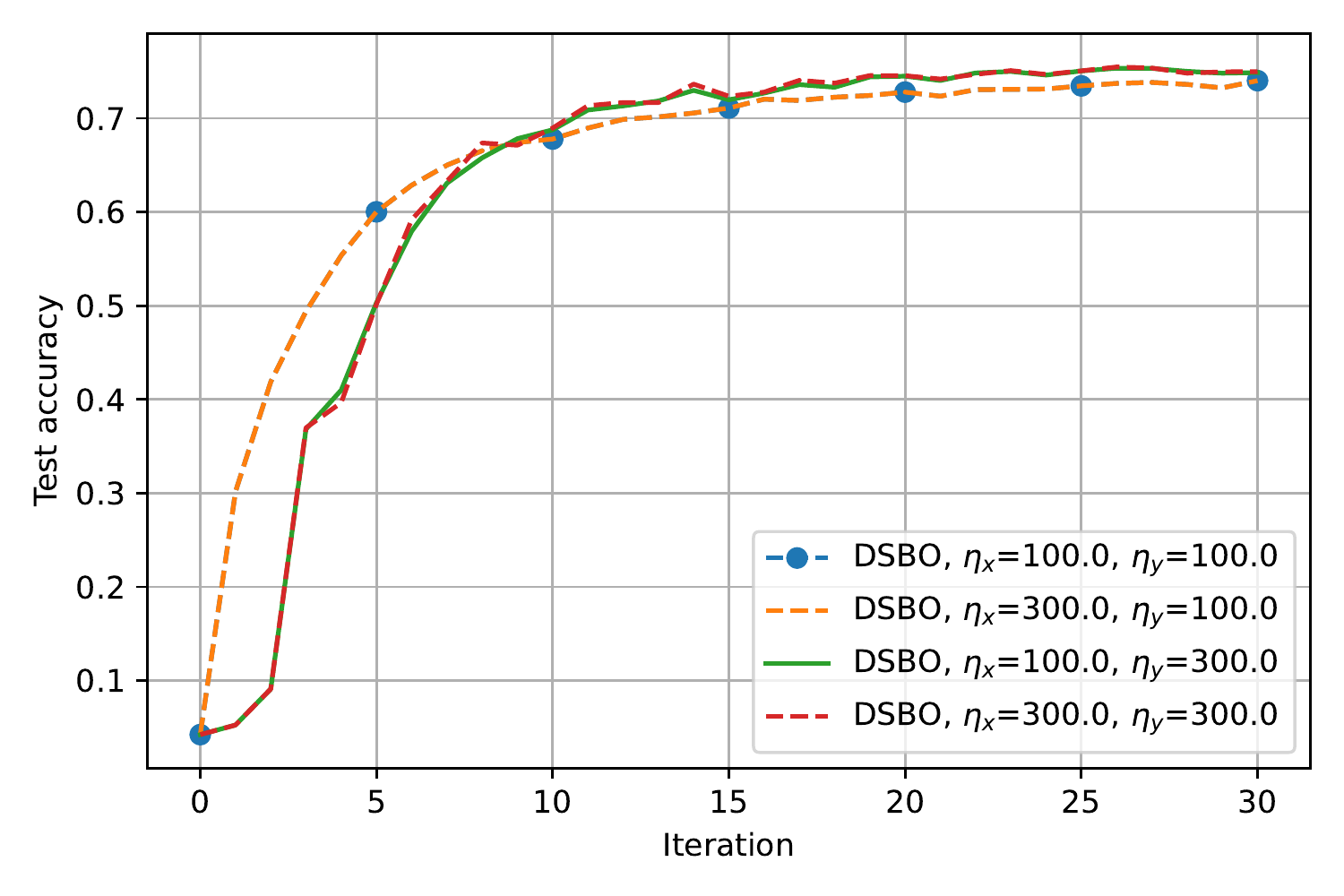} }
	\vspace{-0.2cm}
	\caption{}\label{Figure1}
	\vspace{-0.3cm}
\end{figure*}

We first conduct logistic regression with $l^2$ regularization on synthetic heterogeneous data (e.g., \cite{pedregosa2016hyperparameter}, \cite{grazzi2020iteration}). On node $i$ we have:
\[
    \begin{aligned}
        &f_i(\lambda,\tau^*(\lambda)) = \sum_{(x_e,y_e)\in\mathcal{D}_i'}\psi(y_ex_e^{\mathsf{T}}\tau^*(\lambda)), \\
        &g_i(\lambda,\tau) = \sum_{(x_e,y_e)\in\mathcal{D}_i}\psi(y_ex_e^{\mathsf{T}}\tau) + \frac{1}{2}\tau^{\mathsf{T}}\text{diag}(e^{\lambda})\tau,
    \end{aligned}
\]
where $e^{\lambda}$ is element-wise, $\text{diag}(v)$ denotes the diagonal matrix generated by vector $v$, and $\psi(x) = \log(1+e^{-x})$. $\mathcal{D}_i'$ and $D_i$ represent validation set and training set on node $i$. Following the setup in \cite{grazzi2020iteration}, we first randomly generate $\tau^*\in \mathbb{R}^{p}$ and the noise vector $\epsilon\in \mathbb{R}^{p}$. For the data point $(x_e, y_e)$ on node $i$, each element of $x_e$ is sampled from the normal distribution with mean 0, variance $i^2$. $y_e$ is then set by $y_e = \text{sign}(x_e^{\mathsf{T}}\tau^* + m\epsilon)$, where $\text{sign}$ denotes the sign function and $m=0.1$ denotes the noise rate. In the experiment we choose $p = q = 50,$ and the number of inner-loop and outer-loop iterations as $10$ and $100$ respectively. $N,$ the number of iterations of the JHIP oracle \ref{algo: JHI_oracle} is $20$. The stepsizes are $\eta_x=\eta_y=\gamma = 0.01.$ The number of agents $n$ is chosen as $20,$ and the weight parameter $a=0.4$. We plot the logarithm of the norm of the gradient in Figure \ref{figure:synthetic_data}. From this figure we see that all three algorithms: DBO (Algorithm \ref{algo:DBO}), DBOGT (Algorithm \ref{algo:DBOGT}), and DSBO (Algorithm \ref{algo:DSBO}) can reduce the gradient to an acceptable level. Moreover, DBO and DBOGT have similar performance, and they are both slightly better than DSBO. 

\subsection{Real-world data}
We now conduct the DSBO algorithm on a logistic regression problem on 20 Newsgroup dataset\footnote{\url{http://qwone.com/~jason/20Newsgroups/}}\cite{grazzi2020iteration}. On node $i$ we have:
\[
    \begin{aligned}
        &f_i(\lambda, \tau^*(\lambda))= \frac{1}{|\mathcal{D}^{(i)}_{val}|}\sum_{(x_e,y_e)\in\mathcal{D}^{(i)}_{val}}L(x_e^{\mathsf{T}}\tau^*, y_e), \\
        &g_i(\lambda, \tau)= \frac{1}{|D_{tr}^{(i)}|}\sum_{(x_e,y_e)\in \mathcal{D}^{(i)}_{tr}}L(x_e^{\mathsf{T}}\tau, y_e) + \frac{1}{cp} \sum_{i=1}^{c}\sum_{j=1}^{p}e^{\lambda_j}\tau_{ij}^2,
    \end{aligned}
\]
where $c=20$ denotes the number of topics, $p=101631$ is the feature dimension, $L$ is the cross entropy loss, $\mathcal{D}_{val}$ and $\mathcal{D}_{tr}$ are the validation and training data sets, respectively. Our codes can be seen as decentralized versions of the one provided in \cite{ji2021bilevel}. 
We first set inner and outer stepsizes $\eta_x=\eta_y=100$ (the same as the ones used in \cite{ji2021bilevel}), and then compare its performance with different stepsizes. We set the number of inner-loop iterations $T = 10,$ the number of outer-loop iterations $ K=30,$ the number of agents $ n=20,$ and the weight parameter $a=0.33$. At the end of $j$-th outer-loop iteration we use the average $\overline{\tau_j} = \frac{1}{n}\sum_{i=1}^{n}\tau_{i,j}$ as the model parameter and then do the classification on the test set to get the test accuracy. In Figure \ref{figure:DSBO} we plot the test accuracy of every iteration. From this figure we see that the DSBO algorithm is able to get good test accuracy under different settings of stepsizes.

\begin{figure*}[ht]
	\centering  
	\subfigure[]{\label{figure:DBO_DBOGT1}\includegraphics[width=0.45\textwidth]{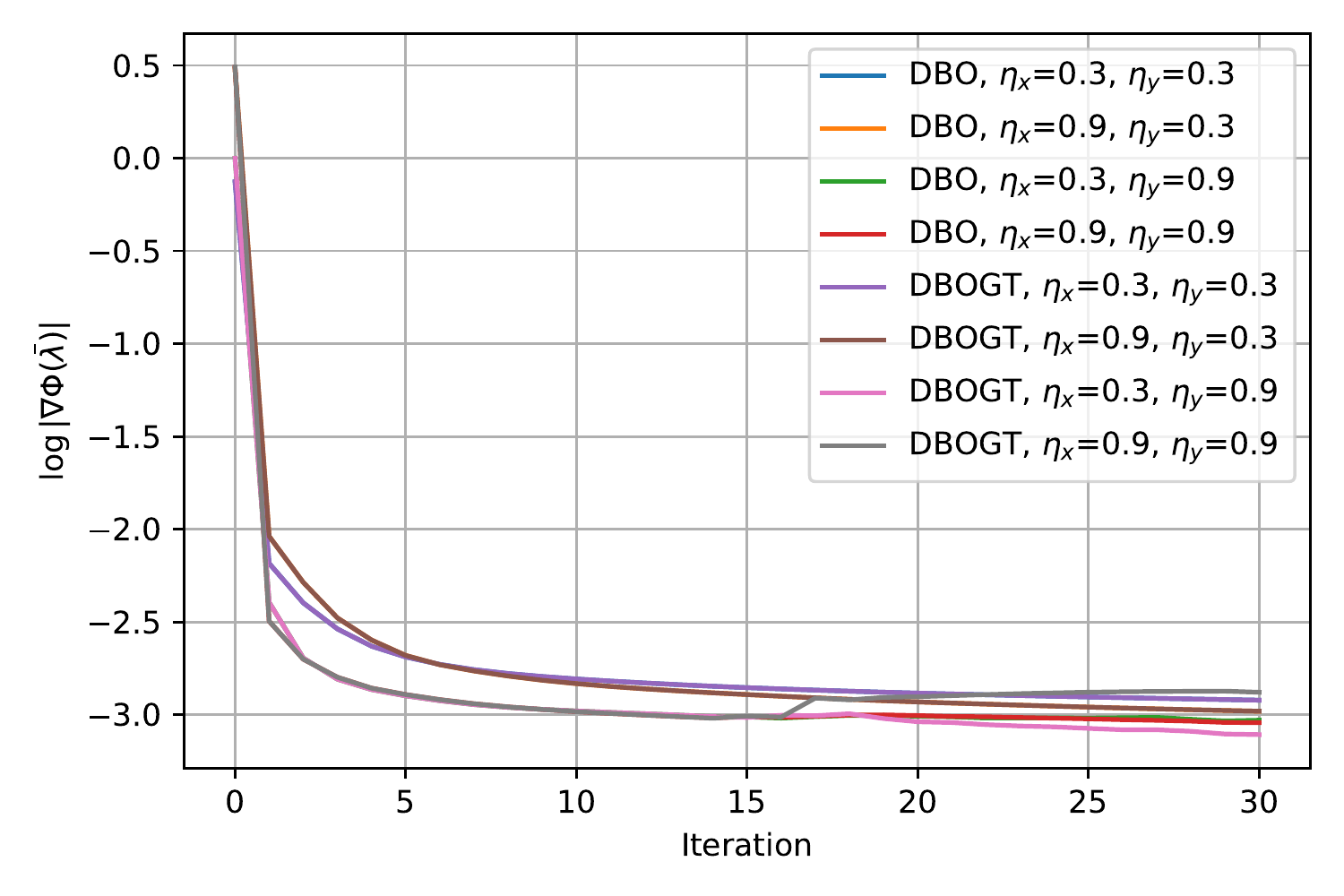} }
	\subfigure[]{\label{figure:DBO_DBOGT2}\includegraphics[width=0.45\textwidth]{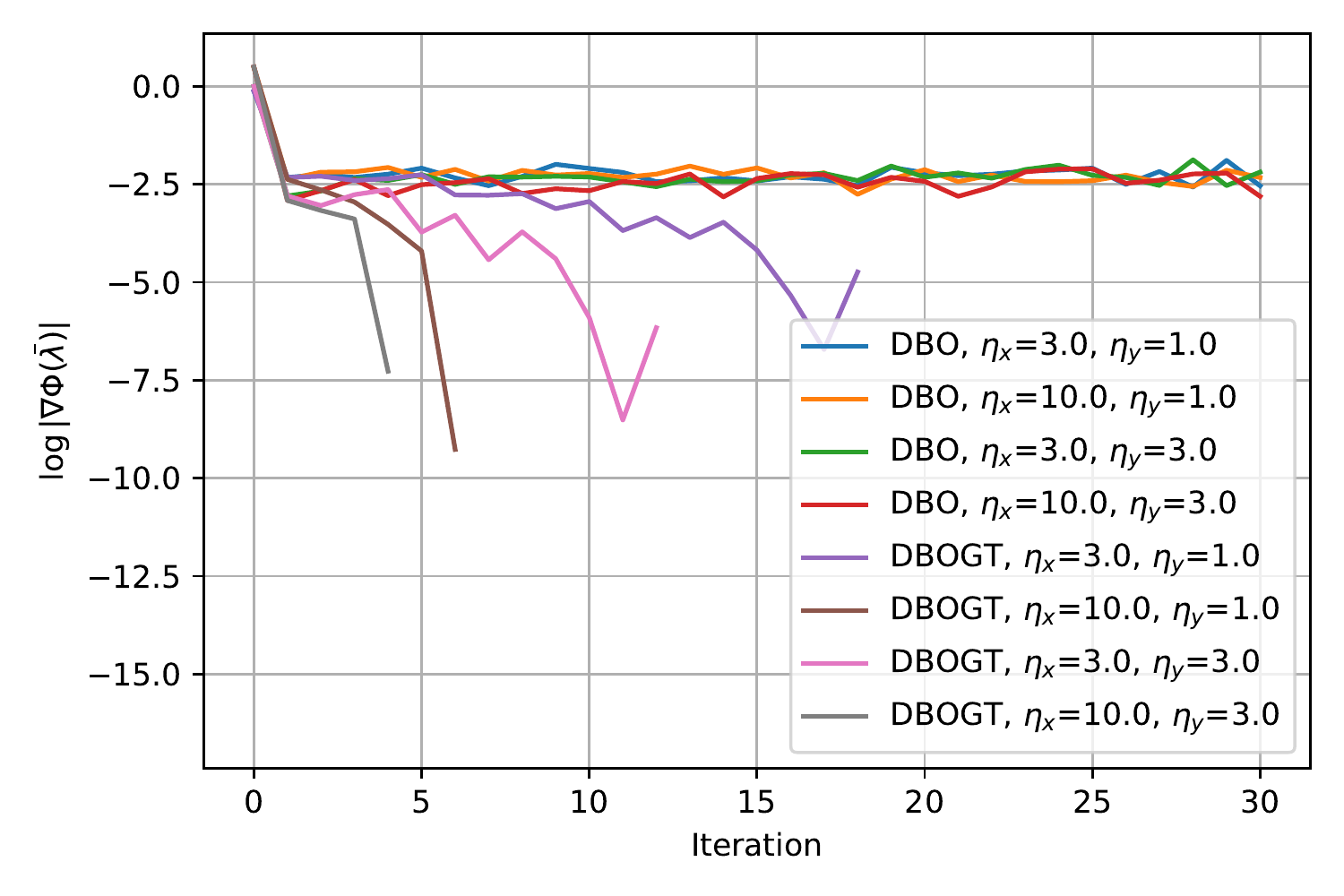} }
	\vspace{-0.2cm}
	\caption{(a) and (b): Data hyper-cleaning on MNIST}\label{Figure2}
	  \vspace{-0.3cm}
\end{figure*}

Finally we apply deterministic DBO and DBOGT algorithms on a data hyper-cleaning problem \cite{shaban2019truncated, ji2021bilevel} for MNIST dataset \cite{lecun1998gradient}. The purpose is to demonstrate the advantage of the gradient tracking technique. On node $i$ we have:
\[
    \begin{aligned}
        f_i(\lambda, \tau)&= \frac{1}{|\mathcal{D}^{(i)}_{val}|}\sum_{(x_e,y_e)\in\mathcal{D}^{(i)}_{val}}L(x_e^{\mathsf{T}}\tau, y_e), \\
        g_i(\lambda, \tau)&= \frac{1}{|D^{(i)}_{tr}|}\sum_{(x_e,y_e)\in \mathcal{D}^{(i)}_{tr}}\sigma(\lambda_e)L(x_e^{\mathsf{T}}\tau, y_e) + C_r\|\tau\|^2,
    \end{aligned}
\]
where $L$ is the cross-entropy loss and $\sigma(x)= (1+e^{-x})^{-1} $ is the sigmoid function. The number of inner-loop iterations $T$ and outer-loop iterations $K$ are set as $10$ and $30$, respectively. The number of agents $n=20$ and the weight parameter $a=0.5$. Following \cite{shaban2019truncated, ji2021bilevel} the regularization parameter $C_r$ is set as $0.001$. We first choose stepsizes similar to those in\cite{ji2021bilevel} and then set larger stepsizes. In each iteration we evaluate the norm of the hypergradient at the average of the hyperparameters $\bar{\lambda}$, and plot the logarithm (base 10) of the norm of the hypergradient versus iteration number in Figure \ref{Figure2}. The Figure \ref{figure:DBO_DBOGT1} shows that the perofromance of DBO and DBOGT are similar when the stepsizes are small. However, the Figure \ref{figure:DBO_DBOGT2} shows that DBOGT converges much faster than DBO when the stepsizes are relatively large. This supports the conclusions in Theorem \ref{thm:dbogt}.


    \section{Conclusion}
In this paper we propose both deterministic and stochastic algorithms for solving decentralized bilevel optimization problems. We obtain sublinear convergence rates when the lower level function is generated by homogeneous data. Moreover, at the price of computing Jacobian matrices, we propose decentralized algorithms with sublinear convergence rates when the lower level function is generated by heterogeneous data. Numerical experiments demonstrate that the proposed algorithms are efficient. It is still an open question whether one can design decentralized optimization algorithms without assuming data homogeneity and deterministic (or stochastic) Jacobian computation. We leave this as a future work.

\newpage
\bibliographystyle{plain}
\bibliography{bibfile.bib}

\appendix

\section{Appendix}
\section*{Convergence analysis}\label{sec: analysis}
In this section we provide the proofs of convergence results. For convenience, we first list the notation below.
\begin{itemize}
    \item $W := (w_{ij})\in \mathbb{R}^{n\times n}$ is symmetric doubly stochastic, and $\rho := \max\{|\lambda_2|, |\lambda_n|\} < 1$.

    \item $X_k := [x_{1,k},x_{2,k},...,x_{n,k}]\in \mathbb{R}^{p\times n}$.
    
    \item $\overline{x_{k}} := \frac{1}{n}\sum_{i=1}^{n}x_{i,k}\in\mathbb{R}^{p}$.
    
    \item $\partial \Phi(X_k) := [\hat{\nabla} f_1(x_{1,k}, y_{1,k}^{(T)}), \hat{\nabla} f_2(x_{2,k}, y_{2,k}^{(T)}),...,\hat{\nabla} f_n(x_{n,k}, y_{n,k}^{(T)}) ] \in \mathbb{R}^{p\times n}$.
    
    \item $\overline{\partial \Phi(X_k)} := \frac{1}{n}\sum_{i=1}^{n}\hat{\nabla}f_{i}(x_{i,k}, y_{i,k}^{(T)})$.
    
    \item $\partial \Phi(X_k;\phi) := [\hat{\nabla} f_1(x_{1,k}, y_{1,k}^{(T)};\phi_{1,k}), \hat{\nabla} f_2(x_{2,k}, y_{2,k}^{(T)};\phi_{2,k}),...,\hat{\nabla} f_n(x_{n,k}, y_{n,k}^{(T)};\phi_{n,k}) ] \in \mathbb{R}^{p\times n}$.
    
    \item $\overline{\partial \Phi(X_k;\phi)} := \frac{1}{n}\sum_{i=1}^{n}\hat{\nabla}f_{i}(x_{i,k}, y_{i,k}^{(T)}; \phi_{i,k})$.
    
    \item $q_{i,k} := x_{i,k} - \overline{x_{k}} = X_k(e_i - \frac{\mathbf{1}_n}{n})$.
    
    \item $Q_k := [q_{1,k}, q_{2,k},...,q_{n,k}]\in\mathbb{R}^{p\times n}$. Assume $x_{i,0} = 0$ and thus $Q_0 = 0$.
    
    \item $r_{i,k} := u_{i,k} - \overline{u_k}$.
    
    \item $R_k := \left[r_{1,k},r_{2,k},...,r_{n,k}\right]\in\mathbb{R}^{p\times n}$.
    
    \item $S_K := \sum_{k=1}^{K}\|Q_k\|^2$.
    
    \item $T_{K} := \sum_{j=0}^{K}\|\nabla\Phi(\overline{x_j})\|^2$.
    
    \item $E_K := \sum_{j=1}^{K}\sum_{i=1}^{n}\|x_{i,j} - x_{i,j-1}\|^2$.
    
    \item When Assumption \ref{assump: data_similarity} holds, define:
    \[
        A_K := \sum_{j=0}^{K}\sum_{i=1}^{n}\|y_{i,j}^{(T)} - y_i^*(x_{i,j})\|^2, B_K := \sum_{j=0}^{K}\sum_{i=1}^{n}\|v_{i,j}^* - v_{i,j}^{(0)}\|^2,
    \]
    where
    \[
    v_{i,j}^* = \left[\nabla_y^2g_i(x_{i,j},y_i^*(x_{i,j}))\right]^{-1}\nabla_y f_i(x_{i,j},y_i^*(x_{i,j})).
    \]
    
    \item $\delta_y^T := (1 - \frac{2\eta_y\mu L}{\mu + L})^T$,  $\delta_{\kappa}^N := (\frac{\sqrt{\kappa} - 1}{\sqrt{\kappa} + 1})^{2N}$.
\end{itemize}

We first introduce a few lemmas that are useful in the proofs. The following lemma gives bounds for matrix $2$-norm.
\begin{lemma}\label{lem: matrix_norm}
    For any given matrix $A = (a_1,a_2,...,a_q)\in\mathbb{R}^{p\times q}$, we have:
    \[
        \|a_j\|^2\leq \|A\|^2 \leq \sum_{i=1}^{q}\|a_i\|^2,\ \forall{j}\in\{1,2,...,q\}.
    \]
\end{lemma}
The following lemma is a common result in decentralized optimization (e.g., \cite[Lemma 4]{lian2017can}).
\begin{lemma}\label{lem: WJ}
Under Assumption \ref{assump: W} we have:
\[
    \left\|W^{k} - \frac{J_n}{n}\right\|\leq \rho^k,\ \forall k\in\{0,1,2,...\}.
\]
\end{lemma}
\begin{proof}
Assume $1=\lambda_1>\lambda_2\geq...\geq \lambda_n>-1$ are eigenvalues of $W$. Since $W^{k}J_n = J_nW^k$, we know $W^k$ and $J_n$ are simultaneously diagonalizable. Hence there exists an orthogonal matrix $P$ such that 
\[
    W^k = P\text{diag}(\lambda_i^k)P^{-1},\quad \frac{J_n}{n}=P\text{diag}(1,0,0,...,0)P^{-1},
\]
and thus:
\[
    \left\|W^{k} - \frac{J_n}{n}\right\|  = \|P(\text{diag}(\lambda_i^{k}) - \text{diag}(1,0,0,...,0)  )P^{-1}\|\leq \max\{|\lambda_2|^k, |\lambda_n|^k\} = \rho^k.
\]
\end{proof}

The following three lemmas are adopted from Lemma 2.2 in \cite{ghadimi2018approximation}: 
\begin{lemma}\label{lem: Phi_lip}
    (Hypergradient) Define $\Phi_i(x):= f_i(x, y^*(x))$, where $y^*(x) = \argmin_{y\in\mathbb{R}^{q}}g(x,y)$. Under Assumptions \ref{assump: lip} and \ref{assump: strong_convexity} we have:
    \[
        \nabla \Phi_i(x) = \nabla_x f_i(x,y^*(x)) - \nabla_{xy} g(x,y^*(x))\left[\nabla_y^2g(x,y^*(x))\right]^{-1}\nabla_y f_i(x,y^*(x)).
    \]
    Moreover, $\nabla\Phi_i$ is Lipschitz continuous:
    \[
        \|\nabla\Phi_i(x_1) - \nabla\Phi_i(x_2)\| \leq L_{\Phi}\|x_1 - x_2\|,
    \]
    with the Lipschitz constant given by:
    \[
        L_{\Phi} = L + \frac{2L^2 + L_{g,2}L_{f,0}^2 }{\mu} + \frac{LL_{f,0}L_{g,2}+L^3 + L_{g,2}L_{f,0}L }{\mu^2} + \frac{L_{g,2}L^2L_{f,0}}{\mu^3} = \Theta(\kappa^3).
    \]
\end{lemma}

\textbf{Remark:} if Assumption \ref{assump: data_similarity} does not hold, then this hypergradient is completely different from the local hypergradient:
\begin{equation}\label{local_hyper}
    \begin{aligned}
        \nabla f_i(x, y_i^*(x)) = \nabla_x f_i(x,y_i^*(x)) - \nabla_{xy} g_i(x,y_i^*(x))\left[\nabla_y^2g_i(x,y_i^*(x))\right]^{-1}\nabla_y f_i(x,y_i^*(x)),
    \end{aligned}
\end{equation}
where $y_i^*(x)=\argmin_{y\in\mathbb{R}^{q}}g_i(x,y)$.
\begin{lemma}\label{lem: hyper_surrogate}
    Define:
    \[
        \bar{\nabla} f_i(x, y) = \nabla_x f_i(x,y) - \nabla_{xy} g(x,y) \left[\nabla_y^2g(x,y)\right]^{-1}\nabla_y f_i(x,y).
    \]
    Under the Assumption \ref{assump: lip} we have:
    \[
        \|\bar{\nabla} f_i(x, y) - \bar{\nabla} f_i(\tilde{x}, \tilde{y}) \|\leq L_{f}\|(x,y) - (\tilde{x}, \tilde{y})\|,
    \]
    where the Lipschitz constant is given by:
    \[
        L_{f} = L + \frac{L^2}{\mu} + L_{f,0}\left[\frac{L_{g,2}}{\mu} + \frac{L_{g,2}L}{\mu^2} \right] = \Theta(\kappa^2).
    \]
\end{lemma}
\begin{lemma}\label{lem: y_star_lip}
    Suppose Assumptions \ref{assump: lip} and \ref{assump: strong_convexity} hold, then we have:
    \[
        \begin{aligned}
            \|y_i^*(x_1) - y_i^*(x_2)\|\leq \kappa \|x_1 - x_2\|,\quad \forall i\in\{1,2,...,n\}.
        \end{aligned}
    \]
\end{lemma}
These lemmas basically reveal some nice properties of functions in bilevel optimization under Assumptions \ref{assump: lip} and \ref{assump: strong_convexity}. We will make use of these lemmas in our theoretical analysis.

\begin{lemma}\label{lem: phi}
    If the iterates satisfy:
    \[
        \overline{x_{k+1}} = \overline{x_{k}} - \eta_x\overline{\partial\Phi(X_k)},\quad \text{where } 0<\eta_x\leq \frac{1}{L_{\Phi}},
    \]
    then we have the following inequality holds:
    \begin{equation}\label{ineq: DBO_T}
    \frac{1}{K+1}\sum_{k=0}^{K}\|\nabla \Phi(\overline{x_k})\|^2 \leq \frac{2}{\eta_x(K+1)}(\Phi(\overline{x_0}) - \inf_x\Phi(x)) + \frac{1}{K+1}\sum_{k=0}^{K}\|\overline{\partial \Phi(X_k)} - \nabla \Phi(\overline{x_{k}})\|^2.
    \end{equation}
\end{lemma}
\begin{proof}
Since $\Phi(x)$ is $L_{\Phi}$-smooth, we have:
\[
    \begin{aligned}
        &\Phi(\overline{x_{k+1}}) - \Phi(\overline{x_{k}})
        \leq \nabla\Phi(\overline{x_{k}})^{\mathsf{T}}(-\eta_x\overline{\partial \Phi(X_k)}) + \frac{L_{\Phi}\eta_x^2}{2}\|\overline{\partial\Phi(X_k)}\|^2 =\frac{L_{\Phi}\eta_x^2}{2}\|\overline{\partial \Phi(X_k)}\|^2 - \eta_x \nabla\Phi(\overline{x_{k}})^{\mathsf{T}}\overline{\partial \Phi(X_k)} \\
        =&\frac{L_{\Phi}\eta_x^2}{2}\|\overline{\partial \Phi(X_k)} - \nabla \Phi(\overline{x_{k}})\|^2 +  (\frac{L_{\Phi}\eta_x^2}{2} - \eta_x)\|\nabla\Phi(\overline{x_{k}})\|^2 + (L_{\Phi}\eta_x^2- \eta_x) \nabla\Phi(\overline{x_{k}})^{\mathsf{T}}(\overline{\partial \Phi(X_k)} - \nabla\Phi(\overline{x_{k}})) \\
        \leq &\frac{L_{\Phi}\eta_x^2}{2}\|\overline{\partial \Phi(X_k)} - \nabla \Phi(\overline{x_{k}})\|^2 +  (\frac{L_{\Phi}\eta_x^2}{2} - \eta_x)\|\nabla\Phi(\overline{x_{k}})\|^2 + (\eta_x - L_{\Phi}\eta_x^2)(\frac{1}{2}\|\overline{\partial \Phi(X_k)} - \nabla\Phi(\overline{x_{k}})\|^2 + \frac{1}{2}\|\nabla\Phi(\overline{x_{k}})\|^2  ) \\
        = &\frac{\eta_x}{2} \|\overline{\partial \Phi(X_k)} - \nabla \Phi(\overline{x_{k}})\|^2 -\frac{\eta_x}{2}\|\nabla\Phi(\overline{x_{k}})\|^2,
    \end{aligned}
\]
where the second inequality is due to Young's inequality and $\eta_x\leq \frac{1}{L_{\Phi}}$.
Therefore, we have:
\begin{equation}\label{eq: Phi_main1}
    \begin{aligned}
        \|\nabla\Phi(\overline{x_{k}})\|^2&\leq \frac{2}{\eta_x}(\Phi(\overline{x_{k}}) - \Phi(\overline{x_{k+1}})) + \|\overline{\partial \Phi(X_k)} - \nabla \Phi(\overline{x_{k}})\|^2.
    \end{aligned}
\end{equation}
Summing \eqref{eq: Phi_main1} over $k=0,\ldots,K$, yields:
\[
    \begin{aligned}
        \sum_{k=0}^{K}\|\nabla \Phi(\overline{x_k})\|^2 \leq \frac{2}{\eta_x}(\Phi(\overline{x_0}) - \Phi(\overline{x_{K+1}})) + \sum_{k=0}^{K}\|\overline{\partial \Phi(X_k)} - \nabla \Phi(\overline{x_{k}})\|^2,
    \end{aligned}
\]
which completes the proof.
\end{proof}
We have the following lemma which provides an upper bound for $E_K$:
\begin{lemma}\label{lem: H}
In each iteration, if we have $\overline{x_{k+1}} = \overline{x_{k}} - \eta_x\overline{\partial\Phi(X_k)}$, then the following inequality holds:
    \[
        E_K \leq 8nS_K + 4n\eta_x^2\sum_{j=0}^{K-1}\|\overline{\partial\Phi(X_{j})} -  \nabla\Phi(\overline{x_j})\|^2 + 4n\eta_x^2T_{K-1}.
    \]
\end{lemma}

\begin{proof}
By the definition of $E_K$, we have:
    \[
        \begin{aligned}
            E_K =& \sum_{j=1}^{K}\sum_{i=1}^{n}\|x_{i,j} - x_{i,j-1}\|^2 = \sum_{j=1}^{K}\sum_{i=1}^{n}\|x_{i,j} - \overline{x_j} + \overline{x_j} - \overline{x_{j-1}} + \overline{x_{j-1}} - x_{i,j-1}\|^2 \\
            = & \sum_{j=1}^{K}\sum_{i=1}^{n}\|q_{i,j} -\eta_x(\overline{\partial\Phi(X_{j-1})} - \nabla\Phi(\overline{x_{j-1}})) - \eta_x\nabla\Phi(\overline{x_{j-1}}) - q_{i,j-1}\|^2 \\ 
            \leq & 4\sum_{j=1}^{K}\sum_{i=1}^{n}(\|q_{i,j}\|^2 + \eta_x^2\|\overline{\partial\Phi(X_{j-1})} -  \nabla\Phi(\overline{x_{j-1}}))\|^2 + \eta_x^2\|\nabla\Phi(\overline{x_{j-1}})\|^2 + \|q_{i,j-1}\|^2) \\
            \leq & 4n\sum_{j=1}^{K}(\|Q_j\|^2 + \|Q_{j-1}\|^2 + \eta_x^2\|\overline{\partial\Phi(X_{j-1})} -  \nabla\Phi(\overline{x_{j-1}})\|^2 + \eta_x^2\|\nabla\Phi(\overline{x_{j-1}})\|^2) \\
            \leq & 8nS_K + 4n\eta_x^2\sum_{j=0}^{K-1}(\|\overline{\partial\Phi(X_{j})} -  \nabla\Phi(\overline{x_{j}})\|^2 + \|\nabla\Phi(\overline{x_{j}})\|^2) \\
            = & 8nS_K + 4n\eta_x^2\sum_{j=0}^{K-1}\|\overline{\partial\Phi(X_{j})} -  \nabla\Phi(\overline{x_{j}})\|^2 + 4n\eta_x^2T_{K-1},
        \end{aligned}
    \]
where the second inequality is by the definition of $Q_j$, the third inequality is by the definition of $S_K$ and $Q_0 = 0$, the last equality is by the definition of $T_{K-1}.$
\end{proof}

Next we give bounds for $A_K$ and $B_K$.
\begin{lemma}\label{lem: AB}
Under Assumption \ref{assump: data_similarity}, if we properly choose $\eta_y, T$ and $N$ in Algorithm \ref{algo:DBO} and \ref{algo:DBOGT} such that:
\begin{equation}\label{condition1}
    \eta_y\leq \frac{2}{\mu + L},\quad \delta_y^T < \frac{1}{3},\quad \delta_{\kappa}^N <\frac{1}{8\kappa},
\end{equation}
then the following inequalities hold:
\[
    A_K\leq 3\delta_y^T(c_1 + 2\kappa^2E_K),\quad B_K\leq 2c_2 + 2d_1A_{K-1} + 2d_2E_K,
\]
where the constants are defined as follows:
\begin{equation}\label{CD}
    \begin{aligned}
        &c_1 = \sum_{i=1}^{n}\|y_{i,0}^{(0)} - y_i^*(x_{i,0})\|^2,\ c_2 = \sum_{i=1}^{n}\|v_{i,0}^* - v_{i,0}^{(0)}\|^2, \\ 
        &d_1 = 4(1+\sqrt{\kappa})^2(\kappa + \frac{L_{g,2} L_{f,0}}{\mu^2})^2=\Theta(\kappa^5),\ d_2 = 2(\kappa^2 + \frac{2L_{f,0}\kappa}{\mu} + \frac{2L_{f,0}\kappa^2}{\mu})^2 = \Theta(\kappa^6).
    \end{aligned}
\end{equation}
\end{lemma}

\begin{proof}
Recall that $A_K := \sum_{j=0}^{K}\sum_{i=1}^{n}\|y_{i,j}^{(T)} - y_i^*(x_{i,j})\|^2$. In order to bound $A_K$, we  bound each term as follows.
\begin{equation}\label{ineq: inner_error}
    \begin{aligned}
        &\|y_{i,j}^{(T)} - y_i^*(x_{i,j})\|^2 = \|y_{i,j}^{(T-1)} - \eta_y\nabla_y g(x_{i,j}, y_{i,j}^{(T-1)}) - y_i^*(x_{i,j}) \|^2 \\
        =&\|y_{i,j}^{(T-1)} - y_i^*(x_{i,j})\|^2 - 2\eta_y\nabla_y g_i(x_{i,j},y_{i,j}^{(T-1)})^{\mathsf{T}}(y_{i,j}^{(T-1)} - y_i^*(x_{i,j})) + \eta_y^2\|\nabla_y g_i(x_{i,j},y_{i,j}^{(T-1)})\|^2 \\
        \leq&(1 - \frac{2\eta_y\mu L}{\mu + L})\|y_{i,j}^{(T-1)} - y_i^*(x_{i,j})\|^2 + \eta_y(\eta_y - \frac{2}{\mu + L})\|\nabla_y g_i(x_{i,j},y_{i,j}^{(T-1)})\|^2 \\
        \leq&(1 - \frac{2\eta_y\mu L}{\mu + L})\|y_{i,j}^{(T-1)} - y_i^*(x_{i,j})\|^2 \\
        \leq&(1 - \frac{2\eta_y\mu L}{\mu + L})^{T} \|y_{i,j}^{(0)} - y_i^*(x_{i,j})\|^2\\
        = & \delta_y^T \|y_{i,j}^{(0)} - y_i^*(x_{i,j})\|^2,
    \end{aligned}
\end{equation}
where the first inequality follows from the strong-convexity and smoothness of the function $g$ \cite[Theorem 2.1.12]{nesterov2018lectures}, and the second inequality holds since we have $\eta_y\leq \frac{2}{\mu + L}$. We further have:
\[
    \begin{aligned}
        &\|y_{i,j}^{(0)} - y_i^*(x_{i,j})\|^2 = \|y_{i,j-1}^{(T)} - y_i^*(x_{i,j-1}) + y_i^*(x_{i,j-1}) - y_i^*(x_{i,j})\|^2 \\
        \leq &2(\|y_{i,j-1}^{(T)} - y_i^*(x_{i,j-1})\|^2 + \|y_i^*(x_{i,j-1}) - y_i^*(x_{i,j})\|^2) \\
        \leq &2\delta_y^T \|y_{i,j-1}^{(0)} - y_i^*(x_{i,j-1})\|^2 + 2\kappa^2\|x_{i,j-1} - x_{i,j}\|^2 \\
        < &\frac{2}{3}\|y_{i,j-1}^{(0)} - y_i^*(x_{i,j-1})\|^2 + 2\kappa^2\|x_{i,j-1} - x_{i,j}\|^2,
    \end{aligned}
\]
where the second inequality is by \eqref{ineq: inner_error} and Lemma \ref{lem: y_star_lip}, and the last inequality is by the condition \eqref{condition1}. Taking summation on both sides, we get
\[
    \begin{aligned}
        &\sum_{j=1}^{K}\sum_{i=1}^{n}\|y_{i,j}^{(0)} - y_i^*(x_{i,j})\|^2 \\ 
        \leq &\frac{2}{3}\sum_{j=1}^{K}\sum_{i=1}^{n}\|y_{i,j-1}^{(0)} - y_i^*(x_{i,j-1})\|^2 + 2\kappa^2\sum_{j=1}^{K}\sum_{i=1}^{n}\|x_{i,j} - x_{i,j-1}\|^2 \\
        \leq &\frac{2}{3} \sum_{j=0}^{K}\sum_{i=1}^{n}\|y_{i,j}^{(0)} - y_i^*(x_{i,j})\|^2 + 2\kappa^2E_K \\
        \leq &\frac{2}{3} c_1 + \frac{2}{3}\sum_{j=1}^{K}\sum_{i=1}^{n}\|y_{i,j}^{(0)} - y_i^*(x_{i,j})\|^2 + 2\kappa^2E_K,
    \end{aligned}   
\]
which directly implies:
\begin{equation}\label{y0-ystar-sum-bound}
    \sum_{j=1}^{K}\sum_{i=1}^{n}\|y_{i,j}^{(0)} - y_i^*(x_{i,j})\|^2 \leq 2c_1 + 6\kappa^2E_K.
\end{equation}
Combining \eqref{ineq: inner_error} and \eqref{y0-ystar-sum-bound} leads to:
\[
    \begin{aligned}
        &A_K = \sum_{j=0}^{K}\sum_{i=1}^{n}\|y_{i,j}^{T} - y_i^*(x_{i,j})\|^2\leq \delta_y^T \sum_{j=0}^{K}\sum_{i=1}^{n}\|y_{i,j}^{(0)} - y_i^*(x_{i,j})\|^2 \\
        &\leq \delta_y^T (c_1 + 2c_1 + 6\kappa^2E_K) = 3\delta_y^T(c_1 + 2\kappa^2E_K).
    \end{aligned}
\]

We then consider the bound for $B_K$. Recall that:
\[
    v_{i,k}^* = \left[\nabla_y^2g_i(x_{i,k},y_i^*(x_{i,k}))\right]^{-1}\nabla_y f_i(x_{i,k},y_i^*(x_{i,k})),
\]
which is the solution of the linear system $\left[\nabla_y^2g_i(x_{i,k},y_i^*(x_{i,k})) \right]v = \nabla_y f_i(x_{i,k},y_i^*(x_{i,k}))$ in the AID-based approach in Algorithm \ref{algo: Hypergrad}. Note that $v_{i,k}^*$ is a function of $x_{i,k}$, and it is $(\kappa^2 + \frac{2L_{f,0}L}{\mu^2} + \frac{2L_{f,0}L\kappa}{\mu^2})$-Lipschitz continuous with respect to $x_{i,k}$ \cite{ji2021bilevel}. For each term in $B_K$, we have:
\[
    \begin{aligned}
        \|v_{i,j}^* - v_{i,j}^{(0)}\|^2 &\leq 2(\|v_{i,j-1}^* - v_{i,j-1}^{(N)}\|^2 + \|v_{i,j}^* - v_{i,j-1}^*\|^2 ) \\
        &\leq 4(1+\sqrt{\kappa})^2(\kappa + \frac{L_{g, 2} L_{f,0}}{\mu^2})^2\|y_{i,j-1}^{(T)} - y_i^*(x_{i,j-1})\|^2 \\
        &+ 4\kappa (\frac{\sqrt{\kappa} - 1}{\sqrt{\kappa} + 1})^{2N}\|v_{i,j-1}^* - v_{i,j-1}^{(0)}\|^2 + 2(\kappa^2 + \frac{2L_{f,0}L}{\mu^2} + \frac{2L_{f,0}L\kappa}{\mu^2})^2 \|x_{i,j} - x_{i,j-1}\|^2,
    \end{aligned}
\]
where the second inequality follows \cite[Lemma 4]{ji2021bilevel}. Taking summation over $i,j$, we get
\begin{equation}\label{B_K}
    \begin{aligned}
        \sum_{j=1}^{K}\sum_{i=1}^{n}\|v_{i,j}^* - v_{i,j}^{(0)}\|^2\leq d_1 A_{K-1} + 4\kappa \delta_\kappa^N B_{K-1} + d_2E_K\leq d_1 A_{K-1} + \frac{1}{2}B_{K} + d_2E_K,
    \end{aligned}
\end{equation}
where the last inequality holds since we pick $N$ such that $4\kappa \delta_\kappa^N <\frac{1}{2}$. Therefore, we can get:
\[
    B_K\leq c_2 + d_1 A_{K-1} + \frac{1}{2}B_{K} + d_2E_K\quad \Rightarrow\quad B_K\leq 2c_2 + 2d_1A_{K-1} + 2d_2E_K,
\]
which completes the proof.
\end{proof}

The following lemmas give bounds on $\sum\|\overline{\partial \Phi(X_k)} - \nabla \Phi(\overline{x_{k}})\|^2$ in \eqref{eq: Phi_main1}.We first consider the case when the Assumption \ref{assump: data_similarity} holds. In this case, the outer loop computes the hypergradient via AID based approach. Therefore, we borrow \cite[Lemma 3]{ji2021bilevel} and restate it as follows. 
\begin{lemma}\cite[Lemma 3]{ji2021bilevel}\label{lem: aid_approx}
    Suppose Assumptions \ref{assump: lip}, \ref{assump: strong_convexity}, and \ref{assump: data_similarity} hold, then we have:
    \[
        \|\hat{\nabla}f_i(x_{i,j}, y_{i,j}^{(T)}) - \nabla f_i(x, y_i^*(x)) \|^2\leq \Gamma\|y_i^*(x_{i,j}) - y_{i, j}^{(T)}\|^2 + 6L^2\kappa\left(\frac{\sqrt{\kappa} -1 }{\sqrt{\kappa} + 1}\right)^{2N}\|v_{i,j}^* - v_{i,j}^{(0)}\|^2,
    \]
    where the constant $\Gamma$ is
    \[
        \Gamma = 3L^2 + \frac{3L_{g,2}^2L_{f,0}}{\mu^2} + 6L^2(1+\sqrt{\kappa})^2\left(\kappa + \frac{L_{g,2}L_{f,0}}{\mu^2}\right)^2=\Theta(\kappa^5).
    \]
\end{lemma}
Next, we are ready to bound $\sum\|\overline{\partial \Phi(X_k)} - \nabla \Phi(\overline{x_{k}})\|^2$ when Assumption \ref{assump: data_similarity} holds. 
\begin{lemma}\label{lem: phiestimateerror}
    Under Assumption \ref{assump: data_similarity}, we have:
    \begin{equation}\label{eq: grad_difference_homo}
    \begin{aligned}
        &\sum_{k=0}^{K}\|\overline{\partial \Phi(X_k)} - \nabla \Phi(\overline{x_{k}})\|^2  \leq  2L_{\Phi}^2S_K + \frac{2\Gamma}{n}A_K + 
        \frac{12L^2\kappa}{n}\delta_\kappa^NB_K.
    \end{aligned}
    \end{equation}
\end{lemma}
\begin{proof}
Under Assumption \ref{assump: data_similarity} we know $g_i = g$, and thus from \eqref{eq: global_phi} and \eqref{eq: local_phi} we have
\[
    \nabla \Phi_i(\overline{x_k}) = \nabla f_i(\overline{x_{k}}, y^*(\overline{x_{k}})).
\]
Therefore, we have
\[
    \begin{aligned}
        &\|\overline{\partial \Phi(X_k)} - \nabla \Phi(\overline{x_{k}})\|^2 = \frac{1}{n^2}\left\|\sum_{i=1}^{n}\left[\hat{\nabla}f_i(x_{i,k}, y_{i,k}^{(T)}) - \nabla f_i(\overline{x_{k}}, y^*(\overline{x_{k}}))\right]\right\|^2 \\
        \leq& \frac{1}{n}\sum_{i=1}^{n}\|\hat{\nabla}f_i(x_{i,k}, y_{i,k}^{(T)}) - \nabla f_i(\overline{x_{k}}, y^*(\overline{x_{k}}))\|^2 \\
        \leq &\frac{2}{n}\sum_{i=1}^n(
        \|\hat{\nabla}f_i(x_{i,k}, y_{i,k}^{(T)}) - \nabla f_i(x_{i,k}, y_i^*(x_{i,k}) )\|^2 + \| \nabla f_i(x_{i,k}, y_i^*(x_{i,k})) - \nabla f_i(\overline{x_{k}}, y^*(\overline{x_{k}})) \|^2 ) \\
        \leq & \frac{2}{n}\sum_{i=1}^n (\Gamma\|y_i^*(x_{i,k}) - y_{i, k}^{(T)}\|^2 + 6L^2\kappa\delta_\kappa^N\|v_{i,k}^* - v_{i,k}^{(0)}\|^2 + \| \nabla f_i(x_{i,k}, y^*(x_{i,k})) - \nabla f_i(\overline{x_{k}}, y^*(\overline{x_{k}})) \|^2) \\
        \leq &\frac{2\Gamma}{n}\sum_{i=1}^n\|y_i^*(x_{i,k}) - y_{i, k}^{(T)}\|^2 + 
        \frac{12L^2\kappa}{n}\delta_\kappa^N\sum_{i=1}^{n}\|v_{i,k}^* - v_{i,k}^{(0)}\|^2 + 2L_{\Phi}^2\|Q_k\|^2 ,
    \end{aligned}
\]
where the first inequality follows from the convexity of $\|\cdot\|^2$, the third inequality follows from Lemma \ref{lem: aid_approx} and Assumption \ref{assump: data_similarity}, the last inequality is by Lemma \ref{lem: Phi_lip}:
\[
    \begin{aligned}
        &\| \nabla f_i(x_{i,k}, y^*(x_{i,k})) - \nabla f_i(\overline{x_{k}}, y^*(\overline{x_{k}})) \|^2=\|\nabla\Phi_i(x_{i,k}) - \nabla\Phi_i(\overline{x_k})\|^2 \\
        &\leq L_{\Phi}^2\|x_{i,k} - \overline{x_k}\|^2 = L_{\Phi}^2\|q_{i,k}\|^2\leq L_{\Phi}^2\|Q_k\|^2.
    \end{aligned}
\]

Taking summation on both sides, we get:
\begin{equation}\label{ineq: Phi_main2}
    \begin{aligned}
        &\sum_{k=0}^{K}\|\overline{\partial \Phi(X_k)} - \nabla \Phi(\overline{x_{k}})\|^2  \leq  2L_{\Phi}^2S_K + \frac{2\Gamma}{n}A_K + 
        \frac{12L^2\kappa}{n}\delta_\kappa^NB_K .
    \end{aligned}
\end{equation}
\end{proof}

We now consider the case when Assumption \ref{assump: data_similarity} does not hold. In this case, our target in the lower level problem is 
\begin{equation}\label{ll-target}
y^*(\overline{x_k}) = \argmin_{y}\frac{1}{n}\sum_{i=1}^{n}g_i(\overline{x_k}, y).
\end{equation}
However, the update in our decentralized algorithm (e.g. line 8 of Algorithm \ref{algo:DBO}) aims at solving
\begin{equation}\label{ll-true}
    \tilde{y}_k^* := \argmin_{y}\frac{1}{n}\sum_{i=1}^{n}g_i(x_{i,k}, y),
\end{equation}
which is completely different from our target \eqref{ll-target}. To resolve this problem, we introduce the following lemma to characterize the difference:

\begin{lemma}\label{lem: y_tilde_and_y_star}
The following inequality holds: 
    \[
        \|\tilde{y}_k^* - y^*(\overline{x_k})\|\leq \frac{\kappa}{n}\sum_{i=1}^{n}\|x_{i,k} - \overline{x_k}\| \leq \kappa \|Q_k\|.
    \]
\end{lemma}
\begin{proof}
    By optimality conditions of \eqref{ll-target} and \eqref{ll-true}, we have:
    \[
        \begin{aligned}
            \frac{1}{n}\sum_{i=1}^{n}\nabla_yg_i(x_{i,k}, \tilde{y}_k^*) = 0,\quad \frac{1}{n}\sum_{i=1}^{n}\nabla_yg_i(\overline{x_k}, y^*(\overline{x_k})) = 0.
        \end{aligned}
    \]
    Combining with the strongly convexity and the smoothness of $g_i$ yields:
    \[
        \begin{aligned}
            &\|\frac{1}{n}\sum_{i=1}^{n}\nabla_yg_i(\overline{x_k}, \tilde{y}_k^*)\| = \|\frac{1}{n}\sum_{i=1}^{n}\nabla_yg_i(\overline{x_k}, \tilde{y}_k^*) - \frac{1}{n}\sum_{i=1}^{n}\nabla_yg_i(\overline{x_k}, y^*(\overline{x_k}))\| \geq \mu \|\tilde{y}_k^* - y^*(\overline{x_k})\|, \\
            &\|\frac{1}{n}\sum_{i=1}^{n}\nabla_yg_i(\overline{x_k}, \tilde{y}_k^*)\| = \|\frac{1}{n}\sum_{i=1}^{n}\nabla_yg_i(\overline{x_k}, \tilde{y}_k^*) - \frac{1}{n}\sum_{i=1}^{n}\nabla_yg_i(x_{i,k}, \tilde{y}_k^*) \|\leq \frac{L}{n}\sum_{i=1}^{n}\|x_{i,k} - \overline{x_k}\|.
        \end{aligned}
    \]
    Therefore, we obtain the following inequality:
    \[
        \|\tilde{y}_k^* - y^*(\overline{x_k})\|\leq \frac{\kappa}{n} \sum_{i=1}^{n}\|x_{i,k} - \overline{x_k}\|=\frac{\kappa}{n} \sum_{i=1}^{n}\|q_{i,k}\| \leq \kappa \|Q_k\|,
    \]
    where the last inequality is by Lemma \ref{lem: matrix_norm}.
\end{proof}

Notice that in the inner loop of Algorithms \ref{algo:DBO}, \ref{algo:DBOGT} and \ref{algo:DSBO}, i.e., Lines 4-11 of Algorithms \ref{algo:DBO} and \ref{algo:DBOGT}, and Lines 4-10 of Algorithm \ref{algo:DSBO}, $y_{i,k}^{T}$ converges to $\tilde{y}_k^*$ and the rates are characterized by \cite{olshevsky2019non, pu2021distributed, xin2020decentralized, qu2017harnessing} (e.g., Corollary 4.7. in \cite{olshevsky2019non}, Theorem 10 in \cite{nedic2017achieving} and Theorem 1 in \cite{pu2021distributed}). We include all the convergence rates here.
\begin{lemma}\label{lem: inner_error}
    Suppose Assumption \ref{assump: data_similarity} does not hold. We have:
    \begin{itemize}
        \item In Algorithm \ref{algo:DBO} there exists $\eta_y^{(t)} = \mathcal{O}(\frac{1}{t})$ such that $\frac{1}{n}\sum_{i=1}^{n}\|y_{i,k}^{(T)} - \tilde{y}_k^*\|^2 \leq \frac{C_1}{T}$;
        \item In Algorithm \ref{algo:DBOGT} there exists a constant $\eta_y$ such that $\frac{1}{n}\sum_{i=1}^{n}\|y_{i,k}^{(T)} - \tilde{y}_k^*\|^2\leq C_2 \alpha_1^T$;
        \item In Algorithm \ref{algo:DSBO} there exists $\eta_y^{(t)} = \mathcal{O}(\frac{1}{t})$ such that $\frac{1}{n}\sum_{i=1}^{n}\mathbb{E}\left[\|y_{i,k}^{(T)} - \tilde{y}_k^*\|^2\right]\leq \frac{C_3}{T}$.
    \end{itemize}
    Here $C_1, C_2, C_3$ are positive constant and $\alpha_1\in(0,1)$.
\end{lemma}
Besides, the JHIP oracle (Algorithm \ref{algo: JHI_oracle}) also performs standard decentralized optimization with gradient tracking in deterministic case (Algorithms \ref{algo:DBO} and \ref{algo:DBOGT}) and stochastic case (Algorithm \ref{algo:DSBO}). We have:
\begin{lemma}\label{lem: JHI_error}
In Algorithm \ref{algo: JHI_oracle}, we have:
    \begin{itemize}
        \item For deterministic case, there exists a constant $\gamma$ such that if $\gamma_t \equiv \gamma$ then $\|Z_i^{(t)}-Z^*\|^2 \leq C_4\alpha_2^t$. (See \cite{qu2017harnessing}).
        
        \item For stochastic case and there exists a diminishing stepsize sequence $\gamma_t = \mathcal{O}(\frac{1}{t})$, such that $\frac{1}{n}\sum_{i=1}^{n}\mathbb{E}\left[\|Z_i^{(t)}-Z^*\|^2\right] \leq \frac{C_5}{t}$. (See \cite{xin2020decentralized}).
    \end{itemize}
    Here $C_4, C_5$ are positive constant, and $\alpha_2\in (0,1)$. Here the optimal solution is denoted by $(Z^*)^{\mathsf{T}} = \left[\sum_{i=1}^{n}J_i\right]\left[\sum_{i=1}^{n}H_i\right]^{-1}$.
\end{lemma}
For simplicity we define:
\[
    C = \max\{C_1,C_2,C_3,C_4,C_5\},\quad \alpha = \max\{\alpha_1,\alpha_2\}.
\]
Since the objective functions mentioned in Lemma  \ref{lem: inner_error} (the lower level function $g$) and \ref{lem: JHI_error} (the objective in \eqref{eq: min_trace}) are strongly convex, we know $C$ and $\alpha$ only depend on $L, \mu, \rho$ and the stepsize (only when it is a constant). For example $\alpha_2$ in Lemma \ref{lem: JHI_error} only depends on the spectral radius of $H_i$, smallest eigenvalue of $H_i$, $\rho$ and $\gamma$.\\
For heterogeneous data on $g$ we have a different error estimation.
\begin{lemma}\label{lem:aid_error}
Suppose Assumption \ref{assump: data_similarity} does not hold. In Algorithm \ref{algo:DBO} and \ref{algo:DBOGT} we have:
\begin{equation}\label{lem: hat_and_tilde_nabla}
    \|\hat{\nabla}f_i(x_{i,k}, y_{i,k}^{(T)}) - \bar{\nabla} f_i(x_{i,k}, \tilde{y}_k^*)\|^2\leq 3L^2(1 + \kappa^2)\|y_{i,k}^{(T)}-\tilde{y}_k^*\|^2 + 3L_{f,0}^2C\alpha^N.
\end{equation}
\end{lemma}
\begin{proof}
    Note that:
    \[
        \begin{aligned}
            \hat{\nabla}f_i(x_{i,k}, y_{i,k}^{(T)}) &= \nabla_x f_i(x_{i,k},y_{i,k}^{(T)}) -\left[Z_{i,k}^{(N)}\right]^{\mathsf{T}} \nabla_y f_i(x_{i,k},y_{i,k}^{(T)}), \\
            \bar{\nabla} f_i(x_{i,k}, \tilde{y}_k^*) &= \nabla_x f_i(x_{i,k}, \tilde{y}_k^*) - \nabla_{xy} g(x_{i,k}, \tilde{y}_k^*)\nabla_y^2g(x_{i,k}, \tilde{y_k}^*)^{-1}\nabla_yf_i(x_{i,k},\tilde{y}_k^*).
        \end{aligned}
    \]
    Then we know:
    \[
        \begin{aligned}
            &\|\hat{\nabla}f_i(x_{i,k}, y_{i,k}^{(T)}) - \bar{\nabla} f_i(x_{i,k}, \tilde{y}_k^*))\|^2 \leq 3\|\nabla_x f_i(x_{i,k},y_{i,k}^{(T)}) - \nabla_x f_i(x_{i,k}, \tilde{y}_k^*)\|^2 \\
            &+ 3\left\|\left[ \left[Z_{i,k}^{(N)}\right]^{\mathsf{T}} - \nabla_{xy} g(x_{i,k}, \tilde{y}_k^*)\nabla_y^2g(x_{i,k}, \tilde{y_k}^*)^{-1} \right]  \nabla_y f_i(x_{i,k},y_{i,k}^{(T)})\right\|^2 \\ 
            &+3\left\|\nabla_{xy} g(x_{i,k}, \tilde{y}_k^*)\nabla_y^2g(x_{i,k}, \tilde{y_k}^*)^{-1}\left[\nabla_y f_i(x_{i,k},y_{i,k}^{(T)}) - \nabla_yf_i(x_{i,k},\tilde{y}_k^*)\right] \right\|^2 \\
            &\leq 3L^2\|y_{i,k}^{(T)} - \tilde{y}_k^*\|^2 + 3L_{f,0}^2C\alpha^N + 3L^2\kappa^2\|y_{i,k}^{(T)}-\tilde{y}_k^*\|^2 \\
            &=3L^2(1 + \kappa^2)\|y_{i,k}^{(T)}-\tilde{y}_k^*\|^2 + 3L_{f,0}^2C\alpha^N.
        \end{aligned}
    \]
    Notice that we have: $[Z_{i,k}^*]^{\mathsf{T}} = \nabla_{xy} g(x_{i,k}, \tilde{y}_k^*)\nabla_y^2g(x_{i,k}, \tilde{y_k}^*)^{-1}$, so the second inequality is by Assumption \ref{assump: lip} and Lemma \ref{lem: JHI_error}.
\end{proof}

\begin{lemma}
    Suppose Assumption \ref{assump: data_similarity} does not hold, then in Algorithms \ref{algo:DBO} and \ref{algo:DBOGT} we have:
    \begin{equation}\label{eq: grad_difference_hetero}
        \begin{aligned}
            \|\overline{\partial \Phi(X_k)} - \nabla \Phi(\overline{x_{k}})\|^2\leq 6(L^2(1 + \kappa^2))C\alpha^T + 6L_{f,0}^2C\alpha^N + 2(L_f^2 + L_f^2\kappa^2)\|Q_k\|^2.
        \end{aligned}
    \end{equation}
    
\end{lemma}

\begin{proof} We have
\[
    \begin{aligned}
        &\|\overline{\partial \Phi(X_k)} - \nabla \Phi(\overline{x_{k}})\|^2 = \frac{1}{n^2}\|\sum_{i=1}^{n}\left[\hat{\nabla}f_i(x_{i,k}, y_{i,k}^{(T)}) - \nabla f_i(\overline{x_{k}}, y^*(\overline{x_{k}}))\right]\|^2 \\
        \leq& \frac{1}{n}\sum_{i=1}^{n}\|\hat{\nabla}f_i(x_{i,k}, y_{i,k}^{(T)}) - \nabla f_i(\overline{x_{k}}, y^*(\overline{x_{k}}))\|^2 \\
        \leq &\frac{2}{n}\sum_{i=1}^n(
        \|\hat{\nabla}f_i(x_{i,k}, y_{i,k}^{(T)}) - \bar{\nabla} f_i(x_{i,k}, \tilde{y}_k^*))\|^2 + \| \bar{\nabla} f_i(x_{i,k}, \tilde{y}_k^*)- \nabla f_i(\overline{x_{k}}, y^*(\overline{x_{k}}))  \|^2) \\
        \leq & \frac{2}{n}\sum_{i=1}^{n}(3(L^2(1 + \kappa^2))\|y_{i,k}^{(T)}-\tilde{y}_k^*\|^2 + 3L_{f,0}^2C\alpha^N + L_{f}^2\|x_{i,k} - \overline{x_k}\|^2 + L_{f}^2\|\tilde{y}_k^* - y^*(\overline{x_k})\|^2)   \\
        \leq & \frac{2}{n}\sum_{i=1}^{n}(3(L^2(1 + \kappa^2))\|y_{i,k}^{(T)}-\tilde{y}_k^*\|^2 + 3L_{f,0}^2C\alpha^N + L_{f}^2\|q_{i,k}\|^2 + L_{f}^2\kappa^2\|Q_k\|^2) \\
        = & \frac{6(L^2(1 + \kappa^2))}{n}\sum_{i=1}^{n}\|y_{i,k}^{(T)}-\tilde{y}_k^*\|^2 + 6L_{f,0}^2C\alpha^N + 2(L_f^2 + L_f^2\kappa^2)\|Q_k\|^2,
    \end{aligned}
\]
where the third inequality is due to Lemma \ref{lem:aid_error} and Lemma \ref{lem: hyper_surrogate}, the last inequality is by Lemma \ref{lem: y_tilde_and_y_star}.
Notice that $\frac{1}{n}\sum_{i=1}^{n}\|y_{i,k}^{(T)} - \tilde{y}_k^*\|^2$ in the first term denotes the error of the inner loop iterates. In both DBO (Algorithm \ref{algo:DBO}) and DBOGT (Algorithm \ref{algo:DBOGT}), the inner loop performs a decentralized gradient descent with gradient tracking. By Lemma \ref{lem: inner_error}, we have the error bound $\frac{1}{n}\sum_{i=1}^{n}\|y_{i,k}^{(T)} - \tilde{y}_k^*\|^2\leq C\alpha^T$, which completes the proof.
\end{proof}

\subsection{Proof of the DBO convergence}\label{subsec: DBO_analysis}
In this section we will prove the following convergence result of the DBO algorithm:
\begin{theorem}\label{thm:dbo_crude}
    In Algorithm \ref{algo:DBO}, suppose Assumptions \ref{assump: lip}, \ref{assump: strong_convexity}, and \ref{assump: W} hold. If Assumption \ref{assump: data_similarity} holds, then by setting $\eta_x\leq\frac{1-\rho}{130\sqrt{n}L_{\Phi}} ,\ \eta_y\leq \frac{2}{\mu + L},\ T=\Theta(\log(\kappa)),\ N=\Theta(\log(\kappa))$, we have:
    \[
        \frac{1}{K+1}\sum_{j=0}^{K}\|\nabla\Phi(\overline{x_j})\|^2\leq \frac{4}{\eta_x(K+1)}(\Phi(\overline{x_0})-\inf_x \Phi(x)) + n\eta_x^2\cdot\frac{1272L_{\Phi}^2L_{f,0}^2(1+\kappa)^2}{(1-\rho)^2} + \frac{C_1}{K+1}.
    \]
    If Assumption \ref{assump: data_similarity} does not hold, then by setting $ \eta_x\leq \frac{1}{L_{\Phi}}, \eta_y^{(t)} = \mathcal{O}(\frac{1}{t})$, we have:
    \[
        \frac{1}{K+1}\sum_{j=0}^{K}\|\nabla\Phi(\overline{x_j})\|^2\leq \frac{2}{\eta_x(K+1)}(\Phi(\overline{x_0}) - \inf_x\Phi(x)) + \eta_x^2\frac{4nL_f^2(1+\kappa^2)L_{f,0}^2}{(1-\rho)^2}((1+\kappa)^2 + C\alpha^N)
        + \tilde{C_1},
    \]
    where 
    $C_1 = \Theta(1), C = \Theta(1)$ and $\tilde{C_1} = \mathcal{O}(\alpha^T + \alpha^N)$.
\end{theorem}

We first consider bounding the consensus error estimation for DBO:
\begin{lemma}\label{lem: S}
In Algorithm \ref{algo:DBO}, we have
\[
    \begin{aligned}
        S_{K} := \sum_{k=1}^{K}\|Q_k\|^2&< \frac{\eta_x^2}{(1-\rho)^2}\sum_{j=0}^{K-1}\sum_{i=1}^{n}\|\hat{\nabla}f_i(x_{i,j}, y_{i,j}^{(T)})\|^2. \\
    \end{aligned}
\]
\end{lemma}

\begin{proof}
Note that the $x$ update can be written as 
\[
    X_k = X_{k-1}W - \eta_x\partial\Phi(X_{k-1}),
\]
which indicates
\[
    \overline{x_k} = \overline{x_{k-1}} - \eta_x\overline{\partial\Phi(X_{k-1})}.
\]
By definition of $q_{i,k}$, we have
\[
    \begin{aligned}
        q_{i,k+1} &= \sum_{j=1}^{n}w_{ij}x_{j,k} - \eta_x \hat{\nabla}f(x_{i,k}, y_{i,k}^{(T)}) - (\overline{x_{k}} - \eta_x\overline{\partial \Phi(X_k)} ) \\
        &=\sum_{j=1}^{n}w_{ij}(x_{j,k} - \overline{x_{k}}) - \eta_x(\hat{\nabla}f(x_{i,k}, y_{i,k}^{(T)}) - \overline{\partial \Phi(X_k)}) \\
        &= Q_kWe_i - \eta_x\partial\Phi(X_k)(e_i - \frac{\mathbf{1}_n}{n})\quad \text{ -- $W$ is symmetric.} \\
    \end{aligned}
\]

Therefore, the update for the matrix $Q_{k+1}$ is
\[
    \begin{aligned}
        Q_{k+1} &= Q_kW - \eta_x\partial\Phi(X_k)(I - \frac{J_n}{n}) \\
        &=(Q_{k-1}W - \eta_x\partial\Phi(X_{k-1})(I - \frac{J_n}{n}))W - \eta_x\partial\Phi(X_k)(I - \frac{J_n}{n}) \\ 
        &=Q_0W^{k+1} - \eta_x\sum_{i=0}^{k}(\partial\Phi(X_{i})(I - \frac{J_n}{n})W^{k-i}) \\ 
        &=-\eta_x\sum_{i=0}^{k}\partial\Phi(X_{i})(W^{k-i} - \frac{J_n}{n}),
    \end{aligned}
\]
where the last equality is obtained by $Q_0 = 0$ and $J_n W = J_n.$ By Cauchy-Schwarz inequality, we have the following estimate
\[
    \begin{aligned}
        \|Q_{k+1}\|^2 &= \eta_x^2 \|\sum_{i=0}^{k} \partial\Phi(X_{i})(W^{k-i} - \frac{J_n}{n}) \|^2 \\
        &\leq \eta_x^2 (\sum_{i=0}^{k}\rho^{k-i}\|\partial \Phi(X_i)\|^2 ) (\sum_{i=0}^{k}\frac{1}{\rho^{k-i}}\|W^{k-i} - \frac{J_n}{n}\|^2) \\
        &\leq \eta_x^2 (\sum_{i=0}^{k}\rho^{k-i}\|\partial \Phi(X_i)\|^2)(\sum_{i=0}^{k}\rho^{k-i})  \\
        &< \frac{\eta_x^2}{1-\rho}(\sum_{i=0}^{k}\rho^{k-i}\|\partial \Phi(X_i)\|^2), \\
    \end{aligned}
\]
where the second inequality is obtained by Lemma \ref{lem: WJ}. Recall the definition of $\partial \Phi(X_k)$:
\[
    \partial \Phi(X_k) = [\hat{\nabla} f_1(x_{1,k}, y_{1,k}^{(T)}), \hat{\nabla} f_2(x_{2,k}, y_{2,k}^{(T)}),...,\hat{\nabla} f_n(x_{n,k}, y_{n,k}^{(T)}) ] \in \mathbb{R}^{p\times n}.
\]
Lemma \ref{lem: matrix_norm} further indicates:
\begin{equation}\label{PartialPhi}
    \begin{aligned}
        \|Q_{k+1}\|^2&< \frac{\eta_x^2}{1-\rho}(\sum_{j=0}^{k}\rho^{k-j}\|\partial \Phi(X_j)\|^2) \leq\frac{\eta_x^2}{1-\rho}(\sum_{j=0}^{k}\rho^{k-j}\sum_{i=1}^{n}\|\hat{\nabla}f_i(x_{i,j}, y_{i,j}^{(T)})\|^2 ) \\
        &=\frac{\eta_x^2}{1-\rho}\sum_{i=1}^{n}\sum_{j=0}^{k}\rho^{k-j}\|\hat{\nabla}f_i(x_{i,j}, y_{i,j}^{(T)})\|^2.
    \end{aligned}
\end{equation}
Summing \eqref{PartialPhi} over $k=0, \ldots, K-1$ yields
\begin{equation}\label{S_K}
    \begin{aligned}
        S_K &= \sum_{k=0}^{K-1}\|Q_{k+1}\|^2 < \frac{\eta_x^2}{(1-\rho)}\sum_{k=0}^{K-1}\sum_{i=1}^{n}\sum_{j=0}^k\rho^{k-j}\|\hat{\nabla}f_i(x_{i,j}, y_{i,j}^{(T)})\|^2 \\
        &=\frac{\eta_x^2}{(1-\rho)}\sum_{j=0}^{K-1}\sum_{i=1}^{n}\sum_{k=j}^{K-1}\rho^{k-j}\|\hat{\nabla}f_i(x_{i,j}, y_{i,j}^{(T)})\|^2
        < \frac{\eta_x^2}{(1-\rho)^2}\sum_{j=0}^{K-1}\sum_{i=1}^{n}\|\hat{\nabla}f_i(x_{i,j}, y_{i,j}^{(T)})\|^2,
    \end{aligned}
\end{equation}
where the second equality holds since we can change the order of summation.
\end{proof}

\subsubsection{Case 1: Assumption \ref{assump: data_similarity} holds}
We first consider the case when Assumption \ref{assump: data_similarity} holds. In this case, we obtain the following lemma.

\begin{lemma}\label{lem: aid_bound}
Suppose Assumptions \ref{assump: lip}, \ref{assump: strong_convexity}, and \ref{assump: data_similarity} hold, then we have:
    \[
        \|\hat{\nabla}f_i(x_{i,j}, y_{i,j}^{(T)})\|^2\leq 2(L^2\kappa\delta_\kappa^N\|v_{i,j}^{(0)}- v_{i,j}^*\|^2 + (1+\kappa)^2L_{f,0}^2).
    \]
\end{lemma}
\begin{proof}
    Notice that we have:
    \[
        \begin{aligned}
            &\|\hat{\nabla}f_i(x_{i,j}, y_{i,j}^{(T)})\|^2 \leq 2\|\hat{\nabla}f_i(x_{i,j}, y_{i,j}^{(T)}) - \bar{\nabla} f_i(x_{i,j}, y_{i,j}^{(T)})\|^2 + 2\|\bar{\nabla} f_i(x_{i,j}, y_{i,j}^{(T)})\|^2 \\
            \leq &2\|\nabla_{xy} g_i(x_{i,j},y_{i,j}^{(T)})(v_{i,j}^{(N)} - v_{i,j}^*)\|^2 \\
            + &2\|\nabla_x f_i(x_{i,j}, y_{i,j}^{(T)}) - \nabla_{xy} g_i(x_{i,j}, y_{i,j}^{(T)}) \left[\nabla_y^2g_i(x_{i,j}, y_{i,j}^{(T)})\right]^{-1}\nabla_y f_i(x_{i,j}, y_{i,j}^{(T)})\|^2 \\
            \leq &2(L^2\|v_{i,j}^{(N)} - v_{i,j}^*\|^2 + (L_{f,0} + L\frac{1}{\mu}L_{f,0})^2) \\
            \leq &2(L^2\kappa\delta_\kappa^N\|v_{i,j}^{(0)}- v_{i,j}^*\|^2 + (1+\kappa)^2L_{f,0}^2),
        \end{aligned}
    \]
    where the second inequality is via the Assumption \ref{assump: lip}, the last inequality is based on the convergence result of CG for the quadratic programming, e.g., eq. (17) in \cite{grazzi2020iteration}.
\end{proof}

Next we obtain the upper bound for $S_K$.
\begin{lemma} \label{lem: SKdbo1}
    Suppose Assumptions \ref{assump: lip}, \ref{assump: strong_convexity}, and \ref{assump: data_similarity} hold, then we have:
    \[
        S_K < \frac{2\eta_x^2}{(1-\rho)^2} (L^2\kappa \delta_\kappa^NB_{K-1} + nK(1+\kappa)^2L_{f,0}^2).
    \]
\end{lemma}
\begin{proof}
    By Lemmas \ref{lem: S} and \ref{lem: aid_bound}, we have:
    \[
        \begin{aligned}
            &S_K < \frac{\eta_x^2}{(1-\rho)^2}\sum_{j=0}^{K-1}\sum_{i=1}^{n}\|\hat{\nabla}f_i(x_{i,j}, y_{i,j}^{(T)})\|^2 \\
            \leq& \frac{\eta_x^2}{(1-\rho)^2}\sum_{j=0}^{K-1}\sum_{i=1}^{n}2(L^2\kappa\delta_\kappa^N\|v_{i,j}^{(0)}- v_{i,j}^*\|^2 + (1+\kappa)^2L_{f,0}^2) \\
            =&\frac{2\eta_x^2}{(1-\rho)^2} (L^2\kappa \delta_\kappa^NB_{K-1} + nK(1+\kappa)^2L_{f,0}^2),
        \end{aligned}
    \]
    which completes the proof.
\end{proof}

We are ready to prove the main results in Theorem \ref{thm:dbo_crude}. We first summarize main results in Lemmas \ref{lem: SKdbo1}, \ref{lem: AB} and \ref{lem: H}:
\begin{equation}\label{SAE}
    \begin{aligned}
        S_{K}&< \frac{2\eta_x^2}{(1-\rho)^2} (L^2\kappa \delta_\kappa^NB_{K-1} + nK(1+\kappa)^2L_{f,0}^2), \\
        A_K&\leq 3\delta_y^T(c_1 + 2\kappa^2E_K)\quad B_K\leq 2c_2 + 2d_1A_{K-1} + 2d_2E_K, \\
        E_K&\leq 8nS_K + 4n\eta_x^2\sum_{j=0}^{K-1}\|\overline{\partial\Phi(X_{j})} -  \nabla\Phi(\overline{x_{j}})\|^2 + 4n\eta_x^2T_{K-1}.
    \end{aligned}
\end{equation}
The next lemma proves the first part of Theorem \ref{thm:dbo_crude}.
\begin{lemma}
    Suppose the assumptions of Lemma \ref{lem: AB} hold. Further Set $N = \Theta(\log(\kappa)), T=\Theta(\log(\kappa)), \eta_x = \mathcal{O}(\frac{1}{\sqrt{n}\kappa^3})$ such that:
    \[
        \begin{aligned}
            &\delta_{\kappa}^N<\min\{ \frac{L_{\Phi}^2}{L^2\kappa(4d_1\kappa^2 + 2d_2)}, \kappa^{-6}\}=\Theta(\kappa^{-6}),\quad \delta_y^T<\min\{\frac{L_{\Phi}^2}{12\Gamma\kappa^2}, \kappa^{-5}\}=\Theta(\kappa^{-5}),\\
            &\eta_x< \frac{1-\rho}{130\sqrt{n}L_{\Phi}}\ (\text{which implies }\eta_x<\frac{(1-\rho)}{4\sqrt{2n}L_{\Phi}},\quad 106L_{\Phi}^2\cdot \frac{16nL_{\Phi}^2\eta_x^4}{(1-\rho)^2} + 52\eta_x^2L_{\Phi}^2 < \frac{1}{3} ),
        \end{aligned}
    \]
    
    we have:
    \[
        \frac{1}{K+1}\sum_{j=0}^{K}\|\nabla\Phi(\overline{x_j})\|^2\leq \frac{4}{\eta_x(K+1)}(\Phi(\overline{x_0})-\inf_x \Phi(x)) + \eta_x^2\cdot\frac{1272nL_{\Phi}^2L_{f,0}^2(1+\kappa)^2}{(1-\rho)^2} + \frac{C_1}{K+1},
    \]
    where the constant is given by:
    \[
        \begin{aligned}
            C_1 &= 106L_{\Phi}^2\cdot\frac{6\eta_x^2}{(1-\rho)^2}L^2\kappa\delta_{\kappa}^N(2c_2 + 2d_1c_1) +\frac{18L^2\kappa\delta_{\kappa}^N(2c_2 + 2d_1c_1) + 9\Gamma c_1\delta_y^T}{n}\\
            &= \Theta(\eta_x^2\delta_{\kappa}^{N}\kappa^{12} + \kappa^5\delta_y^T) = \Theta(1).
        \end{aligned}
    \]
\end{lemma}
\begin{proof}
    For $B_K$ we know:
    \begin{equation}\label{boundB}
        \begin{aligned}
            B_K&\leq 2c_2 + 2d_1A_K + 2d_2E_K \leq 2c_2 + \frac{2}{3}d_1(3c_1 + 6\kappa^2E_K) + 2d_2E_K \\
            &= 2c_2 + 2d_1c_1 + (4d_1\kappa^2 + 2d_2)E_K.
        \end{aligned}
    \end{equation}
    We first eliminate $B_K$ in the upper bound of $S_K$. Pick $N,T$ such that:
    \begin{equation}\label{deltaN}
        \delta_{\kappa}^N\cdot(4d_1\kappa^2 + 2d_2)\cdot L^2\kappa<L_{\Phi}^2\quad\Rightarrow\quad \delta_{\kappa}^N< \frac{L_{\Phi}^2}{L^2\kappa(4d_1\kappa^2 + 2d_2)}.
    \end{equation}
    
    Therefore, we have
    \[
        \begin{aligned}
            S_K&\leq \frac{2\eta_x^2}{(1-\rho)^2}(L^2\kappa\delta_{\kappa}^N(2c_2 + 2d_1c_1) + L^2\kappa\delta_{\kappa}^N(4d_1\kappa^2 + 2d_2)E_K + nK(1+\kappa)^2L_{f,0}^2) \\
            &\leq \frac{2\eta_x^2}{(1-\rho)^2}(L_{\Phi}^2 E_K + L^2\kappa\delta_{\kappa}^N(2c_2 + 2d_1c_1) + nK(1+\kappa)^2L_{f,0}^2).
        \end{aligned}
    \]

    Next we eliminate $E_K$ in this bound. By the definition of $\eta_x$, we know:
    \[
    \begin{aligned}
        \eta_x<\frac{(1-\rho)}{4\sqrt{2n}L_{\Phi}}\quad \Rightarrow
        \quad\frac{16n\eta_x^2L_{\Phi}^2}{(1-\rho)^2}<\frac{1}{2},
    \end{aligned}
    \]
    which, combining with \eqref{SAE}, yields
    \begin{equation}\label{SKafterE}
        \begin{aligned}
            S_{K}&\leq \frac{2\eta_x^2}{(1-\rho)^2}(L_{\Phi}^2(8nS_K + 4n\eta_x^2\sum_{j=0}^{K-1}\|\overline{\partial\Phi(X_{j})} -  \nabla\Phi(\overline{x_{j}})\|^2 + 4n\eta_x^2T_{K-1}) \\
            &+ L^2\kappa\delta_{\kappa}^N(2c_2 + 2d_1c_1) + nK(1+\kappa)^2L_{f,0}^2) \\
            &< \frac{1}{2}S_K + \frac{2\eta_x^2}{(1-\rho)^2}(4n\eta_x^2L_{\Phi}^2\sum_{j=0}^{K-1}\|\overline{\partial\Phi(X_{j})} -  \nabla\Phi(\overline{x_{j}})\|^2 + 4n\eta_x^2L_{\Phi}^2T_{K-1} \\
            &+ L^2\kappa\delta_{\kappa}^N(2c_2 + 2d_1c_1) + nK(1+\kappa)^2L_{f,0}^2). \\
        \end{aligned}
    \end{equation}
    The above inequality indicates
    \begin{equation}\label{SKafterE2}
        S_K\leq \frac{4\eta_x^2}{(1-\rho)^2}(4n\eta_x^2L_{\Phi}^2\sum_{j=0}^{K-1}\|\overline{\partial\Phi(X_{j})} -  \nabla\Phi(\overline{x_{j}})\|^2 + 4n\eta_x^2L_{\Phi}^2T_{K-1} + L^2\kappa\delta_{\kappa}^N(2c_2 + 2d_1c_1) + nK(1+\kappa)^2L_{f,0}^2). \\
    \end{equation}
    Note that we have
    \begin{equation}\label{deltaT}
        \delta_y^T<\frac{L_{\Phi}^2}{12\Gamma\kappa^2}\quad\Rightarrow\quad\delta_y^T\cdot 6\kappa^2\cdot 2\Gamma < L_{\Phi}^2.
    \end{equation}
    By Lemma \ref{lem: phiestimateerror},
    \[
    \begin{aligned}
        &\sum_{k=0}^{K}\|\overline{\partial \Phi(X_k)} - \nabla \Phi(\overline{x_{k}})\|^2\leq 2L_{\Phi}^2S_K + \frac{2\Gamma}{n}A_K + 
        \frac{12L^2\kappa}{n}\delta_\kappa^NB_K \\
        \leq & 2L_{\Phi}^2S_K + (\frac{2\Gamma}{n}\cdot 6\kappa^2\delta_y^T + \frac{12L^2\kappa}{n}\cdot\delta_{\kappa}^N\cdot(4d_1\kappa^2 + 2d_2))E_K + \frac{12L^2\kappa\delta_{\kappa}^N(2c_2 + 2d_1c_1) + 6\Gamma c_1\delta_y^T}{n} \\
        \leq & 2L_{\Phi}^2S_K + (\frac{L_{\Phi}^2}{n} + \frac{12L_{\Phi}^2}{n})E_K + \frac{12L^2\kappa\delta_{\kappa}^N(2c_2 + 2d_1c_1) + 6\Gamma c_1\delta_y^T}{n} \\
        \leq & 2L_{\Phi}^2S_K + \frac{13L_{\Phi}^2}{n}(8nS_K + 4n\eta_x^2\sum_{j=0}^{K-1}\|\overline{\partial\Phi(X_{j})} -  \nabla\Phi(\overline{x_{j}})\|^2 + 4n\eta_x^2T_{K-1}) + \frac{12L^2\kappa\delta_{\kappa}^N(2c_2 + 2d_1c_1) + 6\Gamma c_1\delta_y^T}{n}\\
        < &106L_{\Phi}^2S_K + 52\eta_x^2L_{\Phi}^2 (\sum_{j=0}^{K}\|\overline{\partial\Phi(X_{j})} -  \nabla\Phi(\overline{x_{j}})\|^2 + T_{K}) + \frac{12L^2\kappa\delta_{\kappa}^N(2c_2 + 2d_1c_1) + 6\Gamma c_1\delta_y^T}{n} \\
        \leq & (106L_{\Phi}^2\cdot \frac{16nL_{\Phi}^2\eta_x^4}{(1-\rho)^2} + 52\eta_x^2L_{\Phi}^2 )(\sum_{j=0}^{K}\|\overline{\partial\Phi(X_{j})} -  \nabla\Phi(\overline{x_{j}})\|^2 + T_K) \\
        + &106L_{\Phi}^2\cdot\frac{4\eta_x^2}{(1-\rho)^2}(L^2\kappa\delta_{\kappa}^N(2c_2 + 2d_1c_1) + nK(1+\kappa)^2L_{f,0}^2) + \frac{12L^2\kappa\delta_{\kappa}^N(2c_2 + 2d_1c_1) + 6\Gamma c_1\delta_y^T}{n}, \\
    \end{aligned}
    \]
    where the second inequality is by \eqref{SAE} and \eqref{boundB}, the third inequality is by \eqref{deltaN} and \eqref{deltaT}, the fourth inequality is obtained by \eqref{SAE} and the last inequality is by \eqref{SKafterE2}. Note that the definition of $\eta_x$ also indicates:
    \[
        106L_{\Phi}^2\cdot \frac{16nL_{\Phi}^2\eta_x^4}{(1-\rho)^2} + 52\eta_x^2L_{\Phi}^2 < \frac{1}{3}.
    \]
    Therefore, 
    \[
        \begin{aligned}
            &\sum_{k=0}^{K}\|\overline{\partial \Phi(X_k)} - \nabla \Phi(\overline{x_{k}})\|^2<\frac{1}{3}(\sum_{k=0}^{K}\|\overline{\partial \Phi(X_k)} - \nabla \Phi(\overline{x_{k}})\|^2 + T_K) \\
            +&106L_{\Phi}^2\cdot\frac{4\eta_x^2}{(1-\rho)^2}(L^2\kappa\delta_{\kappa}^N(2c_2 + 2d_1c_1) + nK(1+\kappa)^2L_{f,0}^2) + \frac{12L^2\kappa\delta_{\kappa}^N(2c_2 + 2d_1c_1) + 6\Gamma c_1\delta_y^T}{n},\\
        \end{aligned}
    \]
    which leads to
    \[
        \begin{aligned}
            \sum_{k=0}^{K}\|\overline{\partial \Phi(X_k)} - \nabla \Phi(\overline{x_{k}})\|^2&\leq \frac{1}{2}T_K +  106L_{\Phi}^2\cdot\frac{6\eta_x^2}{(1-\rho)^2}(L^2\kappa\delta_{\kappa}^N(2c_2 + 2d_1c_1) + nK(1+\kappa)^2L_{f,0}^2) \\
             &+\frac{18L^2\kappa\delta_{\kappa}^N(2c_2 + 2d_1c_1) + 9\Gamma c_1\delta_y^T}{n}.
        \end{aligned}
    \]
    
    Combining this bound with \eqref{ineq: DBO_T}, we can obtain
    \[
    \begin{aligned}
        T_{K}&\leq \frac{2}{\eta_x}(\Phi(\overline{x_0}) - \inf_{x}\Phi(x)) + \sum_{k=0}^{K}\|\overline{\partial \Phi(X_k)} - \nabla \Phi(\overline{x_{k}})\|^2 \\
        &\leq \frac{2}{\eta_x}(\Phi(\overline{x_0}) - \inf_{x}\Phi(x)) + \eta_x^2\cdot\frac{636nL_{\Phi}^2L_{f,0}^2(1+\kappa)^2}{(1-\rho)^2}K + \frac{1}{2}T_K + \frac{1}{2}C_1,
    \end{aligned}
    \]
    which implies
    \[
    \frac{1}{K+1}\sum_{j=0}^{K}\|\nabla\Phi(\overline{x_j})\|^2\leq \frac{4}{\eta_x(K+1)}(\Phi(\overline{x_0})-\inf_x \Phi(x)) + \eta_x^2\cdot\frac{1272nL_{\Phi}^2L_{f,0}^2(1+\kappa)^2}{(1-\rho)^2} + \frac{C_1}{K+1}.
    \]
    The constant $C_1$ is defined as
    \[
        \begin{aligned}
            C_1 &= 106L_{\Phi}^2\cdot\frac{6\eta_x^2}{(1-\rho)^2}L^2\kappa\delta_{\kappa}^N(2c_2 + 2d_1c_1) +\frac{18L^2\kappa\delta_{\kappa}^N(2c_2 + 2d_1c_1) + 9\Gamma c_1\delta_y^T}{n}\\
            &= \Theta(\eta_x^2\delta_{\kappa}^{N}\kappa^{12} + \kappa^5\delta_y^T) = \Theta(1).
        \end{aligned}
    \]
    Moreover, we notice that by setting $\eta_x = \Theta(K^{-\frac{1}{3}}n^{-\frac{1}{3}}\kappa^{-\frac{8}{3}})$, we have
    \[
        \frac{1}{K+1}\sum_{j=0}^{K}\|\nabla\Phi(\overline{x_j})\|^2 = \mathcal{O}(\frac{n^{\frac{1}{3}}\kappa^{\frac{8}{3}}}{K^{\frac{2}{3}}}).
    \]
\end{proof}

\subsubsection{Case 2: Assumption \ref{assump: data_similarity} does not hold}
Now we consider the case when Assumption \ref{assump: data_similarity} does not hold. 
\begin{lemma}\label{lem: hetero_S}
    \[
        \begin{aligned}
            S_{K} < \frac{\eta_x^2}{(1-\rho)^2}\sum_{j=0}^{K-1}\sum_{i=1}^{n}\|\hat{\nabla}f_i(x_{i,j}, y_{i,j}^{(T)})\|^2 < \frac{\eta_x^2L_{f,0}^2}{(1-\rho)^2}nK(2(1+\kappa)^2 + 2C\alpha^N).
        \end{aligned}
    \]
\end{lemma}
\begin{proof}
    The first inequality follows from Lemma \ref{lem: S}. For the second one observe that:
    \[
        \begin{aligned}
            &\|\hat{\nabla}f_i(x_{i,j},y_{i,j}^{(T)})\| = \|\nabla_x f_i(x_{i,k},y_{i,k}^{(T)}) -Z_{i,k}^{(N)} \nabla_y f_i(x_{i,k},y_{i,k}^{(T)})\| \\
            \leq &\|\nabla_x f_i(x_{i,k},y_{i,k}^{(T)})\| + \|(Z_{i,k}^{(N)} - Z_{i,k}^*)\nabla_y f_i(x_{i,k},y_{i,k}^{(T)})\| + \|Z_{i,k}^*\nabla_y f_i(x_{i,k},y_{i,k}^{(T)})\| \\
            \leq &(1 + \|Z_{i,k}^{(N)} - Z_{i,k}^*\| + \kappa)L_{f,0},
        \end{aligned}
    \]
    where we use $Z_{i,k}^{(N)}$ to denote the output of Algorithm \ref{algo: JHI_oracle} in outer loop iteration $k$ of agent $i$, and $Z_{i,k}^*$ denotes the optimal solution. By Cauchy-Schwarz inequality we know:
    \[
    \begin{aligned}
        \|\hat{\nabla}f_i(x_{i,j},y_{i,j}^{(T)})\|^2&\leq(1+\kappa+\|Z_{i,k}^{(N)}- Z_{i,k}^*\|)^2L_{f,0}^2\leq (2(1+\kappa)^2 + 2\|Z_{i,k}^{(N)} - Z_{i,k}^*\|^2)L_{f,0}^2\\
        &\leq (2(1+\kappa)^2+2C\alpha^N)L_{f,0}^2,
    \end{aligned}
    \]
    which completes the proof.
\end{proof}

Taking summation on both sides of \eqref{eq: grad_difference_hetero} and applying Lemma \ref{lem: hetero_S} and \ref{lem: inner_error} we know:
\[
    \begin{aligned}
        &\sum_{k=0}^{K}\|\overline{\partial \Phi(X_k)} - \nabla \Phi(\overline{x_{k}})\|^2 \leq 6(L^2(1 + \kappa^2))(K+1)C\alpha^T + (K+1)6L_{f,0}^2C\alpha^N + 2(L_f^2 + L_f^2\kappa^2)S_K \\
        \leq &(K+1)(6(L^2(1 + \kappa^2))C\alpha^T + 6L_{f,0}^2C\alpha^N) + 2L_f^2(1+\kappa^2)\frac{\eta_x^2L_{f,0}^2}{(1-\rho)^2}nK(2(1+\kappa)^2 + 2C\alpha^N).
    \end{aligned}
\]
Together with \eqref{ineq: DBO_T} we have:
\[
    \begin{aligned}
        &\frac{1}{K+1}\sum_{k=0}^{K}\|\nabla \Phi(\overline{x_k})\|^2\leq \frac{2}{\eta_x(K+1)}(\Phi(\overline{x_0}) - \inf_x\Phi(x)) + \frac{1}{K+1}\sum_{k=0}^{K}\|\overline{\partial \Phi(X_k)} - \nabla \Phi(\overline{x_{k}})\|^2 \\ 
        \leq &\frac{2}{\eta_x(K+1)}(\Phi(\overline{x_0}) - \inf_x\Phi(x)) + \eta_x^2\frac{4nL_f^2(1+\kappa^2)L_{f,0}^2 }{(1-\rho)^2}((1+\kappa)^2 + C\alpha^N)
        + \tilde{C_1},
    \end{aligned}
\]
where we define:
\[
    \tilde{C_1} = 6(L^2(1 + \kappa^2))C\alpha^T + 6L_{f,0}^2C\alpha^N = \mathcal{O}(\alpha^T + \alpha^N).
\]
Moreover, if we choose $\eta_x = \Theta(K^{-\frac{1}{3}}n^{-\frac{1}{3}}\kappa^{-\frac{8}{3}}),\ \eta_y^{(t)} = \mathcal{O}(\frac{1}{t}),\ T = \Theta(\log K),\ N=\Theta(\log K)$ then we can get:
\[
    \frac{1}{K+1}\sum_{j=0}^{K}\|\nabla\Phi(\overline{x_j})\|^2 = \mathcal{O}(\frac{n^{\frac{1}{3}}\kappa^{\frac{8}{3}}}{K^{\frac{2}{3}}}).
\]


\subsection{Proof of the convergence of DBOGT}\label{subsec: DBOGT_analysis}
In this section we will prove the following convergence result of Algorithm \ref{algo:DBOGT}
\begin{theorem}\label{thm:dbogt_crude}
    In Algorithm \ref{algo:DBOGT}, suppose Assumptions \ref{assump: lip}, \ref{assump: strong_convexity}, and \ref{assump: W} hold. If Assumption \ref{assump: data_similarity} holds, then by setting $\eta_x<\frac{(1-\rho)^2}{4\sqrt{2n}L_{\Phi}},\ \eta_y\leq \frac{2}{\mu + L},\ T = \Theta(\log(\kappa)),\ N = \Theta(\log(\kappa))$, we have:
    \[
        \frac{1}{K+1}\sum_{j=0}^{K}\|\nabla\Phi(\overline{x_j})\|^2\leq \frac{4}{\eta_x(K+1)}(\Phi(\overline{x_0})-\inf_x \Phi(x)) + \frac{C_2}{K+1}.
    \]
    If Assumption \ref{assump: data_similarity} does not hold, then by setting $\eta_x<\frac{(1-\rho)^2}{4\sqrt{6n}\kappa L_f},\ \eta_y=\Theta(1)$, we have:
    \[
        \frac{1}{K+1}\sum_{j=0}^{K}\|\nabla\Phi(\overline{x_j})\|^2\leq\frac{6}{\eta_x(K+1)}(\Phi(\overline{x_0}) - \inf_x\Phi(x))  + \frac{\|\partial\Phi(X_0)\|^2 }{K+1} +\tilde{C}_2.
    \]
    Here $C_2 = \Theta(1)$ and $\tilde{C}_2 = \Theta(\alpha^T + \alpha^N + \frac{1}{K+1})$.
\end{theorem}

We first bound the consensus estimation error in the following lemma.
\begin{lemma}\label{lem: S_gt}
    In Algorithm \ref{algo:DBOGT}, we have the following inequality holds:
    \[
    \begin{aligned}
        S_K&\leq\frac{\eta_x^2}{(1-\rho)^4}(\sum_{j=1}^{K-1}\sum_{i=1}^{n}\|\hat{\nabla} f_i(x_{i,j}, y_{i,j}^{(T)}) - \hat{\nabla} f_i(x_{i,j-1}, y_{i,j-1}^{(T)}) \|^2 + \|\partial\Phi(X_0)\|^2).\\
    \end{aligned}
    \]
\end{lemma}

\begin{proof}
From the updates of $x$ and $u$, we have:
\[
    \begin{aligned}
        \overline{u_k} = \overline{u_{k-1}} + \overline{\partial\Phi(X_k)} - \overline{\partial\Phi(X_{k-1})},\quad \overline{u_0} = \overline{\partial\Phi(X_0)},\quad
        \overline{x_{k+1}} = \overline{x_{k}} - \eta_x\overline{u_k},
    \end{aligned}
\]
which implies:
\[
    \begin{aligned}
        \overline{u_k} = \overline{\partial\Phi(X_k)},\quad \overline{x_{k+1}} = \overline{x_{k}} - \eta_x\overline{\partial\Phi(X_k)}.
    \end{aligned}
\]
Hence by definition of $q_{i,k+1}$:
\[
    \begin{aligned}
        q_{i,k+1} =& x_{i,k+1} - \overline{x_{k+1}} =\sum_{j=1}^{n}w_{ij}x_{j,k} - \eta_x u_{i,k} - \overline{x_{k}} + \eta_x\overline{u_k}=\sum_{j=1}^{n}w_{ij}(x_{j,k} - \overline{x_{k}}) - \eta_x(u_{i,k} - \overline{u_k})\\
        = &\sum_{j=1}^{n}w_{ij}q_{j,k} - \eta_xr_{i,k} =Q_{k}We_i - \eta_x R_k e_i. \\ 
    \end{aligned}
\]
Therefore, we can write the update of the matrix $Q_{k+1}$ as
\[
    \begin{aligned}
        Q_{k+1} = Q_kW - \eta_x R_k,\quad Q_1 = -\eta_x R_0.
    \end{aligned}
\]


Note that $Q_{k+1}$ takes the form of
\begin{equation}\label{eq: Q}
    Q_{k+1} = (Q_{k-1}W - \eta_x R_{k-1})W - \eta_x R_k = -\eta_x\sum_{i=0}^{k}R_iW^{k-i}.
\end{equation}
We then compute $r_{i,k}$ as following
\[
    \begin{aligned}
        &r_{i,k+1} = u_{i,k+1} - \overline{u_{k+1}} \\
        = &\sum_{j=1}^{n}w_{ij}u_{j,k} + \hat{\nabla} f_i(x_{i,k+1},y_{i,k+1}^{(T)}) -  \hat{\nabla} f_i(x_{i,k},y_{i,k}^{(T)}) - \overline{u_k} - (\overline{\partial\Phi(X_{k+1})} - \overline{\partial\Phi(X_{k})})\\
        =&\sum_{j=1}^{n}w_{ij}(u_{j,k} - \overline{u_k}) + (\partial\Phi(X_{k+1}) - \partial\Phi(X_{k}))(e_i - \frac{\textbf{1}_n}{n}) \\
        = &R_kWe_i + (\partial\Phi(X_{k+1}) - \partial\Phi(X_{k}))(e_i - \frac{\textbf{1}_n}{n}).
    \end{aligned}
\]
The matrix $ R_{k+1} $ can be written as
\begin{equation}\label{Rkform}
    \begin{aligned}
        R_{k+1} &= R_kW + (\partial\Phi(X_{k+1}) - \partial\Phi(X_{k}))(I - \frac{J_n}{n})\\
        &=R_0W^{k+1} + \sum_{j=0}^k(\partial\Phi(X_{j+1}) - \partial\Phi(X_{j}))(I - \frac{J_n}{n})W^{k-j} \\
        &=\partial\Phi(X_0)(I - \frac{J_n}{n})W^{k+1} + \sum_{j=0}^k(\partial\Phi(X_{j+1}) - \partial\Phi(X_{j}))(I - \frac{J_n}{n})W^{k-j} \\
        &=\sum_{j=0}^{k+1}(\partial\Phi(X_{j}) - \partial\Phi(X_{j-1}))(I - \frac{J_n}{n})W^{k+1-j},
    \end{aligned}
\end{equation}
where the third equality holds because of the initialization $u_{i,0} = \hat{\nabla} f_i(x_{i,0},y_{i,0}^{(T)})$ and we denote $\partial\Phi(X_{-1})=0$. Plugging \eqref{Rkform} into \eqref{eq: Q} yields
\[
    \begin{aligned}
        Q_{k+1} &= -\eta_x\sum_{i=0}^{k}\sum_{j=0}^{i}(\partial\Phi(X_{j}) - \partial\Phi(X_{j-1}))(I - \frac{J_n}{n})W^{k-j} \\
        &=-\eta_x\sum_{j=0}^{k}\sum_{i=j}^{k}(\partial\Phi(X_{j}) - \partial\Phi(X_{j-1}))(W^{k-j} - \frac{J_n}{n}) \\
        &=-\eta_x\sum_{j=0}^{k}(k+1-j)(\partial\Phi(X_{j}) - \partial\Phi(X_{j-1}))(W^{k-j} - \frac{J_n}{n}),
    \end{aligned}
\]
where the second equality is obtained by $J_n W = J_n$ and switching the order of the summations. Therefore, we have
\begin{equation}\label{normQgt}
    \begin{aligned}
        \|Q_{k+1}\|^2 &= \eta_x^2\|\sum_{j=0}^{k}(k+1-j)(\partial\Phi(X_{j}) - \partial\Phi(X_{j-1}))(W^{k-j} - \frac{J_n}{n})\|^2 \\
        &\leq \eta_x^2 (\sum_{j=0}^{k} \rho^{k-j}(k+1-j)\|\partial\Phi(X_{j}) - \partial\Phi(X_{j-1})\|^2)(\sum_{j=0}^{k}\frac{(k+1-j)}{\rho^{k-j}}\|W^{k-j}-\frac{J_n}{n}\|^2)\\
        &\leq \eta_x^2(\sum_{j=0}^{k} \rho^{k-j}(k+1-j)\|\partial\Phi(X_{j}) - \partial\Phi(X_{j-1})\|^2)(\sum_{j=0}^{k}(k+1-j)\rho^{k-j}) \\
        &< \frac{\eta_x^2}{(1-\rho)^2}(\sum_{j=0}^{k} \rho^{k-j}(k+1-j)\|\partial\Phi(X_{j}) - \partial\Phi(X_{j-1})\|^2), \\
    \end{aligned}
\end{equation}
where the first inequality is by Cauchy-Schwarz inequality, the second inequality is by Lemma \ref{lem: WJ} the last inequality uses the fact that:
\begin{equation}\label{ineq: geom_2_order}
    \sum_{j=0}^{k}(k+1-j)\rho^{k-j} = \sum_{m=0}^{k}(m+1)\rho^m = \frac{1 - (k+2)\rho^{k+1} + (k+1)\rho^{k+2}}{(1-\rho)^2}<\frac{1}{(1-\rho)^2}.
\end{equation}
Summing \eqref{normQgt} over $k=0,\ldots, K-1$, we get:
\[
    \begin{aligned}
         S_K &= \sum_{k=0}^{K-1}\|Q_{k+1}\|^2\leq
         \frac{\eta_x^2}{(1-\rho)^2}(\sum_{k=0}^{K-1}\sum_{j=0}^{k} \rho^{k-j}(k+1-j)\|\partial\Phi(X_{j}) - \partial\Phi(X_{j-1})\|^2) \\
         &=\frac{\eta_x^2}{(1-\rho)^2}(\sum_{j=0}^{K-1}\sum_{k=j}^{K-1} \rho^{k-j}(k+1-j)\|\partial\Phi(X_{j}) - \partial\Phi(X_{j-1})\|^2) \\
         &<\frac{\eta_x^2}{(1-\rho)^4}\sum_{j=0}^{K-1}\|\partial\Phi(X_{j}) - \partial\Phi(X_{j-1})\|^2 \\
         &\leq \frac{\eta_x^2}{(1-\rho)^4}(\sum_{j=1}^{K-1}\sum_{i=1}^{n}\|\hat{\nabla} f_i(x_{i,j}, y_{i,j}^{(T)}) - \hat{\nabla} f_i(x_{i,j-1}, y_{i,j-1}^{(T)}) \|^2 + \|\partial\Phi(X_0)\|^2), \\
    \end{aligned}
\]
which completes the proof.
\end{proof}

\subsubsection{Case 1: Assumption \ref{assump: data_similarity} holds}
When Assumption \ref{assump: data_similarity} holds, we have the following lemmas.

\begin{lemma}
Under Assumption \ref{assump: data_similarity}, the following inequality holds for Algorithm \ref{algo:DBOGT}:
    \[
        \sum_{j=1}^{K-1}\sum_{i=1}^{n}\|\hat{\nabla} f_i(x_{i,j}, y_{i,j}^{(T)}) - \hat{\nabla} f_i(x_{i,j-1}, y_{i,j-1}^{(T)})\|^2\leq 6\Gamma A_{K-1} + 12L^2\kappa\delta_\kappa^NB_{K-1} + L_{\Phi}^2E_{K-1}.
    \]
    Moreover, we have:
    \[
        S_K \leq \frac{\eta_x^2}{(1-\rho)^4}(6\Gamma A_{K-1} + 12L^2\kappa\delta_\kappa^NB_{K-1} + L_{\Phi}^2E_{K-1} +  \|\partial\Phi(X_0)\|^2).
    \]
\end{lemma}

\begin{proof}
     For each term, we know that for $j\geq 1$:
    \[
        \begin{aligned}
            &\|\hat{\nabla} f_i(x_{i,j}, y_{i,j}^{(T)}) - \hat{\nabla} f_i(x_{i,j-1}, y_{i,j-1}^{(T)}) \|^2 \\
            \leq& 3(\|\hat{\nabla} f_i(x_{i,j}, y_{i,j}^{(T)}) - \nabla \Phi_i(x_{i,j})\|^2 + \|\nabla \Phi_i(x_{i,j}) - \nabla \Phi_i(x_{i,j-1}) \|^2\\
            &+ \|\nabla \Phi_i(x_{i,j-1}) -  \hat{\nabla} f_i(x_{i,j-1}, y_{i,j-1}^{(T)})\|^2) \\
            \leq& 3(\Gamma(\|y_i^*(x_{i,j}) - y_{i, j}^{(T)}\|^2 + \|y_i^*(x_{i,j-1}) - y_{i, j-1}^{(T)}\|^2)  + 6L^2\kappa\delta_\kappa^N(\|v_{i,j}^* - v_{i,j}^{(0)}\|^2 + \|v_{i,j-1}^* - v_{i,j-1}^{(0)}\|^2)\\
            &+ L_{\Phi}^2\|x_{i,j} - x_{i,j-1}\|^2),
        \end{aligned}
    \]
    where the last inequality uses Lemma \ref{lem: aid_approx} and Lemma \ref{lem: Phi_lip}. Taking summation ($j=1,2,...,K-1$ and $i=1,2,...,n$) on both sides, we have:
    \[
        \begin{aligned}
            \sum_{j=1}^{K-1}\sum_{i=1}^{n}\|\hat{\nabla} f_i(x_{i,j}, y_{i,j}^{(T)}) - \hat{\nabla} f_i(x_{i,j-1}, y_{i,j-1}^{(T)})\|^2 \leq 6\Gamma A_{K-1} + 12L^2\kappa\delta_\kappa^NB_{K-1} + L_{\Phi}^2E_{K-1}. \\ 
        \end{aligned}
    \]
    Together with Lemma \ref{lem: S_gt}, we can prove the second inequality for $S_K$.
\end{proof}

Together with Lemma \ref{lem: AB} and \ref{lem: H}, we have:
\begin{equation}\label{SAEgt}
    \begin{aligned}
        S_K&\leq \frac{\eta_x^2}{(1-\rho)^4}(6\Gamma A_{K-1} + 12L^2\kappa\delta_\kappa^NB_{K-1} + L_{\Phi}^2E_{K-1} +  \|\partial\Phi(X_0)\|^2) \\
        A_K &\leq \delta_y^T (3c_1 + 6\kappa^2E_K)\quad B_K\leq 2c_2 + 2d_1A_{K-1} + 2d_2E_K\\
        E_K &\leq 8nS_K + 4n\eta_x^2\sum_{j=0}^{K-1}\|\overline{\partial\Phi(X_{j})} -  \nabla\Phi(\overline{x_{j}})\|^2 + 4n\eta_x^2T_{K-1}.
    \end{aligned}
\end{equation}
Now we can obtain the following result.
\begin{lemma}
    Suppose Assumptions \ref{assump: lip}, \ref{assump: strong_convexity}, \ref{assump: W} and \ref{assump: data_similarity} hold. Set parameters as:
    \[
    \begin{aligned}
        &\delta_y^{T}<\min\{ \frac{L_{\Phi}^2}{72\kappa^2\Gamma}, \kappa^{-5}\}=\Theta(\kappa^{-5}),\quad \delta_{\kappa}^N<\min\{ \frac{L_{\Phi}^2}{24L^2\kappa(4d_1\kappa^2 + 2d_2)}, \kappa^{-4}\}=\Theta(\kappa^{-4}),\\
        &\eta_x<\frac{(1-\rho)^2}{4\sqrt{2n}L_{\Phi}}.
    \end{aligned}
    \]
    For Algorithm \ref{algo:DBOGT}, we have:
    \[
        \frac{1}{K+1}\sum_{k=0}^{K}\|\nabla\Phi(\overline{x_{k}})\|^2\leq \frac{1}{K+1}\sum_{j=0}^{K}\|\nabla\Phi(\overline{x_j})\|^2\leq \frac{4}{\eta_x(K+1)}(\Phi(\overline{x_0})-\inf_x \Phi(x)) + \frac{C_2}{K+1},
    \]
    where the constants are defined as:
    \[
    \begin{aligned}
        \frac{1}{2}C_2 =& \frac{24\eta_x^2L_{\Phi}^2}{(1-\rho)^4}(\|\partial\Phi(X_0)\|^2 + 18\Gamma c_1\delta_y^T + 12L^2\kappa\delta_{\kappa}^N(2c_2 + 2d_1c_1)) \\
        &+\frac{18L^2\kappa\delta_{\kappa}^N(2c_2 + 2d_1c_1) + 9\Gamma c_1\delta_y^T}{n} \\
        =&\Theta(\eta_x^2\kappa^6 + (\eta_x^2\kappa^6 + 1)(\kappa^5\delta_y^T + \kappa^4\delta_{\kappa}^N)) = \Theta(1).
    \end{aligned}
    \]
\end{lemma}
\begin{proof}
We first bound $B_K$ as
\[
\begin{aligned}
    B_K&\leq 2c_2 + 2d_1A_K + 2d_2E_K \leq 2c_2 + \frac{2}{3}d_1(3c_1 + 6\kappa^2E_K) + 2d_2E_K \\
    &= 2c_2 + 2d_1c_1 + (4d_1\kappa^2 + 2d_2)E_K.
\end{aligned}
\]
Next we eliminate $A_K$ and $B_K$ in the upper bound of $S_K$. Choose $N, T$ such that
\[
    \begin{aligned}
        &\delta_y^T\cdot 6\kappa^2\cdot 6\Gamma < \frac{L_{\Phi}^2}{2},\quad \delta_{\kappa}^N\cdot(4d_1\kappa^2 + 2d_2)\cdot 12L^2\kappa<\frac{L_{\Phi}^2}{2}, \\
    \end{aligned}
\]
 \text{which implies}
\begin{equation}\label{deltagt}
    \begin{aligned}
        &\delta_y^{T}< \frac{L_{\Phi}^2}{72\kappa^2\Gamma},\quad \delta_{\kappa}^N< \frac{L_{\Phi}^2}{24L^2\kappa(4d_1\kappa^2 + 2d_2)}.
    \end{aligned}
\end{equation}

By \eqref{SAEgt}, we have 
\[
    \begin{aligned}
        S_K\leq \frac{\eta_x^2}{(1-\rho)^4}(2L_{\Phi}^2E_{K-1} + \|\partial\Phi(X_0)\|^2 + 18\Gamma c_1\delta_y^T + 12L^2\kappa\delta_{\kappa}^N(2c_2 + 2d_1c_1)).
    \end{aligned}
\]
Next we eliminate $E_{K-1}$ in this bound. The definition of $\eta_x$ gives $\eta_x<\frac{(1-\rho)^2}{4\sqrt{2n}L_{\Phi}}$, which implies $\frac{16nL_{\Phi}^2\eta_x^2}{(1-\rho)^4}<\frac{1}{2}.$ 
Together with \eqref{SAEgt} and $E_{K-1}\leq E_K$, we have:
\[
    \begin{aligned}
        S_{K}&\leq \frac{\eta_x^2}{(1-\rho)^4}(2L_{\Phi}^2(8nS_K + 4n\eta_x^2\sum_{j=0}^{K-1}\|\overline{\partial\Phi(X_{j})} -  \nabla\Phi(\overline{x_{j}})\|^2 + 4n\eta_x^2T_{K-1}) \\
        &+ \|\partial\Phi(X_0)\|^2 + 18\Gamma c_1\delta_y^T + 12L^2\kappa\delta_{\kappa}^N(2c_2 + 2d_1c_1) ) \\
        &\leq \frac{1}{2}S_K + \frac{\eta_x^2}{(1-\rho)^4}(2L_{\Phi}^2(4n\eta_x^2\sum_{j=0}^{K}\|\overline{\partial\Phi(X_{j})} -  \nabla\Phi(\overline{x_{j}})\|^2 + 4n\eta_x^2T_{K})  \\
        &+ \|\partial\Phi(X_0)\|^2 + 18\Gamma c_1\delta_y^T + 12L^2\kappa\delta_{\kappa}^N(2c_2 + 2d_1c_1) ),
    \end{aligned}
\]
which immediately implies 
\begin{equation}\label{ineq: DBOGT_S_final}
    \begin{aligned}
        S_K & < \frac{2\eta_x^2}{(1-\rho)^4}(8n\eta_x^2L_{\Phi}^2(\sum_{j=0}^{K}\|\overline{\partial\Phi(X_{j})} -  \nabla\Phi(\overline{x_{j}})\|^2 + T_{K})  \\
        &+ \|\partial\Phi(X_0)\|^2 + 18\Gamma c_1\delta_y^T + 12L^2\kappa\delta_{\kappa}^N(2c_2 + 2d_1c_1) ).
    \end{aligned}
\end{equation}

Moreover, by \eqref{ineq: Phi_main2} we have
\[
    \begin{aligned}
        &\sum_{k=0}^{K}\|\overline{\partial \Phi(X_k)} - \nabla \Phi(\overline{x_{k}})\|^2\leq 2L_{\Phi}^2S_K + \frac{2\Gamma}{n}A_K + 
        \frac{12L^2\kappa}{n}\delta_{\kappa}^NB_K \\
        \leq & 2L_{\Phi}^2S_K + (\frac{L_{\Phi}^2}{6n} + \frac{L_{\Phi}^2}{2n})E_K + \frac{12L^2\kappa\delta_{\kappa}^N(2c_2 + 2d_1c_1) + 6\Gamma c_1\delta_y^T}{n} \\
        \leq & 2L_{\Phi}^2S_K + \frac{2L_{\Phi}^2}{3n}(8nS_K + 4n\eta_x^2\sum_{j=0}^{K-1}\|\overline{\partial\Phi(X_{j})} -  \nabla\Phi(\overline{x_{j}})\|^2 + 4n\eta_x^2T_{K-1})\\
        &+ \frac{12L^2\kappa\delta_{\kappa}^N(2c_2 + 2d_1c_1) + 6\Gamma c_1\delta_y^T}{n}\\
        < &8L_{\Phi}^2S_K + \frac{8\eta_x^2L_{\Phi}^2}{3}(\sum_{j=0}^{K}\|\overline{\partial\Phi(X_{j})} -  \nabla\Phi(\overline{x_{j}})\|^2 + T_{K}) + \frac{12L^2\kappa\delta_{\kappa}^N(2c_2 + 2d_1c_1) + 6\Gamma c_1\delta_y^T}{n} \\
        \leq & (8L_{\Phi}^2\cdot \frac{16nL_{\Phi}^2\eta_x^4}{(1-\rho)^4} + \frac{8\eta_x^2L_{\Phi}^2}{3} )(\sum_{j=0}^{K}\|\overline{\partial\Phi(X_{j})} -  \nabla\Phi(\overline{x_{j}})\|^2 + T_K) \\
        + &8L_{\Phi}^2\cdot\frac{2\eta_x^2}{(1-\rho)^4}(\|\partial\Phi(X_0)\|^2 + 18\Gamma c_1\delta_y^T + 12L^2\kappa\delta_{\kappa}^N(2c_2 + 2d_1c_1)) + \frac{12L^2\kappa\delta_{\kappa}^N(2c_2 + 2d_1c_1) + 6\Gamma c_1\delta_y^T}{n}, \\
    \end{aligned}
\]
where the second inequality is by \eqref{deltagt}, the third inequality uses \eqref{SAEgt}. Note that $\eta_x$ satisfies:
\[
    \eta_x<\frac{(1-\rho)^2}{4\sqrt{2n}L_{\Phi}}\quad \Rightarrow\quad 8L_{\Phi}^2\cdot \frac{16nL_{\Phi}^2\eta_x^4}{(1-\rho)^4} + \frac{8\eta_x^2L_{\Phi}^2}{3} < \frac{1}{3}.
\]
Therefore, we have:
\[
    \begin{aligned}
        &\sum_{k=0}^{K}\|\overline{\partial \Phi(X_k)} - \nabla \Phi(\overline{x_{k}})\|^2\leq \frac{1}{3}\sum_{k=0}^{K}\|\overline{\partial \Phi(X_k)} - \nabla \Phi(\overline{x_{k}})\|^2 + \frac{1}{3}T_K \\
        + &8L_{\Phi}^2\cdot\frac{2\eta_x^2}{(1-\rho)^4}(\|\partial\Phi(X_0)\|^2 + 18\Gamma c_1\delta_y^T + 12L^2\kappa\delta_{\kappa}^N(2c_2 + 2d_1c_1)) + \frac{12L^2\kappa\delta_{\kappa}^N(2c_2 + 2d_1c_1) + 6\Gamma c_1\delta_y^T}{n},
    \end{aligned}
\]
which leads to
\[
    \begin{aligned}
        \sum_{k=0}^{K}\|\overline{\partial \Phi(X_k)} - \nabla \Phi(\overline{x_{k}})\|^2&\leq \frac{1}{2}T_K + \frac{24\eta_x^2L_{\Phi}^2}{(1-\rho)^4}(\|\partial\Phi(X_0)\|^2 + 18\Gamma c_1\delta_y^T + 12L^2\kappa\delta_{\kappa}^N(2c_2 + 2d_1c_1)) \\
        &+\frac{18L^2\kappa\delta_{\kappa}^N(2c_2 + 2d_1c_1) + 9\Gamma c_1\delta_y^T}{n}.
    \end{aligned}
\]
Recall \eqref{ineq: DBO_T}, we have
\[
\begin{aligned}
    \frac{1}{K+1}T_K&\leq\frac{2}{\eta_x(K+1)}(\Phi(\overline{x_0})-\inf_x \Phi(x)) + \frac{1}{K+1}\sum_{k=0}^{K}\|\overline{\partial \Phi(X_k)} - \nabla \Phi(\overline{x_{k}})\|^2 \\
    &\leq \frac{2}{\eta_x(K+1)}(\Phi(\overline{x_0})-\inf_x \Phi(x)) + \frac{1}{2(K+1)}T_K + \frac{1}{2(K+1)}C_2. \\
\end{aligned}
\]
Therefore, we get
\[
    \frac{1}{K+1}\sum_{j=0}^{K}\|\nabla\Phi(\overline{x_j})\|^2\leq \frac{4}{\eta_x(K+1)}(\Phi(\overline{x_0})-\inf_x \Phi(x)) + \frac{C_2}{K+1},
\]
where the constant is defined as following
\[
\begin{aligned}
    \frac{1}{2}C_2 &= \frac{24\eta_x^2L_{\Phi}^2}{(1-\rho)^4}(\|\partial\Phi(X_0)\|^2 + 18\Gamma c_1\delta_y^T + 12L^2\kappa\delta_{\kappa}^N(2c_2 + 2d_1c_1)) 
        +\frac{18L^2\kappa\delta_{\kappa}^N(2c_2 + 2d_1c_1) + 9\Gamma c_1\delta_y^T}{n} \\
        &=\Theta(\eta_x^2\kappa^6 + (\eta_x^2\kappa^6 + 1)(\kappa^5\delta_y^T + \kappa^4\delta_{\kappa}^N)) = \Theta(1).
\end{aligned}
\]
Finally, if we choose $\eta_x = \Theta(\kappa^{-3})$, then we have:
\[
    \frac{1}{K+1}\sum_{j=0}^{K}\|\nabla\Phi(\overline{x_j})\|^2 = \mathcal{O}(\frac{1}{K}).
\]
\end{proof}

\subsubsection{Case 2: Assumption \ref{assump: data_similarity} does not hold}
We first give a bound for $\|\tilde{y}_{j}^* - \tilde{y}_{j-1}^*\|$ in the following lemma.
\begin{lemma}\label{lem: y_tilde_and_y_tilde}
    Recall that $\tilde{y}_j^* = \argmin\frac{1}{n}\sum_{i=1}^{n}g_i(x_{i,j}, y)$. We have:
    \[
        \|\tilde{y}_j^* - \tilde{y}_{j-1}^*\|^2\leq \frac{\kappa^2}{n}\sum_{i=1}^{n}\|x_{i,j} - x_{i,j-1}\|^2.
    \]
\end{lemma}
\begin{proof}
    The proof technique is similar to Lemma \ref{lem: y_tilde_and_y_star}. Consider:
    \[
        \begin{aligned}
            &\|\frac{1}{n}\sum_{i=1}^{n}\nabla_yg_i(x_{i,j-1}, \tilde{y}_j^*)\| = \|\frac{1}{n}\sum_{i=1}^{n}\nabla_yg_i(x_{i,j-1}, \tilde{y}_j^*) - \frac{1}{n}\sum_{i=1}^{n}\nabla_yg_i(x_{i,j-1}, \tilde{y}_{j-1}^* )\| \geq \mu \|\tilde{y}_j^* - \tilde{y}_{j-1}^* \|, \\
            &\|\frac{1}{n}\sum_{i=1}^{n}\nabla_yg_i(x_{i,j-1}, \tilde{y}_j^*)\| = \|\frac{1}{n}\sum_{i=1}^{n}\nabla_yg_i(x_{i,j-1}, \tilde{y}_j^*) - \frac{1}{n}\sum_{i=1}^{n}\nabla_yg_i(x_{i,j}, \tilde{y}_j^*) \|\leq \frac{L}{n}\sum_{i=1}^{n}\|x_{i,j} - x_{i,j-1}\|,
        \end{aligned}
    \]
    which implies:
    \[
        \begin{aligned}
            \|\tilde{y}_j^* - \tilde{y}_{j-1}^* \|^2\leq \frac{\kappa^2}{n^2}(\sum_{i=1}^{n}\|x_{i,j} - x_{i,j-1}\|)^2\leq \frac{\kappa^2}{n}\sum_{i=1}^{n}\|x_{i,j}-x_{i,j-1}\|^2.
        \end{aligned}
    \]
\end{proof}

\begin{lemma}\label{lem: S_gt_hetero_data}
When the Assumption \ref{assump: data_similarity} does not hold, we have for Algorithm \ref{algo:DBOGT}:
    \[
        S_K \leq \frac{\eta_x^2}{(1-\rho)^4}\left[6(L^2(1 + \kappa^2))n(K-1)C\alpha^T + 6n(K-1)L_{f,0}^2C\alpha^N + 3L_f^2(1+\kappa^2)E_{K-1} + \|\partial\Phi(X_0)\|^2\right].
    \]
\end{lemma}

\begin{proof}
    We first consider:
    \[
        \begin{aligned}
            &\|\hat{\nabla} f_i(x_{i,j}, y_{i,j}^{(T)}) - \hat{\nabla} f_i(x_{i,j-1}, y_{i,j-1}^{(T)}) \|^2\leq 3\|\hat{\nabla} f_i(x_{i,j}, y_{i,j}^{(T)}) - \bar{\nabla} f_i(x_{i,j}, \tilde{y}_j^*)\|^2 \\
            &+3\|\bar{\nabla} f_i(x_{i,j}, \tilde{y}_j^*) - \bar{\nabla} f_i(x_{i,j-1}, \tilde{y}_{j-1}^*)\|^2 + 3\|\bar{\nabla} f_i(x_{i,j-1}, \tilde{y}_{j-1}^*) - \hat{\nabla} f_i(x_{i,j-1}, y_{i,j-1}^{(T)})\|^2 \\
            \leq &3(L^2(1 + \kappa^2))\|y_{i,j}^{(T)}-\tilde{y}_j^*\|^2 + 3L_{f,0}^2C\alpha^N + 3L_f^2(\|x_{i,j}-x_{i,j-1}\|^2 + \|\tilde{y}_j^* - \tilde{y}_{j-1}^*\|^2) \\
            &+ 3(L^2(1 + \kappa^2))\|y_{i,j-1}^{(T)}-\tilde{y}_{j-1}^*\|^2 + 3L_{f,0}^2C\alpha^N \\
            \leq& 3(L^2(1 + \kappa^2))(\|y_{i,j}^{(T)}-\tilde{y}_j^*\|^2 + \|y_{i,j-1}^{(T)}-\tilde{y}_{j-1}^*\|^2) + 6L_{f,0}^2C\alpha^N + 3L_f^2(\|x_{i,j}-x_{i,j-1}\|^2 \\&+ \frac{\kappa^2}{n}\sum_{i=1}^{n}\|x_{i,j}-x_{i,j-1}\|^2 ),
        \end{aligned}
    \]
    where the last inequality uses Lemma \ref{lem: y_tilde_and_y_tilde}. Taking summation on both sides, we have:
    \[
        \begin{aligned}
            &\sum_{j=1}^{K-1}\sum_{i=1}^{n}\|\hat{\nabla} f_i(x_{i,j}, y_{i,j}^{(T)}) - \hat{\nabla} f_i(x_{i,j-1}, y_{i,j-1}^{(T)})\|^2\\
            \leq& 3(L^2(1 + \kappa^2))\sum_{j=1}^{K-1}\sum_{i=1}^{n}(\|y_{i,j}^{(T)}-\tilde{y}_j^*\|^2 + \|y_{i,j-1}^{(T)}-\tilde{y}_{j-1}^*\|^2)\\
            & +6n(K-1)L_{f,0}^2C\alpha^N + 3L_f^2(1+\kappa^2)\sum_{j=1}^{K-1}\sum_{i=1}^{n}\|x_{i,j}-x_{i,j-1}\|^2 \\
            \leq& 6(L^2(1 + \kappa^2))n(K-1)C\alpha^T + 6n(K-1)L_{f,0}^2C\alpha^N + 3L_f^2(1+\kappa^2)E_{K-1},
        \end{aligned}
    \]
    which completes the proof together with Lemma \ref{lem: S_gt} and \ref{lem: inner_error}.
\end{proof}

\begin{lemma}\label{lem:phidiffgtnoassump}
    When the Assumption \ref{assump: data_similarity} does not hold, we further have for Algorithm \ref{algo:DBOGT}:
    \[
        \frac{1}{K+1}\sum_{k=0}^{K}\|\overline{\partial \Phi(X_k)} - \nabla  \Phi(\overline{x_{k}})\|^2\leq 6L^2(1+\kappa^2)C\alpha^T + 6L_{f,0}^2C\alpha^N + \frac{2L_f^2(1 + \kappa^2)}{K+1}S_K.
    \]
\end{lemma}

\begin{proof}
    By Lemma \ref{lem: phiestimateerror}, we have
    \[
        \begin{aligned}
            \sum_{k=0}^{K}\|\overline{\partial \Phi(X_k)} - \nabla  \Phi(\overline{x_{k}})\|^2&\leq 6L^2(1+\kappa^2)(K+1)C\alpha^T + 6L_{f,0}^2C\alpha^N(K+1) + 2(L_f^2 + L_f^2\kappa^2)S_K,
        \end{aligned}
    \]
    which yields
    \[
        \frac{1}{K+1}\sum_{k=0}^{K}\|\overline{\partial \Phi(X_k)} - \nabla  \Phi(\overline{x_{k}})\|^2\leq 6L^2(1+\kappa^2)C\alpha^T + 6L_{f,0}^2C\alpha^N + \frac{2L_f^2(1 + \kappa^2)}{K+1}S_K.
    \]
\end{proof}
Now we are ready to provide the convergence rate. Recall that from Lemma \ref{lem: S_gt_hetero_data}, \ref{lem: H} and inequality \eqref{ineq: DBO_T}, we have:
\begin{equation}\label{SEgt}
    \begin{aligned}
        &\frac{1}{K+1}\sum_{k=0}^{K}\|\nabla \Phi(\overline{x_k})\|^2\leq \frac{2}{\eta_x(K+1)}(\Phi(\overline{x_0}) - \inf_x\Phi(x)) + \frac{1}{K+1}\sum_{k=0}^{K}\|\overline{\partial \Phi(X_k)} - \nabla \Phi(\overline{x_{k}})\|^2, \\
        &S_K \leq \frac{\eta_x^2}{(1-\rho)^4}\left[6(L^2(1 + \kappa^2))n(K-1)C\alpha^T + 6n(K-1)L_{f,0}^2C\alpha^N + 3L_f^2(1+\kappa^2)E_{K-1} + \|\partial\Phi(X_0)\|^2\right], \\
        &E_K \leq 8nS_K + 4n\eta_x^2\sum_{j=0}^{K-1}\|\overline{\partial\Phi(X_{j})} -  \nabla\Phi(\overline{x_{j}})\|^2 + 4n\eta_x^2T_{K-1}.\\
    \end{aligned}
\end{equation}

The following lemma proves the convergence results in Theorem \ref{thm:dbogt_crude}.
\begin{lemma}
Suppose the Assumption \ref{assump: data_similarity} does not hold. We set $\eta_x$ as
\begin{equation}\label{etagtnoassum}
    \eta_x<\frac{(1-\rho)^2}{4\sqrt{6n}\kappa L_f}\quad (\text{which implies } \frac{\eta_x^2}{(1-\rho)^4}\cdot 3L_f^2(1+\kappa^2)\cdot 8n<\frac{1}{2}).
\end{equation}
Then we have:
    \[
        \frac{1}{K+1}\sum_{k=0}^{K}\|\nabla \Phi(\overline{x_k})\|^2\leq \frac{6}{\eta_x(K+1)}(\Phi(\overline{x_0}) - \inf_x\Phi(x))  + \frac{\|\partial\Phi(X_0)\|^2 }{K+1} +\tilde{C}_2,
    \]
    where the constant is given by:
\[
    \begin{aligned}
        \frac{\tilde{C}_2}{6} =& 6L^2(1+\kappa^2)C\alpha^T + 6L_{f,0}^2C\alpha^N \\
        &+ 2L_f^2(1 + \kappa^2)\cdot \frac{2\eta_x^2}{(1-\rho)^4}\left[6(L^2(1 + \kappa^2))nC\alpha^T + 6nL_{f,0}^2C\alpha^N + \frac{\|\partial\Phi(X_0)\|^2}{K+1}\right] \\
        =&\Theta(\alpha^T + \alpha^N + \frac{1}{K+1}).
    \end{aligned}
\]
\end{lemma}

\begin{proof}
    We first eliminate $E_{K-1}$ in the upper bound of $S_K$. Since $E_{K-1}\leq E_K$, we can plug $E_K$ into the upper bound of $S_K$ and get
    \begin{equation}\label{Skboundgtnoassum}
        \begin{aligned}
            S_K &\leq \frac{\eta_x^2}{(1-\rho)^4}\left[6(L^2(1 + \kappa^2))n(K-1)C\alpha^T + 6n(K-1)L_{f,0}^2C\alpha^N + 3L_f^2(1+\kappa^2)E_{K-1} + \|\partial\Phi(X_0)\|^2\right] \\
            &\leq \frac{\eta_x^2}{(1-\rho)^4}\left[6(L^2(1 + \kappa^2))n(K-1)C\alpha^T + 6n(K-1)L_{f,0}^2C\alpha^N + \|\partial\Phi(X_0)\|^2\right] \\
            &+ \frac{\eta_x^2}{(1-\rho)^4}\cdot 3L_f^2(1+\kappa^2)\cdot (8nS_K + 4n\eta_x^2\sum_{j=0}^{K-1}\|\overline{\partial\Phi(X_{j})} -  \nabla\Phi(\overline{x_{j}})\|^2 + 4n\eta_x^2T_{K-1}) \\
            &\leq \frac{\eta_x^2}{(1-\rho)^4}\left[6(L^2(1 + \kappa^2))n(K-1)C\alpha^T + 6n(K-1)L_{f,0}^2C\alpha^N + \|\partial\Phi(X_0)\|^2\right] \\ 
             &+ \frac{1}{2}(S_K + \frac{\eta_x^2}{2}\sum_{j=0}^{K-1}\|\overline{\partial\Phi(X_{j})} -  \nabla\Phi(\overline{x_{j}})\|^2 + \frac{\eta_x^2}{2}T_{K-1}), 
        \end{aligned}
    \end{equation}
    where the last step holds because of \eqref{etagtnoassum}. We rewrite \eqref{Skboundgtnoassum} as
    \[
        \begin{aligned}
            S_K&\leq \frac{2\eta_x^2}{(1-\rho)^4}\left[6(L^2(1 + \kappa^2))n(K-1)C\alpha^T + 6n(K-1)L_{f,0}^2C\alpha^N + \|\partial\Phi(X_0)\|^2\right] \\
            &+\frac{\eta_x^2}{2}\sum_{j=0}^{K-1}\|\overline{\partial\Phi(X_{j})} -  \nabla\Phi(\overline{x_{j}})\|^2 + \frac{\eta_x^2}{2}T_{K-1}.
        \end{aligned}
    \]
    By Lemma \ref{lem:phidiffgtnoassump}, we have 
    \begin{equation}\label{eq: DBOGT_main_tem1}
        \begin{aligned}
            &\frac{1}{K+1}\sum_{k=0}^{K}\|\overline{\partial \Phi(X_k)} - \nabla  \Phi(\overline{x_{k}})\|^2\leq 6(L^2+\kappa^2)C\alpha^T + 6L_{f,0}^2C\alpha^N + \frac{2L_f^2(1 + \kappa^2)}{K+1}S_K \\
            &\leq 6(L^2+\kappa^2)C\alpha^T + 6L_{f,0}^2C\alpha^N \\
            &+ \frac{2L_f^2(1 + \kappa^2)}{K+1}\cdot\{
            \frac{2\eta_x^2}{(1-\rho)^4}\left[6(L^2(1 + \kappa^2))n(K-1)C\alpha^T + 6n(K-1)L_{f,0}^2C\alpha^N + \|\partial\Phi(X_0)\|^2\right] \\
            &+\frac{\eta_x^2}{2}\sum_{j=0}^{K-1}\|\overline{\partial\Phi(X_{j})} -  \nabla\Phi(\overline{x_{j}})\|^2 + \frac{\eta_x^2}{2}T_{K-1}\} \\
            &\leq \frac{\tilde{C}_2}{6} + \frac{2L_f^2(1+\kappa^2)}{K+1}(\frac{\eta_x^2}{2}\sum_{j=0}^{K-1}\|\overline{\partial\Phi(X_{j})} -  \nabla\Phi(\overline{x_{j}})\|^2 + \frac{\eta_x^2}{2}T_{K-1}) \\
            &< \frac{\tilde{C}_2}{6} + \frac{1}{2(K+1)}\sum_{j=0}^{K-1}\|\overline{\partial\Phi(X_{j})} -  \nabla\Phi(\overline{x_{j}})\|^2 + \frac{1}{3(K+1)}T_{K-1} + \frac{\|\partial\Phi(X_0)\|^2 }{6(K+1)},\\
        \end{aligned}
    \end{equation}
    where the last inequality holds since we have
    \[
        \eta_x^2L_f^2(1+\kappa^2)=\frac{(1-\rho)^4}{24n}\cdot\frac{\eta_x^2}{(1-\rho)^4}\cdot 3L_f^2(1+\kappa^2)\cdot 8n<\frac{1}{2}\cdot\frac{(1-\rho)^4}{24n}<\frac{1}{48},
    \]
    and the constant is defined as:
    \[
    \begin{aligned}
        \frac{\tilde{C}_2}{6} =& 6L^2(1 +\kappa^2)C\alpha^T + 6L_{f,0}^2C\alpha^N \\
        &+ 2L_f^2(1 + \kappa^2)\cdot \frac{2\eta_x^2}{(1-\rho)^4}\left[6(L^2(1 + \kappa^2))nC\alpha^T + 6nL_{f,0}^2C\alpha^N + \frac{\|\partial\Phi(X_0)\|^2}{K+1}\right] \\
        =&\Theta(\alpha^T + \alpha^N + \frac{1}{K+1}).
    \end{aligned}
    \]
    We then rewrite \eqref{eq: DBOGT_main_tem1} as
    \[
        \frac{1}{K+1}\sum_{k=0}^{K}\|\overline{\partial \Phi(X_k)} - \nabla  \Phi(\overline{x_{k}})\|^2< \frac{\tilde{C}_2}{3} + \frac{2}{3}\frac{1}{K+1}T_{K-1} + \frac{\|\partial\Phi(X_0)\|^2 }{3(K+1)}.
    \]
    By Lemma \ref{lem: phiestimateerror}, we have 
    \[
        \begin{aligned}
            &\frac{1}{K+1}T_K = \frac{1}{K+1}\sum_{k=0}^{K}\|\nabla \Phi(\overline{x_k})\|^2\\
            &\leq \frac{2}{\eta_x(K+1)}(\Phi(\overline{x_0}) - \inf_x\Phi(x)) + \frac{1}{K+1}\sum_{k=0}^{K}\|\overline{\partial \Phi(X_k)} - \nabla \Phi(\overline{x_{k}})\|^2 \\
            &<\frac{2}{\eta_x(K+1)}(\Phi(\overline{x_0}) - \inf_x\Phi(x)) + \frac{\tilde{C}_2}{3} + \frac{2}{3}\frac{1}{K+1}T_{K-1} + \frac{\|\partial\Phi(X_0)\|^2 }{3(K+1)}.
        \end{aligned}
    \]
    Since $T_{K-1}\leq T_K$, we get
    \[
        \frac{1}{K+1}T_K<\frac{6}{\eta_x(K+1)}(\Phi(\overline{x_0}) - \inf_x\Phi(x))  + \frac{\|\partial\Phi(X_0)\|^2 }{K+1} +\tilde{C}_2.
    \]
Furthermore,  by setting $\eta_x=\Theta(\kappa^{-3}),\ \eta_y=\Theta(1),\ N = \Theta(\log K),\ T = \Theta(\log K)$, we have
\[
    \frac{1}{K+1}\sum_{j=0}^{K}\|\nabla\Phi(\overline{x_j})\|^2 = \mathcal{O}(\frac{1}{K}),
\]
which completes the proof.
\end{proof}

\subsection{Proof of the convergence of DSBO}\label{subsec: DSBO_analysis}
In this section we will prove the convergence result of the DSBO algorithm.
\begin{theorem}\label{thm:dsbo_crude}
    In Algorithm \ref{algo:DSBO}, suppose Assumptions \ref{assump: lip}, \ref{assump: strong_convexity}, and \ref{assump: W} hold. If Assumption \ref{assump: data_similarity} holds, then by setting $M = \Theta(\log K),\ T = \Theta(\log(\kappa)),\ \epsilon<\frac{1}{L},\ \eta_x\leq \frac{1}{L_{\Phi}},\ \eta_y\leq \frac{2}{\mu + L}$, we have:
    \[
        \frac{1}{K+1}\sum_{k=0}^{K}\mathbb{E}\left[\|\nabla\Phi(\overline{x_k})\|^2\right]\leq \frac{2}{\eta_x(K+1)}(\mathbb{E}\left[\Phi(\overline{x_0})\right]-\inf_x \Phi(x)) + \frac{3(\mu+ L)}{2\mu L}\eta_y L_f^2 \sigma_{g,1}^2 + \frac{3\eta_x^2nL_{\Phi}^2 }{(1-\rho)^2}\tilde{C}_f^2 + L\eta_x\tilde{\sigma}_f^2 + C_3.
    \]
     If Assumption \ref{assump: data_similarity} does not hold, then by setting $\eta_x\leq \frac{1}{L_{\Phi}},\ \eta_y^{(t)} = \mathcal{O}(\frac{1}{t})$, we have:
    \[
        \begin{aligned}
            \frac{1}{K+1}\sum_{k=0}^{K}\mathbb{E}\left[\|\nabla\Phi(\overline{x_k})\|^2\right]&\leq\frac{2}{\eta_x(K+1)}(\mathbb{E}\left[\Phi(\overline{x_{0}})\right] - \inf_x\Phi(x)) + \eta_x^2\frac{4nL_f^2(1+\kappa^2)}{(1-\rho)^2}((1+\kappa)^2+C\alpha^N) \\
            &+L\eta_x(4\sigma_f^2(1+\kappa^2) + (8L_{f,0}^2 + 4\sigma_f^2)\frac{C}{T})  + \tilde{C}_3.
        \end{aligned}
    \]
    Here $C = \Theta(1),\ C_3 = \Theta(\eta_x^2 + \frac{1}{K+1})$ and $\tilde{C}_3 = \mathcal{O}(\frac{1}{T}+\alpha^N)$.
\end{theorem}

We first define the following filtration:
\[
    \mathcal{F}_k = \sigma(\bigcup_{i=1}^{n}\{x_{i,0}, x_{i,1},...,x_{i,k}\}),\quad \mathcal{G}_{i,j}^{(t)} = \sigma(\{y_{i,l}^{(s)}: 0\leq l\leq j, 0\leq s\leq t\}\bigcup\{x_{i,l}: 0\leq l\leq j\}  ) 
\]
Then in both cases we have the following lemma.
\begin{lemma}
If $\eta_x\leq \frac{1}{L_{\Phi}}$, then we have:
    \[
    \begin{aligned}
        \mathbb{E}\left[\|\nabla\Phi(\overline{x_{k}})\|^2\right]\leq &\frac{2}{\eta_x}(\mathbb{E}\left[\Phi(\overline{x_{k}})\right] - \mathbb{E}\left[\Phi(\overline{x_{k+1}})\right]) + \mathbb{E}\left[\|\mathbb{E} \left[ \overline{\partial \Phi(X_k; \phi)}|\mathcal{F}_k\right] - \nabla\Phi(\overline{x_{k}})\|^2\right]\\
        +&L\eta_x \mathbb{E}\|\left[ \overline{\partial \Phi(X_k; \phi)}\right] - \mathbb{E} \left[ \overline{\partial \Phi(X_k; \phi)}|\mathcal{F}_k\right]\|^2.
    \end{aligned}
    \]
\end{lemma}

\begin{proof}
In each iteration of Algorithm \ref{algo:DSBO}, we have:
\begin{equation}\label{average_update}
    \overline{x_{k+1}} = \overline{x_k} - \eta_x \overline{\partial \Phi(X_k;\phi)}.
\end{equation}

The $L_{\Phi}$-smoothness of $\Phi$ indicates that
\[
    \Phi(\overline{x_{k+1}}) - \Phi(\overline{x_{k}}) \leq \nabla\Phi(\overline{x_{k}})^{\mathsf{T}}(-\eta_x\overline{\partial \Phi(X_k;\phi)}) + \frac{L_{\Phi}\eta_x^2}{2}\|\overline{\partial\Phi(X_k;\phi)}\|^2.
\]

Taking conditional expectation with respect to $\mathcal{F}_k$ on both sides, we have the following
\[
    \begin{aligned}
        &\mathbb{E}\left[\Phi(\overline{x_{k+1}})|\mathcal{F}_k \right] - \Phi(\overline{x_{k}}) \leq \nabla\Phi(\overline{x_{k}})^{\mathsf{T}}(-\eta_x\mathbb{E} \left[ \overline{\partial \Phi(X_k; \phi)}|\mathcal{F}_k\right]) + \frac{L_{\Phi}\eta_x^2}{2} \mathbb{E}\left[\|\overline{\partial\Phi(X_k;\phi)}\|^2|\mathcal{F}_k\right] \\
        =&-\frac{\eta_x}{2}(\|\nabla\Phi(\overline{x_{k}})\|^2 + \|\mathbb{E} \left[ \overline{\partial \Phi(X_k; \phi)}|\mathcal{F}_k\right]\|^2 - \|\mathbb{E} \left[ \overline{\partial \Phi(X_k; \phi)}|\mathcal{F}_k\right] - \nabla\Phi(\overline{x_{k}})\|^2) \\
        + &\frac{L_{\Phi}\eta_x^2}{2}(\|\mathbb{E} \left[ \overline{\partial \Phi(X_k; \phi)}|\mathcal{F}_k\right]\|^2 + \mathbb{E}\left[\| \overline{\partial \Phi(X_k; \phi)} - \mathbb{E} \left[ \overline{\partial \Phi(X_k; \phi)}|\mathcal{F}_k\right]\|^2|\mathcal{F}_k\right]   ) \\
        = & (\frac{L_{\Phi}\eta_x^2}{2} - \frac{\eta_x}{2})\|\mathbb{E} \left[ \overline{\partial \Phi(X_k; \phi)}|\mathcal{F}_k\right]\|^2 + \frac{L_{\Phi}\eta_x^2}{2} \mathbb{E}\left[\| \overline{\partial \Phi(X_k; \phi)} - \mathbb{E} \left[ \overline{\partial \Phi(X_k; \phi)}|\mathcal{F}_k\right]\|^2|\mathcal{F}_k\right] \\ 
        -& \frac{\eta_x}{2}(\|\nabla\Phi(\overline{x_{k}})\|^2 - \|\mathbb{E} \left[ \overline{\partial \Phi(X_k; \phi)}|\mathcal{F}_k\right] - \nabla\Phi(\overline{x_{k}})\|^2) \\
        \leq &\frac{L_{\Phi}\eta_x^2}{2} \mathbb{E}\left[\| \overline{\partial \Phi(X_k; \phi)} - \mathbb{E} \left[ \overline{\partial \Phi(X_k; \phi)}|\mathcal{F}_k\right]\|^2|\mathcal{F}_k\right] - \frac{\eta_x}{2}(\|\nabla\Phi(\overline{x_{k}})\|^2 - \|\mathbb{E} \left[ \overline{\partial \Phi(X_k; \phi)}|\mathcal{F}_k\right] - \nabla\Phi(\overline{x_{k}})\|^2),
    \end{aligned}
\]
where the second inequality holds since we pick $\eta_x\leq \frac{1}{L}$. Thus we can take expectation again and use tower property to obtain:
\begin{equation}\label{ineq: Phi_main1_stoc}
    \begin{aligned}
        \frac{\eta_x}{2}\mathbb{E}\left[\|\nabla\Phi(\overline{x_{k}})\|^2\right]\leq &\mathbb{E}\left[\Phi(\overline{x_{k}})\right] - \mathbb{E}\left[\Phi(\overline{x_{k+1}})\right] + \frac{\eta_x}{2}\mathbb{E}\left[\|\mathbb{E} \left[ \overline{\partial \Phi(X_k; \phi)}|\mathcal{F}_k\right] - \nabla\Phi(\overline{x_{k}})\|^2\right]\\
        +&\frac{L_{\Phi}\eta_x^2}{2} \mathbb{E}\|\left[ \overline{\partial \Phi(X_k; \phi)}\right] - \mathbb{E} \left[ \overline{\partial \Phi(X_k; \phi)}|\mathcal{F}_k\right]\|^2.
    \end{aligned}
\end{equation}
which completes the proof.
\end{proof}

\subsubsection{Case 1: Assumption \ref{assump: data_similarity} holds}
\begin{lemma}\label{lem: neumann_error}
    Under Assumption \ref{assump: data_similarity}, we have:
    \[
        \begin{aligned}
            \|\mathbb{E}\left[\hat{\nabla} f_i(x_{i,k},y_{i,k}^{(T)};\phi_{i,k})|\mathcal{F}_k\right] - \bar{\nabla} f_i(x_{i,k},y_{i,k}^{(T)})\| \leq L_{f,0}(1-\epsilon L)^M\kappa.
        \end{aligned}
    \]
\end{lemma}
\begin{proof}
    We first consider the expectation
    \begin{equation}\label{expectation_neumann}
        \begin{aligned}
            &\mathbb{E}\left[\hat{\nabla} f_i(x_{i,k},y_{i,k}^{(T)};\phi_{i,k})|\mathcal{F}_k\right] \\
            = &\nabla_xf_i(x_{i,k}, y_{i,k}^{(T)}) - \nabla_{xy}g(x_{i,k}, y_{i,k}^{(T)})\left[\epsilon \sum_{k=0}^{M-1} (I - \epsilon\nabla_y^2g(x_{i,k}, y_{i,k}^{(T)}) )^k\right]\nabla_yf_i(x_{i,k}, y_{i,k}^{(T)}),
        \end{aligned}
    \end{equation}
    where we use the fact that $M'$ in Algorithm \ref{algo: Hypergrad} is uniformly sampled from $[M]$. Notice that for the finite sum we have:
    \[
        \begin{aligned}
            \epsilon \sum_{k=0}^{M-1} (I - \epsilon\nabla_y^2g(x_{i,k}, y_{i,k}^{(T)}) )^k =& \epsilon \left[\epsilon\nabla_y^2g(x_{i,k}, y_{i,k}^{(T)})\right]^{-1}(I - (I - \epsilon\nabla_y^2g(x_{i,k}, y_{i,k}^{(T)}))^M) \\
            =&\left[\nabla_y^2g(x_{i,k}, y_{i,k}^{(T)})\right]^{-1}(I - (I - \epsilon\nabla_y^2g(x_{i,k}, y_{i,k}^{(T)}))^M),
        \end{aligned}
    \]
    which implies:
    \begin{equation}\label{Nserieserror}
        \begin{aligned}
            \|\epsilon \sum_{k=0}^{M-1} (I - \epsilon\nabla_y^2g(x_{i,k}, y_{i,k}^{(T)}) )^k - \left[\nabla_y^2g(x_{i,k}, y_{i,k}^{(T)})\right]^{-1}\|\leq \frac{ (1-\epsilon L)^M}{\mu}.
        \end{aligned}
    \end{equation}
    These inequalities yields the following error bound
    \[
        \begin{aligned}
            \|\mathbb{E}\left[\hat{\nabla} f_i(x_{i,k},y_{i,k}^{(T)};\phi_{i,k})|\mathcal{F}_k\right] - \bar{\nabla} f_i(x_{i,k},y_{i,k}^{(T)})\| \leq L_{f,0}(1-\epsilon L)^M\kappa,
        \end{aligned}
    \]
    which completes the proof.
\end{proof}

\begin{lemma}\label{lem: main_term1_stoc}
    Under the assumption \ref{assump: data_similarity}, we have:
    \begin{equation}\label{eq: main_term1_stoc}
        \sum_{k=0}^{K}\|\mathbb{E} \left[ \overline{\partial \Phi(X_k; \phi)}|\mathcal{F}_k\right] - \nabla\Phi(\overline{x_{k}})\|^2\leq 3((K+1)L_{f,0}^2(1-\epsilon L)^{2M}\kappa^2 + \frac{L_f^2}{n}A_K + L_{\Phi}^2S_K).
    \end{equation}
\end{lemma}
\begin{proof}
    We first bound each component of the gradient error as
    \[
        \begin{aligned}
            &\|\mathbb{E} \left[\hat{\nabla}f_i(x_{i,k}, y_{i,k}^{(T)};\phi_{i,k})|\mathcal{F}_k\right] - \nabla\Phi_i(\overline{x_{k}})\|^2 \\
            \leq & 3(\|\mathbb{E} \left[\hat{\nabla}f_i(x_{i,k}, y_{i,k}^{(T)};\phi_{i,k})|\mathcal{F}_k\right] - \bar{\nabla} f_i(x_{i,k},y_{i,k}^{(T)})\|^2 +\|\bar{\nabla} f_i(x_{i,k},y_{i,k}^{(T)}) - \nabla f_i(x_{i,k}, y_i^*(x_{i,k}))\|^2 \\
            + &\|\nabla f_i(x_{i,k}, y_i^*(x_{i,k})) - \nabla\Phi_i(\overline{x_{k}})\|^2) \\
            \leq & 3(L_{f,0}^2(1-\epsilon L)^{2M}\kappa^2 + L_f^2\|y_{i,k}^{(T)} - y_i^*(x_{i,k})\|^2 + L_{\Phi}^2\|x_{i,k} - \overline{x_{k}}\|^2),
        \end{aligned}
    \]
    where the second inequality is obtained by Lemma \ref{lem: neumann_error} and Lemma \ref{lem: Phi_lip}. Taking summation on both sides over $i=1,\ldots,n$, we have:
    \[
        \begin{aligned}
            &\|\mathbb{E} \left[ \overline{\partial \Phi(X_k; \phi)}|\mathcal{F}_k\right] - \nabla\Phi(\overline{x_{k}})\|^2 \leq \frac{1}{n}\sum_{i=1}^{n}\|\mathbb{E} \left[\hat{\nabla}f_i(x_{i,k}, y_{i,k}^{(T)};\phi_{i,k})|\mathcal{F}_k\right] - \nabla\Phi_i(\overline{x_{k}})\|^2 \\
            \leq & 3(L_{f,0}^2(1-\epsilon L)^{2M}\kappa^2 + \frac{L_f^2}{n} \sum_{i=1}^{n}\|y_{i,k}^{(T)} - y_i^*(x_{i,k})\|^2 + \frac{L_{\Phi}^2}{n} \sum_{i=1}^{n}\|x_{i,k} - \overline{x_{k}}\|^2). \\
        \end{aligned}
    \]
    Taking summation on both sides over $k=0,\ldots,K$, we know
    \[
        \begin{aligned}
            \sum_{k=0}^{K}\|\mathbb{E} \left[ \overline{\partial \Phi(X_k; \phi)}|\mathcal{F}_k\right] - \nabla\Phi(\overline{x_{k}})\|^2\leq 3((K+1)L_{f,0}^2(1-\epsilon L)^{2M}\kappa^2 + \frac{L_f^2}{n}A_K + L_{\Phi}^2S_K),
        \end{aligned}
    \]
    which completes the proof.
\end{proof}
The following Lemma is adopted from \cite[Lemma 5]{chen2021closing}.
\begin{lemma}\label{lem: main_term2_stoc}
    Under Assumptions \ref{assump: lip} - \ref{assump: stoc_derivatives}, we have:
    \[
    \begin{aligned}
        \mathbb{E}\|\mathbb{E}\left[\hat{\nabla} f_i(x_{i,k},y_{i,k}^{(T)};\phi_{i,k})|\mathcal{F}_k\right] -  \hat{\nabla} f_i(x_{i,k},y_{i,k}^{(T)};\phi_{i,k})\|^2&\leq \tilde{\sigma}_f^2, \\
        \mathbb{E}\|\left[ \overline{\partial \Phi(X_k; \phi)}\right] - \mathbb{E} \left[ \overline{\partial \Phi(X_k; \phi)}|\mathcal{F}_k\right]\|^2&\leq \tilde{\sigma}_f^2,
    \end{aligned}
    \]
    where the constants are defined as
    \[
        \tilde{\sigma}_f^2 = \sigma_f^2 + \frac{3}{\mu^2}\left[(\sigma_f^2 + L_{f,0}^2)(\sigma_{g,2}^2 + 2L^2) + \sigma_f^2L^2 \right] = \mathcal{O}(\kappa^2).
    \]
\end{lemma}

The following lemmas give the estimation bound of $A_K$ and $S_K$ in the stochastic case.
\begin{lemma}\label{lem: S_stoc}
    In Algorithm \ref{algo:DSBO}, we have
    \[
    \begin{aligned}
        \mathbb{E}\left[S_{K}\right]&< \frac{\eta_x^2}{(1-\rho)^2}\sum_{j=0}^{K-1}\sum_{i=1}^{n}\mathbb{E}\left[ \|\hat{\nabla}f_i(x_{i,j}, y_{i,j}^{(T)};\phi_{i,j})\|^2\right] \leq \frac{\eta_x^2nK}{(1-\rho)^2}\tilde{C}_f^2,
    \end{aligned}
\]
where the constant is defined as
\[
    \tilde{C}_f^2 = (L_{f,0} + \frac{LL_{f,1}}{\mu} +\frac{ LL_{f,1}}{\mu} )^2 + \tilde{\sigma}_f^2 = \mathcal{O}(\kappa^2).
\]
\end{lemma}
\begin{proof}
    Observe that in this stochastic case, we can replace $\hat{\nabla}f_i(x_{i,j}, y_{i,j}^{(T)})$ with $\hat{\nabla}f_i(x_{i,j}, y_{i,j}^{(T)};\phi_{i,j})$ in Lemma \ref{lem: S} to get the first inequality. For the second inequality, we adopt the bound in Lemma 2 of \cite{chen2021closing}.
\end{proof}
\begin{lemma}\label{lem: A_stoc}
    Set parameters in Algorithm \ref{algo:DSBO} as 
    \begin{equation}\label{etaydsbo1}
        \eta_y\leq \frac{2}{\mu + L},\quad \delta_y^{T}\leq \frac{1}{3}. 
    \end{equation}
    Then we have the following inequalities
    \[
        \begin{aligned}
            &\mathbb{E}\left[A_K\right] \leq \delta_y^{T}(2\mathbb{E}\left[c_1\right] + 6\kappa^2\mathbb{E}\left[E_K\right]) + \frac{\mu+ L}{2\mu L}\eta_y nK\sigma_{g,1}^2, \\
            &\mathbb{E}\left[E_K\right] \leq (\frac{6n}{(1-\rho)^2}+3)  \eta_x^2nK\tilde{C}_f^2.
        \end{aligned}
    \]
\end{lemma}
\begin{proof}
The proof is based on Lemma \ref{lem: AB}. Taking conditional expectation with respect to the filtration $\mathcal{G}_{i,j}^{(t-1)}$, we get
\[
    \begin{aligned}
        &\mathbb{E}\left[\|y_{i,j}^{(t)} - y_i^*(x_{i,j})\|^2|\mathcal{G}_{i,j}^{(t-1)}\right] = \mathbb{E}\left[\|y_{i,j}^{(t-1)} - \eta_y\nabla_y g(x_{i,j}, y_{i,j}^{(t-1)};\xi_{i,j}^{(t-1)}) - y_i^*(x_{i,j}) \|^2|\mathcal{G}_{i,j}^{(t-1)}\right] \\
        =& \|y_{i,j}^{(t-1)} - y_i^*(x_{i,j})\|^2 - 2\eta_y\nabla_y g_i(x_{i,j},y_{i,j}^{(t-1)})^{\mathsf{T}}(y_{i,j}^{(t-1)} - y_i^*(x_{i,j})) \\
        + &\eta_y^2\mathbb{E}\left[\|\nabla_y g_i(x_{i,j},y_{i,j}^{(t-1)};\xi_{i,j}^{(t-1)} )\|^2|\mathcal{G}_{i,j}^{(t-1)}\right] \\
        \leq&(1 - \frac{2\eta_y\mu L}{\mu + L})\|y_{i,j}^{(t-1)} - y_i^*(x_{i,j})\|^2 + \eta_y(\eta_y - \frac{2}{\mu + L})\|\nabla_y g_i(x_{i,j},y_{i,j}^{(t-1)})\|^2 + \eta_y^2\sigma_{g,1}^2\\
        \leq&(1 - \frac{2\eta_y\mu L}{\mu + L})\|y_{i,j}^{(t-1)} - y_i^*(x_{i,j})\|^2 + \eta_y^2\sigma_{g,1}^2,\\
    \end{aligned}
\]
where the first inequality is obtained by the smoothness and the strongly convexity of function $g$ and the second inequality is by \eqref{etaydsbo1}. Taking expectation on both sides and using the tower property, we have
\begin{equation}\label{ineq: inner_error_stoc}
    \begin{aligned}
        &\mathbb{E}\left[\|y_{i,j}^{(T)} - y_i^*(x_{i,j})\|^2\right] \leq (1 - \frac{2\eta_y\mu L}{\mu + L})\mathbb{E}\left[ \|y_{i,j}^{(T-1)} - y_i^*(x_{i,j})\|^2\right] + \eta_y^2\sigma_{g,1}^2 \\
        \leq &(1 - \frac{2\eta_y\mu L}{\mu + L})^{T}\mathbb{E}\left[\|y_{i,j}^{(0)} - y_i^*(x_{i,j})\|^2\right] + \eta_y^2\sigma_{g,1}^2\sum_{s=0}^{T-1}(1 - \frac{2\eta_y\mu L}{\mu + L})^s \\ 
        \leq &\delta_y^{T}\mathbb{E}\left[\|y_{i,j}^{(0)} - y_i^*(x_{i,j})\|^2\right] + \frac{\mu+ L}{2\mu L}\eta_y\sigma_{g,1}^2.
    \end{aligned}
\end{equation}
Moreover, by the warm-start strategy, we have
\begin{equation}\label{ytystardsbo}
    \begin{aligned}
        \mathbb{E}\left[\|y_{i,j}^{(0)} - y_i^*(x_{i,j})\|^2\right] &= \mathbb{E}\left[\|y_{i,j-1}^{T} - y_i^*(x_{i,j-1}) + y_i^*(x_{i,j-1}) - y_i^*(x_{i,j})\|^2\right] \\
        &\leq 2(\mathbb{E}\left[\|y_{i,j-1}^{T} - y_i^*(x_{i,j-1})\|^2\right] + \mathbb{E}\left[\|y_i^*(x_{i,j-1}) - y_i^*(x_{i,j})\|^2\right]) \\
        &\leq 2\delta_y^{T} \mathbb{E}\left[\|y_{i,j-1}^{(0)} - y_i^*(x_{i,j-1})\|^2\right] + 2\kappa^2\mathbb{E}\left[\|x_{i,j-1} - x_{i,j}\|^2\right] \\
        &\leq \frac{2}{3}\mathbb{E}\left[\|y_{i,j-1}^{(0)} - y_i^*(x_{i,j-1})\|^2\right] + 2\kappa^2\mathbb{E}\left[\|x_{i,j-1} - x_{i,j}\|^2\right],
    \end{aligned}
\end{equation}
where the second inequality is by Lemma \ref{lem: y_star_lip} and \eqref{ytystardsbo} and the last inequality is by \eqref{etaydsbo1}. Taking summation over $i,j$, we have:
\[
    \begin{aligned}
        \sum_{j=1}^{K}\sum_{i=1}^{n}\mathbb{E}\left[\|y_{i,j}^{(0)} - y_i^*(x_{i,j})\|^2\right] &\leq \frac{2}{3} \sum_{j=1}^{K}\sum_{i=1}^{n}\mathbb{E}\left[\|y_{i,j-1}^{(0)} - y_i^*(x_{i,j-1})\|^2\right] + 2\kappa^2\mathbb{E}\left[E_K\right] \\
        &\leq\frac{2}{3}\mathbb{E}\left[c_1\right] + \frac{2}{3}\sum_{j=1}^{K}\sum_{i=1}^{n}\mathbb{E}\left[\|y_{i,j}^{(0)} - y_i^*(x_{i,j})\|^2\right] + 2\kappa^2\mathbb{E}\left[E_K\right],
    \end{aligned}
\]
which leads to
\begin{equation}\label{ineq: initial_y_stoc}
    \sum_{j=1}^{K}\sum_{i=1}^{n}\mathbb{E}\left[\|y_{i,j}^{(0)} - y_i^*(x_{i,j})\|^2\right]\leq 2\mathbb{E}\left[c_1\right] + 6\kappa^2\mathbb{E}\left[E_K\right].
\end{equation}
Combining \eqref{ineq: initial_y_stoc} with \eqref{ineq: inner_error_stoc} and taking summation over $i,j$, we have
\[
    \begin{aligned}
        \mathbb{E}\left[A_K\right]&\leq \delta_y^{T}\sum_{j=1}^{K}\sum_{i=1}^{n}\mathbb{E}\left[\|y_{i,j}^{(0)} - y_i^*(x_{i,j})\|^2\right] + \frac{\mu+ L}{2\mu L}\eta_y nK\sigma_{g,1}^2 \\
        &\leq \delta_y^{T}(2\mathbb{E}\left[c_1\right] + 6\kappa^2\mathbb{E}\left[E_K\right]) + \frac{\mu+ L}{2\mu L}\eta_y nK\sigma_{g,1}^2.
    \end{aligned}
\]
Recall that for $E_K$ we have:
\[
        \begin{aligned}
            &E_K = \sum_{j=1}^{K}\sum_{i=1}^{n}\|x_{i,j} - x_{i,j-1}\|^2 = \sum_{j=1}^{K}\sum_{i=1}^{n}\|x_{i,j} - \overline{x_j} + \overline{x_j} - \overline{x_{j-1}} + \overline{x_{j-1}} - x_{i,j-1}\|^2 \\
            = & \sum_{j=1}^{K}\sum_{i=1}^{n}\|q_{i,j} -\eta_x\overline{\partial\Phi(X_{j-1};\phi)} - q_{i,j-1}\|^2 \\ 
            \leq & 3\sum_{j=1}^{K}\sum_{i=1}^{n}(\|q_{i,j}\|^2 + \eta_x^2\|\overline{\partial\Phi(X_{j-1};\phi)}\|^2 + \|q_{i,j-1}\|^2) \\
            \leq & 3n\sum_{j=1}^{K}(\|Q_j\|^2 + \|Q_{j-1}\|^2 + \eta_x^2\|\overline{\partial\Phi(X_{j-1};\phi)}\|^2 \\
            \leq & 6nS_K + 3n\eta_x^2\sum_{j=0}^{K-1}\|\overline{\partial\Phi(X_{j};\phi)}\|^2 \\
            \leq & 6nS_K + 3\eta_x^2\sum_{j=0}^{K-1}\sum_{i=1}^{n}\|\hat{\nabla} f_i(x_{i,j}, y_{i,j}^{(T)};\phi_{i,j})\|^2.
        \end{aligned}
\]
Taking expectation on both sides yields
\[
    \begin{aligned}
        \mathbb{E}\left[E_K\right]&\leq 6n\mathbb{E}\left[S_K\right] + 3\eta_x^2\sum_{j=0}^{K-1}\sum_{i=1}^{n}\mathbb{E}\left[\|\hat{\nabla} f_i(x_{i,j}, y_{i,j}^{(T)};\phi_{i,j})\|^2 \right] \\
        &\leq 6n\frac{\eta_x^2nK}{(1-\rho)^2}\tilde{C}_f^2 + 3\eta_x^2nK \tilde{C}_f^2 = (\frac{6n}{(1-\rho)^2}+3)  \eta_x^2nK\tilde{C}_f^2,\\
    \end{aligned}
\]
which completes the proof.
\end{proof}

Next, we prove the main convergence results in Theorem \ref{thm:dsbo_crude}. Taking expectation on both sides in \eqref{eq: main_term1_stoc}, we have:
\begin{equation}\label{ineq: term1_final_stoc}
    \begin{aligned}
        &\frac{1}{K+1}\sum_{k=0}^{K}\mathbb{E}\left[\|\mathbb{E} \left[ \overline{\partial \Phi(X_k; \phi)}|\mathcal{F}_k\right] - \nabla\Phi(\overline{x_{k}})\|^2\right]\\
        \leq &3(L_{f,0}^2(1-\epsilon L)^{2M}\kappa^2 + \frac{L_f^2}{n(K+1)}\mathbb{E}\left[A_K\right] + \frac{L_{\Phi}^2}{K+1}\mathbb{E}\left[S_K\right]) \\
        \leq &C_3 + \frac{3(\mu+ L)}{2\mu L}\eta_y L_f^2 \sigma_{g,1}^2 + \frac{3\eta_x^2nL_{\Phi}^2 }{(1-\rho)^2}\tilde{C}_f^2,
    \end{aligned}
\end{equation}
where the constant is defined as:
\[
    \begin{aligned}
        C_3 &= 3L_{f,0}^2(1-\epsilon L)^{2M}\kappa^2 +  \frac{L_f^2}{n(K+1)}\delta_y^{T}(2\mathbb{E}\left[c_1\right] + 6\kappa^2\mathbb{E}\left[E_K\right]) \\
        &\leq 3L_{f,0}^2(1-\epsilon L)^{2M}\kappa^2 +  \frac{L_f^2}{n(K+1)}\delta_y^{T}(2\mathbb{E}\left[c_1\right] + 6\kappa^2(\frac{6n}{(1-\rho)^2}+3)  \eta_x^2nK\tilde{C}_f^2)\\
        &=\Theta(\delta_{\epsilon}^{M}\kappa^2 + \eta_x^2\delta_y^T\kappa^8).
    \end{aligned}
\]
Here we denote $\delta_{\epsilon} = (1-\epsilon L)^2$ for simplicity. Therefore, we set $M = \Theta(\log K)$ and $T = \Theta(\log\kappa)$ such that $C_3 = \Theta(\eta_x^2 + \frac{1}{K+1})$. Recall that \eqref{ineq: Phi_main1_stoc} yields:
\[
    \begin{aligned}
        \mathbb{E}\left[\|\nabla\Phi(\overline{x_{k}})\|^2\right]\leq &\frac{2}{\eta_x}\mathbb{E}\left[\Phi(\overline{x_{k}})\right] - \mathbb{E}\left[\Phi(\overline{x_{k+1}})\right] + \mathbb{E}\left[\|\mathbb{E} \left[ \overline{\partial \Phi(X_k; \phi)}|\mathcal{F}_k\right] - \nabla\Phi(\overline{x_{k}})\|^2\right]\\
        +&L\eta_x\mathbb{E}\|\left[ \overline{\partial \Phi(X_k; \phi)}\right] - \mathbb{E} \left[ \overline{\partial \Phi(X_k; \phi)}|\mathcal{F}_k\right]\|^2.
    \end{aligned}
\]
If we take sum on both sides and calculate the average, together with inequality \eqref{ineq: term1_final_stoc} and Lemma \ref{lem: main_term2_stoc}, we have
\[
    \begin{aligned}
        \frac{1}{K+1}\sum_{k=0}^{K}\mathbb{E}\left[\|\nabla\Phi(\overline{x_{k}})\|^2\right]\leq & \frac{2}{\eta_x(K+1)}(\mathbb{E}\left[\Phi(\overline{x_0})\right]-\inf_x \Phi(x))\\
        &+ \frac{3(\mu+ L)}{2\mu L}\eta_y L_f^2 \sigma_{g,1}^2 + \frac{3\eta_x^2nL_{\Phi}^2 }{(1-\rho)^2}\tilde{C}_f^2 + L\eta_x\tilde{\sigma}_f^2 + C_3.
    \end{aligned}
\]
By setting $\eta_x = \Theta(K^{-\frac{1}{2}}),\ \eta_y = \Theta(K^{-\frac{1}{2}})$ we have
\[
    \frac{1}{K+1}\sum_{k=0}^{K}\mathbb{E}\left[\|\nabla\Phi(\overline{x_{k}})\|^2\right] = \mathcal{O}(\frac{1}{\sqrt{K}}).
\]

\subsubsection{Case 2: Assumption \ref{assump: data_similarity} does not hold}

We first consider the bound.
\begin{lemma}\label{lem: main_term1_stoc_case2}
Assume the Assumption \ref{assump: data_similarity} does not hold in Algorithm \ref{algo:DSBO}, we have
    \[
        \begin{aligned}
            &\frac{1}{K+1}\sum_{k=0}^{K}\mathbb{E}\left[\|\mathbb{E} \left[ \overline{\partial \Phi(X_k; \phi)}|\mathcal{F}_k\right] - \nabla\Phi(\overline{x_{k}})\|^2\right] \\
            \leq &\frac{6C(L^2(1 + \kappa^2))}{T} + 6L_{f,0}^2C\alpha^N + 2L_f^2(1+\kappa^2)\frac{\eta_x^2}{(1-\rho)^2}n(2(1+\kappa)^2 + 2C\alpha^N).
        \end{aligned}
    \]
\end{lemma}

\begin{proof}
    Denote by $\hat{Z}_{i,k}^{(N)}$ the output of each stochastic JHIP oracle \ref{algo: JHI_oracle} in Algorithm \ref{algo:DSBO}. Then
    \[
        \mathbb{E}\left[\hat{Z}_{i,k}^{(N)}\right] =Z_{i,k}^{(N)},
    \]
    which implies 
    \[
        \mathbb{E} \left[ \overline{\partial \Phi(X_k; \phi)}|\mathcal{F}_k\right] = \overline{\partial \Phi(X_k)}.
    \]
    
    Hence we can follow the same process in case 2 of DBO to get:
    \[
    \begin{aligned}
        &\sum_{k=0}^{K}\|\mathbb{E} \left[ \overline{\partial \Phi(X_k; \phi)}|\mathcal{F}_k\right] - \nabla\Phi(\overline{x_{k}})\|^2 = \sum_{k=0}^{K}\|\overline{\partial \Phi(X_k)} - \nabla \Phi(\overline{x_{k}})\|^2 \\
        \leq &(K+1)(\frac{6C(L^2(1 + \kappa^2))}{T} + 6L_{f,0}^2C\alpha^N) + 2L_f^2(1+\kappa^2)\frac{\eta_x^2}{(1-\rho)^2}nK(2(1+\kappa)^2 + 2C\alpha^N).
    \end{aligned}
    \]
    Taking expectation again on both sides and multiplying by $\frac{1}{K+1}$ on both sides, we complete the proof.
\end{proof}

The next lemma characterizes the variance of the gradient estimation.
\begin{lemma}\label{lem: main_term2_stoc_case2}
Assume the Assumption \ref{assump: data_similarity} does not hold in Algorithm \ref{algo:DSBO}, we have
    \[
        \mathbb{E}\|\left[ \overline{\partial \Phi(X_k; \phi)}\right] - \mathbb{E} \left[ \overline{\partial \Phi(X_k; \phi)}|\mathcal{F}_k\right]\|^2 \leq 4\sigma_f^2(1+\kappa^2) + (8L_{f,0}^2 + 4\sigma_f^2)\frac{C}{N}.
    \]
\end{lemma}

\begin{proof}
Recall that we have:
\[
    \begin{aligned}
        &\hat{\nabla}f_i(x_{i,k}, y_{i,k}^{(T)};\phi) = \nabla_x f_i(x_{i,k},y_{i,k}^{(T)};\phi_{i,k}^{(0)}) -\left[\hat{Z}_{i,k}^{(N)}\right]^{\mathsf{T}} \nabla_y f_i(x_{i,k},y_{i,k}^{(T)};\phi_{i,k}^{(0)}) \\
        &\hat{\nabla}f_i(x_{i,k}, y_{i,k}^{(T)}) = \nabla_x f_i(x_{i,k},y_{i,k}^{(T)}) -\left[Z_{i,k}^{(N)}\right]^{\mathsf{T}} \nabla_y f_i(x_{i,k},y_{i,k}^{(T)}). \\
    \end{aligned}
\]
By introducing intermediate terms we can have the following inequality
\[
    \begin{aligned}
        &\|\hat{\nabla}f_i(x_{i,k}, y_{i,k}^{(T)};\phi) - \hat{\nabla}f_i(x_{i,k}, y_{i,k}^{(T)})\|^2 \leq 4\|\nabla_x f_i(x_{i,k},y_{i,k}^{(T)};\phi_{i,k}^{(0)}) - \nabla_x f_i(x_{i,k},y_{i,k}^{(T)})\|^2\\
        +&4\|\left[\hat{Z}_{i,k}^{(N)}\right]^{\mathsf{T}} - \nabla_{xy} g(x_{i,k},y_{i,k}^{(T)} )\left[\nabla_y^2g(x_{i,k},y_{i,k}^{(T)})\right]^{-1}) \nabla_y f_i(x_{i,k},y_{i,k}^{(T)};\phi_{i,k}^{(0)})\|^2 \\
        +&4\|\nabla_{xy} g(x_{i,k},y_{i,k}^{(T)} )\left[\nabla_y^2g(x_{i,k},y_{i,k}^{(T)})\right]^{-1}(\nabla_y f_i(x_{i,k},y_{i,k}^{(T)};\phi_{i,k}^{(0)}) - \nabla_y f_i(x_{i,k},y_{i,k}^{(T)}))\|^2\\
        +&4\|(\nabla_{xy} g(x_{i,k},y_{i,k}^{(T)} )\left[\nabla_y^2g(x_{i,k},y_{i,k}^{(T)})\right]^{-1} - \left[Z_{i,k}^{(N)}\right]^{\mathsf{T}}) \nabla_y f_i(x_{i,k},y_{i,k}^{(T)})\|^2. \\
    \end{aligned}
\]
    For the first term and the third term we use $\mathbb{E}\left[\|\nabla f_i(x,y;\phi) - \nabla f_i(x,y)\|^2\right]\leq \sigma_f^2$. For the second term (and the fourth term) we use the fact that stochastic (and deterministic) decentralized algorithm achieves sublinear rate (Lemma \ref{lem: inner_error}). Without loss of generality we can set $C$ such that: $\max\{\frac{1}{n}\sum_{i=1}^{n}\mathbb{E}\left[\|\hat{Z}_{i,k}^{(N)} - Z_{i,k}^*\|^2\right], \|Z_{i,k}^{N}-Z_{i,k}^*\|^2\} \leq \frac{C}{N}$. For partial gradients in the second and fourth terms, we have:
    \[
        \begin{aligned}
            &\|\nabla_y f_i(x_{i,k},y_{i,k}^{(T)})\|^2\leq L_{f,0}^2 \\
            &\mathbb{E}\left[\|\nabla_y f_i(x_{i,k},y_{i,k}^{(T)};\phi_{i,k}^{(0)})\|^2\right] \\ 
            = &\mathbb{E}\left[\|\nabla_y f_i(x_{i,k},y_{i,k}^{(T)};\phi_{i,k}^{(0)}) - \nabla_y f_i(x_{i,k},y_{i,k}^{(T)})\|^2\right] + \|\nabla_y f_i(x_{i,k},y_{i,k}^{(T)})\|^2 \\
            \leq &\sigma_f^2 + L_{f,0}^2.
        \end{aligned}
    \]

    Taking summation and expectation on both sides, we have
    \[
        \begin{aligned}
            \frac{1}{n}\sum_{i=1}^{n}\mathbb{E}\left[\|\hat{\nabla}f_i(x_{i,k}, y_{i,k}^{(T)};\phi) - \hat{\nabla}f_i(x_{i,k}, y_{i,k}^{(T)})\|^2\right]\leq 4\sigma_f^2 +
            4(L_{f,0}^2 + \sigma_f^2)\frac{C}{N} + 4\frac{L^2}{\mu^2}\sigma_f^2 + 
            4L_{f,0}^2\frac{C}{N},
        \end{aligned}
    \]
    which, together with
    \[
        \mathbb{E}\|\left[ \overline{\partial \Phi(X_k; \phi)}\right] - \mathbb{E} \left[ \overline{\partial \Phi(X_k; \phi)}|\mathcal{F}_k\right]\|^2\leq  \frac{1}{n}\sum_{i=1}^{n}\mathbb{E}\left[\|\hat{\nabla}f_i(x_{i,k}, y_{i,k}^{(T)};\phi) - \hat{\nabla}f_i(x_{i,k}, y_{i,k}^{(T)})\|^2\right],
    \]
    proves the lemma.
\end{proof}

Now we are ready to give the final proof. Taking summation on both sides of \eqref{ineq: Phi_main1_stoc} and putting Lemma \ref{lem: main_term1_stoc_case2} and \ref{lem: main_term2_stoc_case2} together we know:
\[
    \begin{aligned}
        &\frac{1}{K+1}\sum_{k=0}^{K}\mathbb{E}\left[\|\nabla\Phi(\overline{x_{k}})\|^2\right]\leq \frac{2}{\eta_x(K+1)}(\mathbb{E}\left[\Phi(\overline{x_{0}})\right] - \inf_x\Phi(x)) \\ +&\frac{1}{K+1}\sum_{k=0}^{K}\mathbb{E}\left[\|\mathbb{E} \left[ \overline{\partial \Phi(X_k; \phi)}|\mathcal{F}_k\right] - \nabla\Phi(\overline{x_{k}})\|^2\right] + \frac{L\eta_x}{K+1} \sum_{k=0}^{K}\mathbb{E}\|\left[ \overline{\partial \Phi(X_k; \phi)}\right] - \mathbb{E} \left[ \overline{\partial \Phi(X_k; \phi)}|\mathcal{F}_k\right]\|^2 \\
        \leq& \frac{2}{\eta_x(K+1)}(\mathbb{E}\left[\Phi(\overline{x_{0}})\right] - \inf_x\Phi(x)) + \frac{6C(L^2(1 + \kappa^2))}{T} + 6L_{f,0}^2C\alpha^N + 2L_f^2(1+\kappa^2)\frac{\eta_x^2}{(1-\rho)^2}n(2(1+\kappa)^2 + 2C\alpha^N) \\
        + &L\eta_x(4\sigma_f^2(1+\kappa^2) + (8L_{f,0}^2 + 4\sigma_f^2)\frac{C}{N}) \\
        = &\frac{2}{\eta_x(K+1)}(\mathbb{E}\left[\Phi(\overline{x_{0}})\right] - \inf_x\Phi(x)) + \eta_x^2\frac{4nL_f^2(1+\kappa^2)}{(1-\rho)^2}((1+\kappa)^2+C\alpha^N) \\
        + &L\eta_x(4\sigma_f^2(1+\kappa^2) + (8L_{f,0}^2 + 4\sigma_f^2)\frac{C}{N})  + \tilde{C}_3,
    \end{aligned}
\]
which completes the proof. Here $\tilde{C}_3 = \frac{6C(L^2(1 + \kappa^2))}{T} + 6L_{f,0}^2C\alpha^N = \mathcal{O}(\frac{1}{T}+\alpha^N)$. By setting $\eta_x = \Theta(K^{-\frac{1}{2}}),\ \eta_y = \Theta(K^{-\frac{1}{2}}),\ T=\Theta(K^{\frac{1}{2}}),\ N=\Theta(\log K)$, we have:
\[
    \frac{1}{K+1}\sum_{k=0}^{K}\mathbb{E}\left[\|\nabla\Phi(\overline{x_{k}})\|^2\right] = \mathcal{O}(\frac{1}{\sqrt{K}}).
\]
\end{document}